\documentstyle[epsf]{article}

\newcommand{\ee}{\end{equation}}
\newcommand{\be}{\begin{equation}}

\begin{document}

\centerline {\bf\huge On the Critical Behavior, the Connection Problem  }
\vskip 0.3 cm 
\centerline {\bf\huge  and the Elliptic Representation}
 \vskip 0.3 cm 
\centerline {\bf\huge 
of a  Painlev\'e 6 Equation   (Spring 2001)}
\vskip 0.3 cm
\centerline{\bf   Davide Guzzetti, May 2001}
\vskip 1 cm
\centerline{~~\bf Abstract.} 

\vskip 0.3 cm
 In this paper we find a  class of solutions of the 
sixth Painlev\'e equation appearing in the theory of WDVV equations. This
class 
covers almost all  the {\it monodromy data} associated to the
equation, except one point in the space of the data.   
We describe the critical behavior close to
the critical points in terms of two parameters and we 
find the relation  among the parameters at the
different critical points (connection problem). 
   We also study the critical behavior of 
Painlev\'e  transcendents in the Elliptic representation. 
%


\vskip 1 cm 
AMS CLASSIFICATION: 34M55

\section{Introduction}\label{Introduction}
 
This paper was completed in May 2001, at RIMS, Kyoto University. It is
 published in J. Math. Phys. Anal. Geom. {\bf 4}: 293-377 (2001). I put it on the archive with 9
 years delay, for completeness sake. In the meanwhile, there have been
 several progresses. Nevertheless, this paper contains an important detailed
 analysis which has not been repeated in subsequent papers. 

This work is devoted to the study of the critical behavior of the solutions
 of a Painlev\'e 6 equation given by a particular choice of the
 four parameters  $\alpha$, $\beta$, $\gamma$ $\delta$ of the equation (the
 notations are the standard ones, see \cite{IN}):
$$ 
\alpha={(2\mu-1)^2\over 2},~~~\beta=\gamma=0,~~~\delta={1 \over
 2},~~~~~\mu\in {\bf C}.
$$ 
The equation is 
$$
{d^2y \over dx^2}={1\over 2}\left[ 
{1\over y}+{1\over y-1}+{1\over y-x}
\right]
           \left({dy\over dx}\right)^2
-\left[
{1\over x}+{1\over x-1}+{1\over y-x}
\right]{dy \over dx}$$
\be
+{1\over 2}
{y(y-1)(y-x)\over x^2 (x-1)^2}
\left[
(2\mu-1)^2+{x(x-1)\over (y-x)^2}
\right]
,~~~~~\mu \in {\bf C}.
\label{painleveintroduzione}
\ee
Such an equation will be denoted $PVI_{\mu}$ in the paper. 
The motivation of our work is that (\ref{painleveintroduzione}) 
 is equivalent to
the WDVV equations of associativity in 2-D topological field theory introduced
by Witten \cite{Witten}, Dijkgraaf, Verlinde E., Verelinde H. \cite{DVV}. Such
an equivalence is discussed in \cite{Dub1} and it is a consequence
of the 
theory of Frobenius manifolds.  Frobenius manifolds 
 are the geometrical setting for  the WDVV equations
and were introduced by Dubrovin in \cite{Dub4}. They  are an 
important object in   many branches of
mathematics  like singularity theory and 
 reflection groups \cite{Saito} \cite{SYS}  \cite{Dub6}  \cite{Dub1},
algebraic and 
enumerative geometry \cite{KM} \cite{Manin}. 

\vskip 0.2 cm 

The six classical Painlev\'e  equations were discovered by Painlev\'e \cite{pain} 
 and 
Gambier \cite{gamb}, 
who classified all the second order ordinary differential equations
of the type 
$$
         {d^2 y \over dx^2}= {\cal R}\left(x,y,{dy\over dx}\right)
$$
where ${\cal R}$ is rational in ${dy\over dx}$,  $x$ and
 $y$. The Painlev\'e  equations 
satisfy the {\it Painlev\'e property} of  absence of movable branch points and essential singularities. These singularities will be called {\it critical points}; for $PVI_{\mu}$ 
they are 0,1,$\infty$. The behavior of a solution close to a critical point is called {\it critical 
behavior}.   The 
 general solution of the sixth Painlev\'e equation  
 can be analytically continued to a meromorphic function on the universal
   covering of ${\bf P}^1\backslash \{ 0,1 ,\infty \}$. 
 For generic values of the integration constants and of the parameters in
   the equation, it can not be expressed via elementary or classical
   transcendental functions. For this reason, it  is called a {\it
   Painlev\'e transcendent}.

\vskip 0.2 cm 
The critical behavior  for a class of solutions to  
the Painlev\'e 6 equation was found  by
Jimbo in  \cite{Jimbo}  
for the general Painlev\'e  equation  with generic values of 
 $\alpha$, $\beta$, $\gamma$ $\delta$ (we refer to  \cite{Jimbo} for a precise definition of 
{\it generic}).    
 A transcendent in this class has behavior: 
\be
y(x)= a^{(0)} x^{1-\sigma^{(0)}}(1+O(|x|^{\delta})),~~~~x\to 0,
\label{loc1introduzione}
 \ee
\be
y(x)= 1-a^{(1)}(1-x)^{1-\sigma^{(1)}} (1+O(|1-x|^{\delta})),~~~~x\to 1,
\label{loc2introduzione}
 \ee
\be
y(x)= a^{(\infty)}
 x^{-\sigma^{(\infty)}}(1+O(|x|^{-\delta})),~~~~x\to \infty,
\label{loc3introduzione}
\ee
where $\delta$ is a small positive number, $a^{(i)}$ and
 $\sigma^{(i)}$ are complex numbers such that $a^{(i)}\neq 0$ and 
 \be
0\leq \Re \sigma^{(i)}<1.
\label{RESTOacasa}
\ee 
We remark that  $x$ converges
 to the critical points {\it inside a sector} with vertex
 on the corresponding critical point.

The {\it connection problem}, i.e.  the problem of 
finding  the relation among the three pairs $(\sigma^{(i)},a^{(i)})$,
$i=0,1,\infty$, was solved  by Jimbo in \cite{Jimbo} for the above class of transcendents 
  using the
isomonodromy deformations theory. He considered a  
fuchsian system 
$$
   {dY\over dz}=\left[ {A_0(x)\over z}+{A_x(x) \over z-x}+{A_1(x)\over
z-1}\right] Y$$
such  that the  $2\times 2$ matrices  $A_i(x)$ ($i=0,x,1$ are labels) satisfy Schlesinger equations. 
This ensures that the dependence on $x$ is isomonodromic, 
according to the isomonodromic deformation theory developed in \cite{JMU}. Moreover, for a special choice of the matrices, the Schlesinger  equations are equivalent to the sixth Painlev\'e equation, as 
it is explained in \cite{JM1}. In particular, 
 the local behaviors (\ref{loc1introduzione}), (\ref{loc2introduzione}),
(\ref{loc3introduzione})  were
obtained  using  
a result on the asymptotic
 behavior of a class of  solutions of Schlesinger equations 
 proved by Sato, Miwa, Jimbo in \cite{SMJ}.  The
 connection problem was solved because the parameters
 $\sigma^{(i)}$, $a^{(i)}$ were expressed as functions of 
 the monodromy data of the fuchsian system. For  studies on the 
 asymptotic behavior of the
coefficients of  Fuchsian systems and Schlesinger
equations  see also \cite{Boli}.

\vskip 0.2 cm

Later, Dubrovin and Mazzocco \cite{DM} applied Jimbo's procedure to $PVI_{\mu}$, with the 
restriction that $2\mu \not \in {\bf Z}$. We remark that this case was not studied by Jimbo, being 
  a non-generic case.   Dubrovin and Mazzocco obtained a class of transcendents with behaviors (\ref{loc1introduzione}), (\ref{loc2introduzione}), (\ref{loc3introduzione})  (again, 
$x$ converges to a critical point inside a sector) 
and restriction (\ref{RESTOacasa}). They also solved the connection problem. 

In the case of $PVI_{\mu}$, the  monodromy data of the Fuchsian system, to be introduced later,  
turn out to be expressed in terms
of  a  triple of complex numbers  $(x_0,x_1,x_{\infty})$. 
The two integration constants in $y(x)$ and the parameter $\mu$ are contained
  in the  triple. The following relation holds: 
\be
   x_0^2+x_1^2 +x_{\infty}^2-x_0 x_1 x_{\infty}= 4 \sin^2(\pi \mu).
\label{10iiNUOVA}
\ee
  There
 exists a one-to-one correspondence between  triples (define up to the change of two signs) 
and  branches of 
 the Painlev\'e transcendents. $^1$
%
%
 In other words, any branch $y(x)$ is parameterized by
 a triple:  
 $$
y(x)=y(x;x_0,x_1,x_{\infty}).
$$ 
  As it is proved in \cite{DM}, the
 transcendents (\ref{loc1introduzione}),  (\ref{loc2introduzione}),
 (\ref{loc3introduzione}) are parameterized by a triple according to the
 formulae 
$$ 
   x_i^2= 4 \sin^2\left({\pi\over 2}
   \sigma^{(i)}\right),~~~i=0,1,\infty,~~~~~~~~0\leq \Re \sigma^{(i)}<1.
$$
A more complicated expression gives $a^{(i)}=a^{(i)}(x_0,x_1,x_{\infty})$ in \cite{DM}.  We recall that a branch is defined by the choice of branch cuts, 
like $|\arg(x)|<\pi$, $|\arg(1-x)|<\pi$.  
The analytic 
continuation of a branch when $x$ crosses the cuts 
is obtained by an action of the braid group on the triple. This is explained in \cite{DM} and in section  \ref{Analytic Continuation of a
Branch}.

\vskip 0.3 cm
 As we mentioned above, it is very important to concentrate on $PVI_{\mu}$ due to its equivalence 
to WDVV equations in 2-D topological field theory, and due to its central role in the construction of three-dimensional Frobenius 
manifolds. It is known \cite{Dub1}  that the structure of a local chart of
 a Frobenius manifold can in 
principle be constructed from a set of monodromy data.  To any  manifold a $PVI_{\mu}$ equation is associated and the monodromy data of the local chart are 
contained in $\mu$ and in  the triple $(x_0,x_1, x_{\infty}) $ of a   
 Painlev\'e transcendent. The mentioned action of the braid group, which gives the analytic continuation of the transcendent, allows to pass from one local chart to another.   
  
 The local structure of a Frobenius manifold is explicitly constructed in 
\cite{guz1} starting from the 
Painlev\'e transcendents. In \cite{guz1} it is shown that in order to 
obtain a local chart from its monodromy data we need  
 to know the critical behavior of the corresponding transcendent in terms of the 
 triple $(x_0,x_1,x_{\infty})$ (note that this is equivalent to solving the connection problem). 
 
  Recently, 
 Frobenius manifolds   have become important in enumerative geometry and quantum cohomology 
\cite{KM}, \cite{Manin}. As it is shown in \cite{guz1}, it is possible to compute Gromov-Witten invariants for the quantum cohomology of the two-dimensional projective space 
starting from a special $PVI_{\mu}$, with $\mu=-1$. In this case the triple is 
 $(x_0,x_1,x_{\infty})=(3,3,3)$, as it is proved in \cite{Dub2} and \cite{guz}. 
Due to the  restriction $0\leq \Re \sigma^{(i)}<1$, the formulae for the critical behavior and the 
connection problem obtained by 
Dubrovin-Mazzocco do not apply  
if at least one $x_i$ ($i=0,1,\infty$) is real and $|x_i|\geq
 2$. Thus, they do not apply in  the case of quantum cohomology, because 
 $x_i=3$ and $\Re \sigma^{(i)}=1$.

Therefore, the motivation of our paper becomes clear: in the attempt  to extend the results of \cite{DM} to the case of quantum cohomology, 
we actually extended them to almost all monodromy data, namely we found  the critical behavior and we solved the
 connection problem for  all the triples satisfying
$$
    x_i\neq \pm 2 ~~\Longrightarrow ~~ \sigma^{(i)} \neq 1,~~~~i=0,1,\infty.
$$
 In order to do this, we  extended  Jimbo and Dubrovin-Mazzocco's method and 
we analyzed the elliptic representation of the Painlev\'e 6 equation.

\subsection{ Our results}

 We observe that the branch $y(x;x_0,x_1,x_{\infty})$ has analytic 
continuation on the universal covering of ${\bf P}^1\backslash \{0,1,\infty\}$. We still denote 
this continuation by  $y(x;x_0,x_1,x_{\infty})$, where $x$ is now a point in the universal covering. 
Therefore:

\vskip 0.2 cm
{\it There is a one to one correspondence 
 between triples of monodromy data $(x_0,x_1,x_{\infty})$ (defined up to the change of two signs)   and Painlev\'e transcendents, 
namely $y(x)=y(x;x_0,x_1,x_{\infty})$, $x\in \widetilde{{\bf P}^1\backslash} \{0,1,\infty\}$ }

\vskip 0.2 cm
 We  mentioned that if we fix a branch, namely if  we 
choose branch cuts like $|\arg x|<\pi$,  $|\arg(1-x)| <\pi$, then the branch of $y(x;x_0,x_1,x_{\infty})$ 
has analytic continuation $y(x;x_0^{\prime},x_1^{\prime},x_{\infty}^{\prime})$ in the cut plane, where 
$( x_0^{\prime},x_1^{\prime},x_{\infty}^{\prime})$ is obtained from $(x_0,x_1,x_{\infty})$ by an action of the braid group (see section \ref{Analytic Continuation of a Branch} for details).

\vskip 0.3 cm
 We obtained the following results
\vskip 0.2 cm
1){\it A transcendent $y(x;x_0,x_1,x_{\infty})$ such that $|x_i| \neq 2$ has behaviors    
(\ref{loc1introduzione}), (\ref{loc2introduzione}),  (\ref{loc3introduzione}) in suitable domains, 
to be defined below, 
 contained in $\widetilde{{\bf C}\backslash\{0\}}$, $\widetilde{{\bf C}\backslash\{1\}}$,
 $\widetilde{{\bf P}^1\backslash\{\infty\}}$ respectively. The exponent are restricted by the condition 
$\sigma^{(i)}\not \in (-\infty,0)\cup[1,+\infty)$, which extends (\ref{RESTOacasa}). 
}

\vskip 0.2 cm
2) {\it The parameters $\sigma^{(i)}$, $a^{(i)}$ are computed  as functions of $(x_0,x_1,x_{\infty})$, and vice versa, by explicit formulae which extend those of \cite{DM} }

\vskip 0.2 cm

3) {\it If we enlarge the domains where   (\ref{loc1introduzione}), (\ref{loc2introduzione}),  (\ref{loc3introduzione}) hold, the behavior of $y(x;x_0,x_1,x_{\infty})$ becomes 
oscillatory. The movable poles of the 
transcendent lie outside the enlarged domains. In proving this, we investigated the 
elliptic representation of the transcendent, providing a general result 
stated in Theorem 3  below. }

\vskip 0.3 cm 

 We state result 1) in more detail.    
 Let $\sigma^{(0)} $ be a complex number 
such that $\sigma^{(0)}\not \in (-\infty,0)\cup 
[1,+\infty)$. 
  We introduce additional parameters 
$\theta_1$, $\theta_2 \in {\bf R}$, $0<\tilde{\sigma} 
<1$ to  define a domain  
$$
   D(\epsilon;\sigma^{(0)};\theta_1,\theta_2,\tilde{\sigma}):= 
\{{ x}\in \widetilde{{\bf C}_0}
   \hbox{ s.t. } |x|<\epsilon ,~~ e^{-\theta_1 \Im \sigma^{(0)}}
   |x|^{\tilde{\sigma}} \leq |{x}^{\sigma^{(0)}}| 
\leq e^{-\theta_2 \Im
   \sigma^{(0)}} ,~~~0<\tilde{\sigma}<1
   \},  
$$  
 which can can be rewritten as 
$$|x|<\epsilon,~~~~~ 
        \Re \sigma^{(0)} \log|x|+\theta_2 \Im \sigma^{(0)} \leq \Im \sigma^{(0)} \arg(x)
        \leq (\Re \sigma^{(0)} - \tilde{\sigma}) \log|x| + \theta_1 \Im
        \sigma^{(0)} 
$$
 For real $\sigma^{(0)}$ the domain is more simply defined as  
\be
D(\sigma^{(0)};\epsilon):= \{ { x}\in \widetilde{{\bf C}_0}
   \hbox{ s.t. } |x|<\epsilon \},~~~~~\hbox{ for  } 0\leq \sigma^{(0)}<1
\label{DOMINORIMS}
\ee

 For semplicity, we study the critical behavior of the transcendent for  $x\to 0$ along the family of path defined below.  
Such paths start at some point $x_0$ belonging to the domain. 
 If $\Im
\sigma=0$ any regular path will be allowed. If $\Im
\sigma\neq 0$, we considered the family 
\be
|x|\leq|x_0|<\epsilon,~~~~\arg x= \arg x_0 +{\Re \sigma^{(0)} -\Sigma \over 
                      \Im \sigma^{(0)} } \ln{|x| \over |x_0|}, ~~~~~0\leq \Sigma \leq 
\tilde{\sigma}
\label{spirale}
\ee
The condition $0\leq \Sigma \leq 
\tilde{\sigma}$ ensures that the paths remain in the domain as $x\to 0$. In general, these paths are spirals.

\vskip 0.2 cm 
\noindent
{\bf Theorem 1: } {\it Let $\mu \neq 0$.  For any $\sigma^{(0)} \not \in
(-\infty,0)\cup [1,+\infty)$, for
any $a^{(0)} \in {\bf 
C}$, $a^{(0)}\neq 0$, for any $\theta_1,\theta_2  \in {\bf R}$  and for any
$0<\tilde{\sigma }<1$, 
there exists a sufficiently small positive $\epsilon$ 
and a small positive number $\delta$ 
such that  the  equation (\ref{painleveintroduzione}) has a solution     
\be
    y(x;\sigma^{(0)},a^{(0)})=a(x)~ x^{1-\sigma^{(0)}} \left(
1+O(|x|^{\delta})
\right),~~~~0<\delta<1,
\label{fujisan0}
\ee
as $x\to 0$ along (\ref{spirale}) in the domain $
D(\epsilon;\sigma^{(0)};\theta_1,\theta_2,\tilde{\sigma}) 
  $ defined for non-real $\sigma^{(0)}$,  or along any regular path in  $
D(\epsilon;\sigma^{(0)})$ defined for real $0\leq\sigma^{(0)}<1$. The amplitude $a(x)$ is 
$$
  a(x):= a^{(0)}, ~~~~\hbox{ for $0<\Sigma\leq \tilde{\sigma}$, or for real $\sigma^{(0)}$  }
$$
\be
  a({x}):= a^{(0)}~\left(1 +{1\over 2a^{(0)}}~ |x_0^{\sigma^{(0)}}| e^{i \alpha({x})} 
             + {1\over 16 [a^{(0)}]^2}~ |x_0^{\sigma^{(0)}}|^2 e^{2 i \alpha({x})} 
 \right) = O(1),~~~~\hbox{ for $\Sigma=0$}
\label{cecilia}
\ee
where we have used the notaion $\alpha(x)$ to denote the real phase of $x^{\sigma^{(0)}}= |x^{\sigma^{(0)}}|~e^{i\alpha(x)}\equiv 
|x_0^{\sigma^{(0)}}|~e^{i\alpha(x)}
$, when $\Sigma=0$.
}

\vskip 0.3 cm

Note that in the case (\ref{cecilia})  we can rewrite 
\be
  y({x};\sigma^{(0)},a^{(0)})= \sin^2\left({i\sigma^{(0)} \over 2} \ln x -{i \over 2} \ln (4a^{(0)}) 
-{\pi \over 2} \right)~x~ (1+O(|x|^{\delta}))
\label{annaalbertostellacecilia}
\ee
For brevity, we will sometime denote the domain by $
D(\sigma^{(0)})$. The condition $\mu \neq 0$ is not restrictive because $PVI_{\mu=0}$ 
coincides with $PVI_{\mu=1}$. 
  
From Theorem 1 and the symmetries of (\ref{painleveintroduzione}), 
we  
prove the existence of solutions with the following  local behaviors 
$$y(x,\sigma^{(1)},a^{(1)})= 1-a^{(1)} (1-x)^{1-\sigma^{(1)}}
(1+O(|1-x|^{\delta }))~~~~x\to 1 
$$
$$
   a^{(1)}\neq 0,~~~~\sigma^{(1)}\not \in (-\infty,0)\cup [1,+\infty)
$$
and $$
    y(x; \sigma^{(\infty)},
 a^{(\infty)}) = a^{(\infty)} {x}^{ \sigma^{(\infty)} }\left(
 1+O({1\over |x|^{\delta}}) \right)~~~~~x \to \infty
$$
$$
  a^{(\infty)}\neq 0,~~~~\sigma^{(\infty)}\not \in (-\infty,0)\cup [1,+\infty)
$$
in  domains $D(\sigma^{(1)})$, $D(\sigma^{(\infty)})$ given by (\ref{ANCORAACASAuno}), 
 (\ref{ANCORAACASAinf}) respectively.

The critical behaviors above coincide 
 with (\ref{loc1introduzione}), (\ref{loc2introduzione}) and (\ref{loc3introduzione})  for $0\leq \Re \sigma^{(i)}<1$, $i=0,1,\infty$. 
But  our result is more general because it  extends 
 the range of $\sigma^{(i)}$ to  $\Re\sigma^{(i)}<0$ and
$\Re \sigma^{(i)}\geq 1$. For this larger range,  $x$ may tend to $x=i$ 
($i=0,1,\infty$)  
 along a
spiral, according to the shape of  $
D(\sigma^{(i)})$.  For more comments see sections \ref{Local Behaviour -- Theorem 1} and 
\ref{Singular Points x=1, x=infty  (Connection Problem)}.

\vskip 0.3 cm 
 Result 2) is stated in the theorem below  -- where 
we  write $\sigma,a$ instead of $\sigma^{(0)}, 
a^{(0)}$ -- and in its comment.  

\vskip 0.3 cm
\noindent
{\bf Theorem 2:} {\it Let $\mu$ be any non zero complex number. 

The transcendent $y(x;\sigma,a)$ of Theorem 1, 
defined for  $\sigma \not \in (-\infty,0)\cup [1,+\infty)$  
and $a\neq 0$,  is the representation of a transcendent $y(x;x_0,x_1,x_{\infty})$ in $D(\sigma)$. The 
triple   
($x_0,x_1,x_{\infty}$) is uniquely determined (up to the change of two signs) 
by the following formulae: 
 \vskip 0.15 cm
\noindent
 i) $\sigma \neq 0, \pm 2\mu+2m$ for any $m\in {\bf Z}$. 
$$
   \left\{ \matrix{x_0=2 \sin({\pi \over 2}\sigma) \cr\cr
                  x_1=i\left({1\over
f(\sigma,\mu)G(\sigma,\mu)}~\sqrt{a}-G(\sigma,\mu) ~{1\over \sqrt{a}}\right)
                      \cr\cr
                  x_{\infty}= {1\over f(\sigma,\mu)G(\sigma,\mu)
e^{-{i\pi\sigma\over 2}} }~\sqrt{a}+ G(\sigma,\mu) e^{-i{\pi\sigma\over
2}}~{1\over \sqrt{a}}\cr
                         }\right. 
$$
where 
$$f(\sigma,\mu)={2 \cos^2({\pi \over 2} \sigma) \over \cos(\pi \sigma)- 
                  \cos(2\pi\mu)}, ~~~
~~~ G(\sigma,\mu)= {1\over 2}{4^{\sigma}\Gamma({\sigma+1\over 2})^2\over
                  \Gamma(1-\mu+ {\sigma\over
                  2})\Gamma(\mu+{\sigma\over 2})}$$
Any sign  of $\sqrt{a}$  is good (changing the sign of $\sqrt{a}$ 
                is equivalent
                  to changing the sign of both $x_1$, $x_{\infty}$).
\vskip 0.15 cm
\noindent
 ii) $\sigma =0$
$$\left\{ 
              \matrix{
                         x_0=0  \cr\cr
                        x_1^2=2\sin(\pi\mu)~\sqrt{1-a} \cr\cr
                        x_{\infty}^2=2\sin(\pi\mu)~\sqrt{a}
                        \cr
                       }
 \right.  
$$
We can take any sign of the square roots 
 
\vskip 0.15 cm
\noindent
 iii) $\sigma= \pm 2\mu +2 m$.  

\vskip 0.15 cm
  iii1) $\sigma=2\mu+2m$, $m=0,1,2,...$
 $$ 
         \left\{ 
   \matrix{ 
              x_0=2\sin(\pi\mu)\cr \cr
              x_1= -{i\over 2}{16^{\mu+m} \Gamma(\mu+m+{1\over 2})^2 \over
   \Gamma(m+1) \Gamma(2\mu+m)}~{1\over\sqrt{a}} 
                                           \cr \cr
               x_{\infty}=i~x_1~e^{-i\pi\mu}\cr
          }
                  \right.
$$

 iii2) $\sigma=2\mu+2m$, $m=-1,-2,-3,...$
$$
  \left\{ 
\matrix{
            x_0=2\sin(\pi\mu)\cr \cr
            x_1=2i{\pi^2\over \cos^2(\pi\mu)}
                {1\over 16^{\mu+m} \Gamma(\mu+m+{1\over 2})^2
\Gamma(-2\mu-m+1) \Gamma(-m)}~\sqrt{a}
\cr\cr
x_{\infty}=-ix_1 e^{i\pi\mu}
}
\right.
$$

iii3) $\sigma=-2\mu +2m$, $m=1,2,3,...$
$$
   \left\{ 
\matrix{
           x_0=-2\sin(\pi\mu) \cr\cr
         x_1= -{i\over 2}{16^{-\mu+m} \Gamma(-\mu+m+{1\over 2})^2\over
\Gamma(m-2\mu+1) \Gamma(m)}~{1\over \sqrt{a}} \cr\cr
        x_{\infty}= ix_1e^{i\pi\mu}\cr
}
\right.
$$

iii4) $\sigma= -2\mu+2m$, $m=0,-1,-2,-3,...$
$$
  \left\{
\matrix{
            x_0=-2\sin(\pi\mu) \cr\cr
            x_1= 2i {\pi^2\over \cos^2(\pi\mu)}{1\over 16^{-\mu+m}
\Gamma(-\mu+m+{1\over 2})^2\Gamma(2\mu-m) \Gamma(1-m)}  ~\sqrt{a}
\cr\cr
x_{\infty}=-ix_1 e^{-i\pi \mu}\cr
}
 \right.
$$
In all the above formulae the relation
$x_0^2+x_1^2+x_{\infty}^2-x_0x_1x_{\infty}=4 \sin^2(\pi \mu)$  is
automatically satisfied. Note that $\sigma \neq 1$ implies $x_0\neq \pm
2$. Changes of two signs in the triple of  the formulae above are allowed.


\vskip 0.2 cm 

 Conversely, a transcendent $y(x;x_0,~x_1,~x_{\infty})$,  
 such that $x_0^2+x_1^2+x_{\infty}^2-x_0x_1x_{\infty}=4 \sin^2(\pi
 \mu)$, $x_i\neq \pm 2$,  has representation $y(x;\sigma,a)$ in $D(\sigma)$ of Theorem 1 
 with  parameters   $\sigma$  and $a$
 obtained as follows:  
 
\vskip 0.15 cm 
\noindent
 I) Generic case 
$$
\cos(\pi \sigma)= 1-{x_0^2\over 2}
$$
$$                  %
  a={iG(\sigma,\mu)^2\over 2 \sin(\pi \sigma)} 
      \Bigl[
 2(1+e^{-i\pi\sigma})-f(x_0,x_1,x_{\infty})(x_{\infty}^2+e^{-i\pi\sigma} x_1^2)
\Bigr]  ~ f(x_0,x_1,x_{\infty})
$$
where  
$$
   f(x_0,x_1,x_{\infty}):=f(\sigma(x_0),\mu)={4-x_0^2\over 2-x_0^2-
2\cos(2\pi\mu)}={4-x_0^2\over
    x_1^2+x_{\infty}^2-x_0x_1x_{\infty}}. $$
$\sigma$ 
 is determined up to the ambiguity $\sigma \mapsto \pm 
\sigma +2n $, $n\in {\bf Z}$ [see Remark below]. If $\sigma$ is real we can only choose the 
solution satisfying $0\leq \sigma<1$. 
Any solution $\sigma $ of the first equation must satisfy the additional 
restriction
$\sigma \neq  \pm 2\mu+2m$ for any $m\in {\bf Z}$, otherwise we
encounter the singularities in $G(\sigma,\mu)$ and in $f(\sigma,\mu)$.

\vskip 0.15 cm 
\noindent
 II) $x_0=0$. 
           $$ \sigma=0,$$
$$
    a= {x_{\infty}^2\over x_1^2+x_{\infty}^2}.
$$
provided that $x_1\neq 0$ and $x_{\infty}\neq 0$, namely  $\mu
\not\in {\bf Z}$.  
\vskip 0.15 cm 
\noindent
III) $x_0^2=4 \sin^2(\pi\mu)$.  Then (\ref{10iiNUOVA}) implies
 $x_{\infty}^2=-x_1^2 ~\exp(\pm 2\pi i \mu)$ . 
 Four cases which yield the values of $\sigma$ non included in I)
and II) must be considered 

\vskip 0.15 cm
   III1) If  $x_{\infty}^2=-x_1^2 e^{- 2
\pi i \mu}$ then 
$$\sigma=  2\mu + 2m,~~~~m=0,1,2,...$$ 
 $$ a=-{1\over 4 x_1^2} 
               { 16^{2\mu+2m}
 \Gamma(\mu+m+{1\over 2})^4\over \Gamma(m+1)^2 \Gamma(2\mu+m)^2  } 
$$

\vskip 0.15 cm
  III2) If $x_{\infty}^2=-x_1^2 e^{2\pi i \mu}$ then 
$$\sigma=2\mu+2m,~~~~m=-1,-2,-3,...$$
    $$
            a=-{\cos^4(\pi\mu)\over 4 \pi^4}
                 16^{2\mu+2m} \Gamma(\mu+m+{1\over 2})^4
\Gamma(-2\mu-m+1)^2 \Gamma(-m)^2~ x_1^2  
$$

\vskip 0.15 cm
  III3)If  $x_{\infty}^2=-x_1^2e^{2\pi i\mu}$ then 
$$\sigma=-2\mu+2m,~~~~m=1,2,3,...$$
$$a=-{1\over 4 x_1^2}
         {
  16^{-2\mu+2m} \Gamma(-\mu+m+{1\over 2})^4\over \Gamma(m-2\mu+1)^2
\Gamma(m)^2  }
$$

\vskip 0.15 cm
III4) If  $x_{\infty}^2=-x_1^2 e^{-2\pi i \mu}$ then 
$$\sigma=-2\mu+2m,~~~~m=0,-1,-2,-3,...$$
$$
  a= -{\cos^4(\pi\mu)\over 4 \pi^4} 16^{-2\mu+2m}
\Gamma(-\mu+m+{1\over 2})^4\Gamma(2\mu-m)^2 \Gamma(1-m)^2 ~ x_1^2
$$

  }

\vskip 0.3 cm

 Let us restore the notation $\sigma^{(0)},a^{(0)}$. At  $x=1,\infty$ the 
exponents $\sigma^{(i)}$, $i=1,\infty$ are given by $\cos (\pi \sigma^{(i)})=1-{x_i^2\over 2}$ and 
  the coefficients  $a^{(1)}$, $a^{(\infty)}$ 
 are obtained from the formula of $a=a^{(0)}$ of Theorem 2,  provided that we do the
substitutions $(x_0,x_1,x_{\infty})\mapsto ( x_1,x_0,x_0
   x_1-x_{\infty})$,  $\sigma^{(0)} \mapsto
\sigma^{(1)}$ and  $(x_0,x_1,x_{\infty})\mapsto (x_{\infty},
-x_1,x_0-x_1x_{\infty})$, $\sigma^{(0)} \mapsto
\sigma^{(\infty)}$ respectively. 

 This also solves the  {\it connection problem} 
 for the transcendents $y(x;\sigma^{(i)}, 
a^{(i)})$, because we are able to compute  $
(\sigma^{(i)}, 
a^{(i)})$ for $i=0,1,\infty$ in terms of a fixed triple $(x_0,x_1,x_{\infty})$.

\vskip 0.2 cm
\noindent
{\it Remark : } 
Let $(x_0,x_1,x_{\infty})$ be given and let us compute $\sigma$ and $a$ 
by the formulae of Theorem 2.  The  equation
\be
 \cos(\pi \sigma) = 1-{x_0^2\over2}
\label{LERBAVOGLIO}
\ee
 does not 
determine $\sigma$ uniquely. We can choose $\sigma$ such that 
$$
0\leq \Re \sigma \leq 1.
$$
This convention will be  assumed in the paper.   
Therefore, all the solutions of (\ref{LERBAVOGLIO})
 are 
$$
\pm \sigma +2 n, ~~~~n\in {\bf Z}.
$$ 
If $\sigma$ is real, we can only 
choose $ 
    0\leq \sigma<1  
$. 
With this convention, there is a one to one correspondence between 
$(\sigma,a)$ and 
(a class of equivalence, defined by the change of two signs, of) an admissible triple $(x_0,x_1,x_{\infty})$. 

We observe that 
$\sigma =\sigma(x_0)$ and $a= a(\sigma;x_0,
x_1,x_{\infty})$; namely,
 the transformation $\sigma  \mapsto \pm \sigma +2 n$ affects $a$.  
The  transcendent $y(x;x_0,x_1,x_{\infty})$ has representation $ 
y\bigl(x;~\sigma(x_0),~a(\sigma;x_0,
x_1,x_{\infty})~\bigr)$ in $D(\sigma)$. If we choose another  solution $\pm \sigma +2n$ 
we again have $y(x;x_0,x_1,x_{\infty})=
 y\bigl(x;~\pm \sigma(x_0)+2n,~a(\pm\sigma(x_0)+2n;x_0,
x_1,x_{\infty})~\bigr)$ in the new domain $ D(\pm \sigma + 2n)$.  Hence -- and this is very 
important! -- the transcendent  $y(x;x_0,x_1,x_{\infty})$ 
has different representations and different critical behaviors in different domains. Outside the union of these domains  we are not able to describe 
 the transcendents and we believe that the movable poles lie there (we show this in  one example in the paper).

\vskip 0.3 cm

 According to the above remark, we  restrict to the case $0\leq \Re
\sigma^{(i)}\leq 1$, $\sigma^{(i)}\neq 1$. So the critical behaviors
 of $y(x;\sigma^{(i)}, a^{(i)})$ coincide with (\ref{loc1introduzione}),
(\ref{loc2introduzione}), (\ref{loc3introduzione}) when  $0\leq\Re
\sigma^{(i)} <1$. But for $\Re \sigma^{(i)}=1$ 
 the critical behaviors (\ref{loc1introduzione}),
(\ref{loc2introduzione}), (\ref{loc3introduzione})   hold true {\it only} if $x$
converges to a critical point { \it along spirals}. 

\vskip 0.3 cm 
 We finally describe the third result. 
In the case $\Re \sigma^{(i)}=1$, we  obtained the critical behaviors 
 {\it along radial paths} using the elliptic 
representation of Painlev\'e transcendents. We only 
consider now  the critical point $x=0$, because the symmetries of
(\ref{painleveintroduzione}), to be discussed in section 
\ref{Singular Points x=1, x=infty  (Connection Problem)}, 
 yield the behavior close to the other critical
points.  

The elliptic representation was introduced by R.Fuchs in \cite{fuchs}:
$$ 
  y(x)=\wp\left({u(x)\over 2};\omega_1(x),\omega_2(x)\right)+{1+x\over 3}
$$
Here  $u(x)$ solves a non-linear second order differential equation to be 
studied  later  and
  $\omega_1(x)$, $\omega_2(x)$ are two elliptic integrals, expanded for
  $|x|<1$ in terms of hyper-geometric functions: 
$$
\omega_1(x)= {\pi \over 2} 
\sum_{n=0}^{\infty}{ \left[\left({1 \over 2}\right)_n\right]^2
  \over (n!)^2 } x^n
$$
$$
\omega_2(x)= -{i\over 2}\left\{ 
\sum_{n=0}^{\infty}{ \left[\left({1 \over 2}\right)_n\right]^2
  \over (n!)^2 } x^n ~\ln(x) +   
\sum_{n=0}^{\infty}{ \left[\left({1 \over 2}\right)_n\right]^2
  \over (n!)^2 } 2\left[ \psi(n+{1\over 2}) - \psi(n+1)\right]
x^n \right\}
$$
where $\psi(z):= {d\over dz} \ln \Gamma(z)$.   We introduce a new domain, depending on 
two complex numbers $\nu_1,\nu_2$ and on the small real number $r$: 
 $$
  {\cal D}(r;\nu_1,\nu_2):= \left\{ x\in \tilde{\bf C}_0~ \hbox{ such that }
  |x|<r, \left|{e^{-i\pi \nu_1}\over 16^{2-\nu_2}} x^{2-\nu_2} \right|<r,
\left| 
{e^{i\pi \nu_1} \over 16^{\nu_2}} x^{\nu_2}\right|<r \right\}
$$
The domain can be also written as follows:
$$
|x|<r,~~~~   \Re \nu_2 \ln|x|+ C_1-\ln r < \Im \nu_2
   \arg x < (\Re \nu_2 -2)\ln|x| +C_2 + \ln r,
$$
$$
C_1:= -\bigl[4 \ln 2 ~\Re \nu_2 + \pi~ \Im \nu_1\bigr],~~~~C_2:=C_1+8\ln 2,
$$
if $\Im \nu_2\neq 0$. If $\Im \nu_2 =  0$, the domain is simply $
|x|<r$.

\vskip 0.2 cm
\noindent
{\bf Theorem 3:} 
 {\it For any complex $\nu_1$, $\nu_2$ such that 
$$ 
  \nu_2\not \in (-\infty,0]\cup [2,+\infty)
$$
there exists a sufficiently small $r<1$ such that $PVI_{\mu}$ has a solution of the form
$$
  y(x)= \wp(\nu_1 \omega_1(x)+\nu_2 \omega_2(x)
  +v(x);\omega_1(x),\omega_2(x))+{1+x\over 3}
$$
in the domain $ {\cal D}(r;\nu_1,\nu_2)$ defined above. 
The function $v(x)$ is holomorphic  in $ {\cal D}(r;\nu_1,\nu_2)$ and has
convergent expansion 
\be
  v(x)= \sum_{n\geq 1} a_n x^n +\sum_{n\geq 0,~m\geq 1} b_{nm} x^n
  \left({e^{-i\pi \nu_1}\over 16^{2-\nu_2}} x^{2-\nu_2}\right)^m +\sum_{n\geq
  0,~m\geq 1}c_{nm} x^n \left( 
{e^{i\pi \nu_1} \over 16^{\nu_2}} x^{\nu_2}\right)^m    
\label{vdix}
\ee
where $a_n$, $b_{nm}$, $c_{nm}$ are certain rational functions of  $\nu_2$. 
Moreover, there exists a constant $M(\nu_2)$ depending on $\nu_2$ such that 
 $v(x)\leq M(\nu_2) \left(|x|+\left|{e^{-i\pi \nu_1}\over 16^{2-\nu_2}} x^{2-\nu_2} \right|+\left| 
{e^{i\pi \nu_1} \over 16^{\nu_2}} x^{\nu_2}\right| \right)$ in  ${\cal
D}(r;\nu_1,\nu_2)$ .
}

\vskip 0.2 cm

 Theorem 3 allows to compute the critical behavior. We  consider a family of paths along which $x$ may tend to zero, contained in the domain of the theorem. If $0<\nu_2<2$, any regular path is allowed. If $\nu_2$ is any non-real number, we consider the following family, starting at $x_0\in  {\cal D}(r;\nu_1,\nu_2)$:
\be|x|\leq |x_0|<r,~~~~
 \arg(x)= \arg(x_0)+{\Re \nu_2-{\cal V}\over \Im \nu_2} \ln{|x|\over |x_0|},~~~~~0\leq {\cal V}
\leq 2.
\label{takaragaike}
\ee
The restriction $0\leq {\cal V}\leq 2$ ensures that the paths remain in the domain as $x\to 0$. 
\vskip 0.2 cm 
\noindent
{\bf Corollary:} {\it Consider a transcendent $
  y(x)= \wp(\nu_1 \omega_1(x)+\nu_2 \omega_2(x)
  +v(x);\omega_1(x),\omega_2(x)) +{1+x\over 3}
$ of Theorem 3. Its critical behavior for $x\to 0$ in 
 $ {\cal D}(r;\nu_1,\nu_2)$ along (\ref{takaragaike}) if  $\Im \nu_2\neq 0$ and $0<{\cal V}<2$, or along any regular path if $0<\nu_2<2$ is:
\be 
 y(x)= \left[{1\over 2} x - {1\over 4} \left[{e^{i\pi \nu_1} \over
  16^{\nu_2-1} }  \right] x^{\nu_2} - {1\over 4} \left[{e^{i\pi \nu_1} \over
  16^{\nu_2-1} }  \right]^{-1} x^{2-\nu_2}
  \right](1+O(x^{\delta})),
\label{takaidesu}
\ee
for some  $0<\delta<1$. 
If $\Im \nu_2\neq 0$ and ${\cal V}=0$ the behavior along (\ref{takaragaike}) is:
$$ 
y(x)  = {1\over \sin^2\left( -i{\nu_2\over 2} \ln x +\left[i {\nu_2\over 2} \ln 16
   +{\pi 
   \nu_1\over 2}\right] + \sum_{m=1}^{\infty} c_{0m}(\nu_2) \left[{e^{i\pi \nu_1}\over 16^{\nu_2}}x^{\nu_2}\right]^m 
\right)}~\bigl(1+O(x)\bigr). 
$$
If $\Im \nu_2\neq 0$ and ${\cal V}=2$ the behavior along (\ref{takaragaike}) is:
$$
y(x)= {1\over \sin^2\left( i{2-\nu_2\over 2} \ln x +\left[i {\nu_2-2\over 2} 
\ln 16
   +{\pi 
   \nu_1\over 2}\right] + \sum_{m=1}^{\infty} b_{0m}(\nu_2) \left[{
e^{-i\pi \nu_1}\over 16^{2-\nu_2}}x^{2-\nu_2}\right]^m 
\right)}~\bigl(1
+
O(x)\bigr).
$$
} 
\vskip 0.2 cm 

 Note that (\ref{takaidesu}) is 
\be 
  y(x)=  - {1\over 4} \left[{e^{i\pi \nu_1} \over
  16^{\nu_2-1} }  \right] x^{\nu_2} ~(1+O(x^{\delta})), ~~~~~\hbox{ if } 0<{\cal V}<1, ~\hbox{ or } 
0<\nu_2<1.
\label{fujisan1}
\ee
\be
y(x)= - {1\over 4} \left[{e^{i\pi \nu_1} \over
  16^{\nu_2-1} }  \right]^{-1} x^{2-\nu_2}~(1+O(x^{\delta})), ~~~~~\hbox{ if }1<{\cal V}<2, ~\hbox{ or } 
1<\nu_2<2.
\label{fujisan2}
\ee
and 
\be
y(x) = \sin^2 \left(i{1-\nu_2\over 2} \ln{|x|\over 16} + {\pi \nu_1\over 2} \right) ~x~(1+O(x)),~~~~~\hbox{ if }{\cal V}=1, ~\hbox{ or } 
\nu_2=1.
\label{fujisan3}
\ee
 The elliptic 
representation has been studied from the point of view of algebraic geometry 
in \cite{Manin2}, but to our knowledge  Theorem 3 and its Corollary are 
 the first  general 
result   on its  critical behavior. 
We however note that for the very special value 
 $\mu={1\over 2}$ the function $v(x)$ vanishes; the transcendents are called  {\it Picard solutions} in \cite{M}, because they were known to Picard \cite{Picard}. Their 
 critical behavior is
studied in \cite{M} and agrees with the Corollary.  

\vskip 0.2 cm 
 Comparing (\ref{fujisan0}) with (\ref{fujisan1}) we   prove in section 
 \ref{Elliptic Representation} 
 that the  transcendent of 
Theorem 3 coincides with $y(x;\sigma^{(0)},a^{(0)})$ of Theorem
1 on the domain  $D(\epsilon,\sigma^{(0)})
\cap  {\cal D}(r;\nu_1,\nu_2)$ with 
 the identification $\sigma^{(0)}=1-\nu_2$ and $a^{(0)}= -{1\over 4} \left[
{e^{i\pi\nu_1} \over 16^{\nu_2-1}}\right]$ (note also that  
(\ref{annaalbertostellacecilia}) is (\ref{fujisan3})).     The identification of $a^{(0)}$
and $\sigma^{(0)}$  makes it possible 
 to connect $\nu_1$ and $\nu_2$ to the monodromy data 
$(x_0,x_1,x_{\infty})$ according to Theorem
2 and to solve the connection problem for the elliptic representation.

\vskip 0.2 cm 
The Corollary provides the behavior of the transcendents when 
 $\Re \sigma^{(0)}=1$ ($\sigma^{(0)}\neq 1$) and $x\to 0$ along a radial path. 
This corresponds to the case  $\Re \nu_2=0$ ($\nu_2\neq 0$), with the identification $\sigma^{(0)}=1-\nu_2$ and $a^{(0)}= -{1\over 4} \left[
{e^{i\pi\nu_1} \over 16^{\nu_2-1}}\right]$. The critical behavior along a
radial path is then: 
\be
 y(x)= { 1\over \sin^2\left({\nu\over 2} \ln x - \nu \ln 16 +{\pi
\nu_1\over 2} + \sum_{m=1}^{\infty} c_{0m}(\nu) \left[\left({e^{i\pi \nu_1}
\over 16^{i\nu}}\right)x^{i\nu}\right]^m\right)}~\bigl(1+O(x)\bigr),  ~~~~x\to 0.
\label{fuchsshimomuraintroduzione}
\ee
The number $\nu$ is real, $ \nu \neq 0$ and  $\sigma^{(0)}=1-i\nu$.  
The series $ \sum_{m=1}^{\infty} c_{0m}(\nu)
  \left[\left({e^{i\pi \nu_1} 
\over e^{i\nu}}\right)x^{i\nu}\right]^m$ converges and defines a
holomorphic and bounded function in the domain ${\cal D}(r;\nu_1, i\nu)$ 
$$ 
|x|<r,~~~~  C_1 - \ln r < \nu \arg x < -2 \ln |x| +C_2 + \ln r 
$$ 
 Note   
that not all the values of $\arg x$
are allowed, namely 
$ C_1-\ln r<\nu \arg(x)$.   
 Our belief is that 
 $y(x)$ may have  movable poles   if we extend the range of $\arg x$.  
We are
 not able to prove it in general, but we will give  
an example in section \ref{beyond}.
 
\vskip 0.2 cm 

 We finally remark that the critical behavior of Painlev\'e transcendents can
 also be investigated using  a representation due to S.Shimomura \cite{Sh}
 \cite{IKSY}. 
We will review 
this representation in the paper. However, the connection problem in this
representation was not solved.

\vskip 0.2 cm 
To summarize, 
 in this paper   we give an extended and unified 
picture of both elliptic and Shimomura's representations  and
Dubrovin-Mazzocco's  works, showing that the transcendents obtained in these 
three different ways all are included in the wider 
class  of 
Theorem 1. In this way we solve the
 connection problem for elliptic and Shimomura's representations by virtue of 
Theorem 2. Finally, Theorem 3 provides the oscillatory behavior along 
radial paths when $\Re \sigma^{(0)}=1$.

\section{Monodromy Data and Review of Previous Results}\label{Monodromy Data and Review of Previous Results}

Before giving further details about the result stated above, we review the 
connection between $PVI_{\mu}$  and the theory of isomonodromic deformations. We also 
give a detailed expositions of the results of \cite{DM} \cite{M}. 

    The equation $PVI_{\mu}$ is equivalent to the equations of 
    isomonodromy deformation (Schlesinger equations) of the fuchsian system 
\be{dY\over dz}
= A(z;u)~Y,~~~~~~A(z;u):=
\left[
{A_1(u)\over z-u_1} +{A_2(u)\over z-u_2}+
{A_3(u)\over z-u_3}\right]
\label{fuchsMJ}
\ee
$$u:=(u_1,u_2,u_3), 
~~~~\hbox{tr}(A_i)=\det{A_i}=0,~~~~\sum_{i=1}^3~A_i=-\hbox{diag}(\mu,-\mu)
$$
The dependence of  the system (\ref{fuchsMJ}) on $u$ is isomonodromic, 
as it is explained below.  
 From the system we obtain a transcendent 
$y(x)$ of    $PVI_{\mu}$ as follows: 
$$
    x={u_2-u_1\over u_3-u_1},~~~y(x)= {q(u)-u_1\over u_3-u_1}
$$
where $q(u_1,u_2,u_3)$ is the root of 
$$
    [A(q;u_1,u_2,u_3)]_{12}=0~~~~~~\hbox{if }\mu \neq 0
$$
The case $\mu=0$ is disregarded, because $PVI_{\mu=0} \equiv
PVI_{\mu=1}$. 

Conversely, given a transcendent $y(x)$  the 
system (\ref{fuchsMJ}) associated to it is obtained 
as follows.  Let's define 
$$
   k=k(x,u_3-u_1):= {k_0~\hbox{exp}\left\{ (2\mu-1)~ \int^x d\zeta~
                                                     {y(\zeta)-\zeta\over
   \zeta(\zeta-1) } \right\} \over (u_3-u_1)^{2\mu-1}},~~~~k_0\in {\bf
                                                     C}\backslash\{0\}. 
$$
We have 
 \be
     A_i=-\mu \pmatrix{ \phi_{i1}\phi_{i3} & -\phi_{i3}^2 \cr
                      \phi_{i1}^2         & \phi_{i1}\phi_{i3} \cr
                    } ,~~~~i=1,2,3,
\label{definite in extremis}
\ee
where 
$$
\phi_{13}= i{\sqrt{k} \sqrt{y} \over \sqrt{u_3-u_1}  \sqrt{x}}
$$
$$
\phi_{23}=- {\sqrt{k} \sqrt{y-x} \over \sqrt{u_3-u_1} \sqrt{x} \sqrt{1-x}}
$$
$$
\phi_{33}=i {\sqrt{k} \sqrt{y-1} \over \sqrt{u_3-u_1}\sqrt{1-x}}
$$

$$
 \phi_{11}= {i\over 2 \mu^2} {\sqrt{u_3-u_1}
 \sqrt{y} \over \sqrt{k(x)} \sqrt{x} }
 \left[ A \left(B +{2\mu\over y}\right)+\mu^2(y-1-x) \right]
$$
$$
 \phi_{21}= 
-{1\over 2 \mu^2} {\sqrt{u_3-u_1} \sqrt{y-x} \over \sqrt{k(x)} \sqrt{x} 
\sqrt{1-x} }
 \left[ A \left(B +{2\mu\over y-x}\right)+\mu^2(y-1+x) \right]
$$
$$
 \phi_{31}= {i\over 2 \mu^2} {\sqrt{u_3-u_1} 
\sqrt{y-1} \over \sqrt{k(x)} \sqrt{1-x} }
 \left[ A \left(B +{2\mu\over y-1}\right)+\mu^2(y+1-x) \right]
$$
$$
 A=A(x):= {1\over 2} \left[ {dy\over dx} x(x-1)-y(y-1)\right]
,~~~~B=B(x):= {A\over y(y-1)(y-x)}
$$
Any branch of the square roots can be chosen. 
For a derivation of the above formulae, see \cite{JM1}, \cite{Dub1} and \cite{guz1}. 

\vskip 0.2 cm

\begin{figure}
\epsfxsize=12cm
\centerline{\epsffile{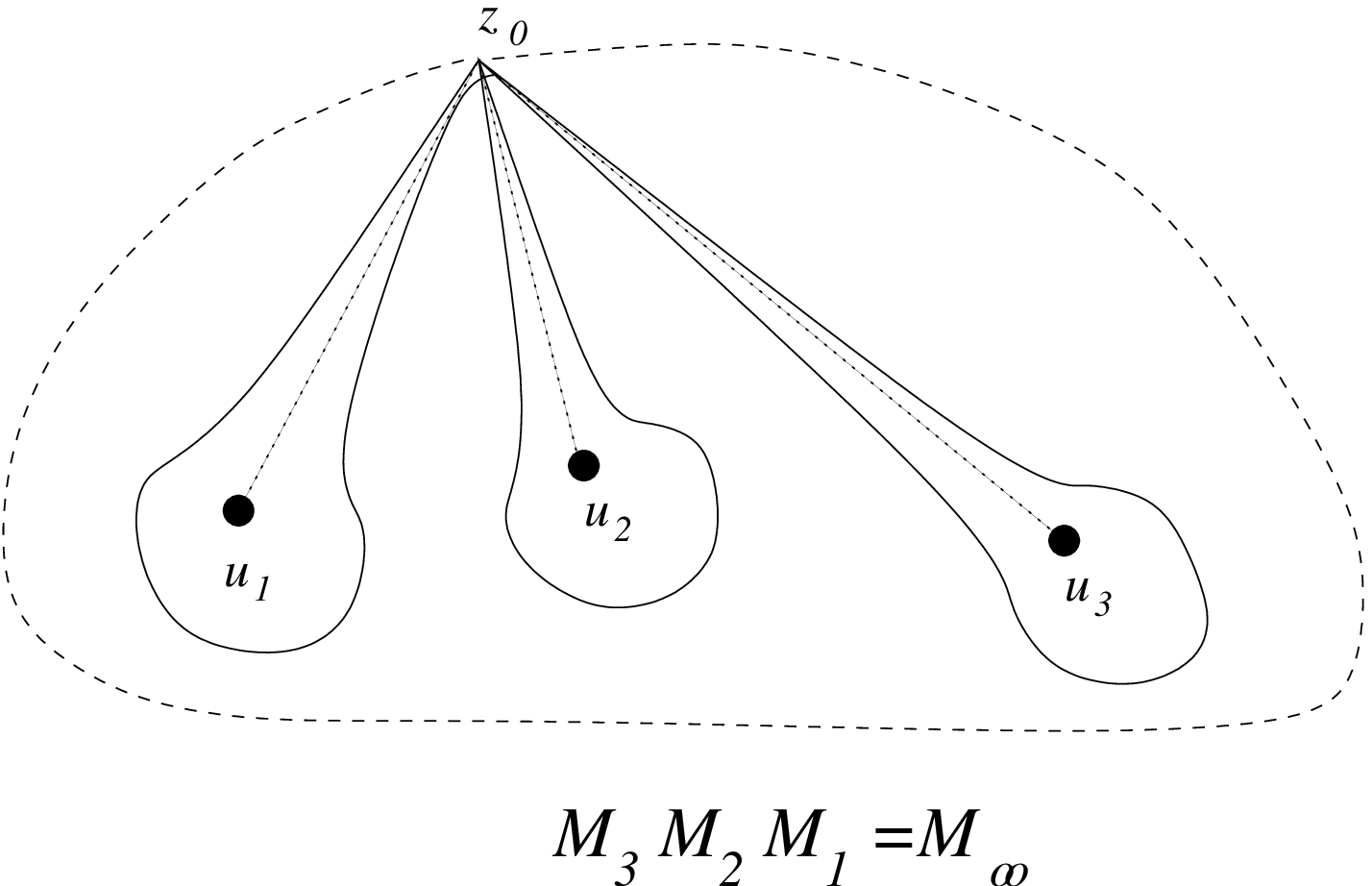}}
\caption{Choice of a basis in $\pi_0({\bf C}\backslash\{u_1,u_2,u_3\}) $}
\label{figura8}
\end{figure}

The system (\ref{fuchsMJ}) has fuchsian singularities at $u_1$, $u_2$,
 $u_3$. Let us fix a branch $Y(z,u)$ of a fundamental matrix solution by
 choosing branch cuts  in the $z$  plane and a basis of loops in
  $\pi({\bf C}\backslash \{u_1,u_2,u_3\};z_0)$, where $z_0$ is a
 base-point. Let   $\gamma_i$ be a basis of loops 
encircling  counter-clockwise the point $u_i$, $i=1,2,3$. See figure
 \ref{figura8}. Then  
 $$ 
Y(z,u) \mapsto Y(z,u)M_i,~~~ i=1,2,3,
~~~~\det M_i\neq 0, 
$$
if $z$ goes around a loop  $\gamma_i$. Along the 
loop $\gamma_{\infty}:=\gamma_1 \cdot \gamma_2 \cdot
\gamma_3$ we have $Y\mapsto Y M_{\infty}$,  $ 
              M_{\infty}= M_3 M_2 M_1
$.  The $2\times 2$ matrices  
$M_i$ are the {\it monodromy matrices}, and they give a
representation of the fundamental group called {\it monodromy 
representation}.   The transformations $ 
 Y^{\prime}(z,u)= Y(z,u) B$, $\det(B)\neq 0
$ 
yields all  possible fundamental matrices, hence the monodromy matrices
of (\ref{fuchsMJ}) are defined up to conjugation 
$$
M_i \mapsto M_i^{\prime}= B^{-1}
M_i B.
$$  
From the standard theory of fuchsian systems it
follows that we can choose  a fundamental solution behaving as follows
\be
 Y(z;u)= \left\{ \matrix{ 
                         \left[I+O({1\over z})\right] ~z^{-\hat{\mu}}
                        z^{R}~C_{\infty},~~~~z\to \infty \cr \cr
                        G_i~\left[I+O(z-u_i)\right]~(z-u_i)^J~C_i,~~~~z\to
                        u_i,~~~i=1,2,3 \cr
                                     } \right.
\label{staraRIMS}
\ee
where $J=\pmatrix{ 0 & 1 \cr 0 & 0 \cr}$,  
$ \hat{\mu}= $diag($\mu,-\mu)$, $G_i J G_i^{-1}= A_i$  and  
$$
     R=\left\{\matrix{ 
~~~~~~
                        0~~~~, ~~~~~~~~~~~~~~~~\hbox{ if } 2\mu \not
\in {\bf Z } \cr \cr
             \left. \matrix{  
           \pmatrix{ 0 & R_{12} \cr 0 & 0 \cr } ,~~~~~~~~\mu>0\cr
                \pmatrix{ 0 & 0 \cr R_{21} & 0 \cr } ,~~~~~~~~\mu<0 \cr 
     } \right\}
                      \hbox{ if } 2\mu\in{\bf Z} \cr   
         }\right. 
$$
The entries  $R_{12}$, $R_{21} $ are  determined by the  matrices $A_i$.  
Then $M_i=C_i^{-1} e^{2\pi i J} C_i$, $M_{\infty}= C_{\infty}^{-1}
e^{-2\pi i \hat{\mu}} e^{2\pi i R} C_{\infty}$.

The dependence of the fuchsian system on $u$ is isomonodromic. This
 means that for small deformations of $u$ the monodromy matrices do
 not change \cite{JM1} \cite{IN}. 
 Small deformation means that $x=(u_3-u_1)/(u_2-u_1)$
 can move in the $x$-plane provided it does not go around  complete
 loops around $0,1,\infty$. 
If the  deformation is not small, the monodromy matrices
 change according to an action of the pure braid group, as it is discussed in \cite{DM}. 

\vskip 0.2 cm 
We consider a branch $y(x)$ of a transcendent and we 
 associate to it the fuchsian
system through the formulae (\ref{definite in extremis}). A branch is fixed  
by the choice of branch cuts, like $\alpha<\arg(x)<\alpha+2\pi$ and
 $\beta<\arg(1-x)<\beta+2\pi$, $\alpha,\beta\in{\bf R}$. Therefore, the 
monodromy matrices of the fuchsian system do not change as $x$ moves in the 
cut plane.  In other words, it is well defined a correspondence which 
associates a monodromy representation to a branch of a transcendent.

 Conversely, the problem of finding 
 a branch of a 
 transcendent  for given monodromy matrices (up to
 conjugation) is the problem of finding a fuchsian system (\ref{fuchsMJ})
 having the given monodromy matrices. This problem is called {\it
 Riemann-Hilbert problem}, or { \it $21^{th}$ Hilbert problem}. For a
 given $PVI_{\mu}$ (i.e. for a fixed $\mu$) there
 is a one-to-one correspondence between a monodromy representation and a
 branch of a transcendent if and only if 
the Riemann-Hilbert problem has a unique
 solution.

\vskip 0.2 cm 
\noindent
{\bf $\bullet$ Riemann-Hilbert problem (R.H.)}: find the 
 coefficients $A_i(u)$, $i=1,2,3$ from the following monodromy data:

a) the matrices  
    $$ 
        \hat{\mu} = \hbox{diag}(\mu,-\mu), ~~~~~\mu\in {\bf
        C}\backslash \{0\}
$$
\vskip 0.2 cm 
$$
     R=\left\{\matrix{ 
~~~~~~
                        0~~~~, ~~~~~~~~~~~~~~~~\hbox{ if } 2\mu \not
\in {\bf Z } \cr \cr
             \left. \matrix{  
           \pmatrix{ 0 & b \cr 0 & 0 \cr } ,~~~~~~~~\mu>0\cr
                \pmatrix{ 0 & 0 \cr b & 0 \cr } ,~~~~~~~~\mu<0 \cr 
     } \right\}
                      \hbox{ if } 2\mu\in{\bf Z} \cr   
         }\right. ,~~~~b \in {\bf C} $$  

b)  three poles $u_1$, $u_2$, $u_3$, a base-point and a base of loops in
 $\pi({\bf C}\backslash \{u_1,u_2,u_3\};z_0)$. See figure \ref{figura8}.

c) three  monodromy matrices $M_1$, $M_2$, $M_3$ relative to the loops
 (counter-clockwise) 
 and a matrix $M_{\infty}$ similar to $e^{-2\pi i \hat{\mu}} e^{2\pi iR}$, 
 and satisfying 
$$ 
\hbox{tr}(M_i)=2,~~~\hbox{det}(M_i)=1,~~~~i=1,2,3 $$
$$ 
  M_3~M_2~M_1 = M_{\infty}
$$
\be
M_{\infty}= C_{\infty}^{-1}~e^{-2\pi i \hat{\mu}}e^{2 \pi i  R}~C_{\infty}
\label{conn1}
\ee
where $C_{\infty}$ realizes the similitude. We also choose the indices
of the problem, namely we fix  ${1 \over 2\pi i}\log M_i$ as
follows: let
$$ 
   J:=\pmatrix{0 & 1 \cr 0 & 0 \cr}
$$ 
We require there exist 
three {\it connection matrices} $C_1$, $C_2$, $C_3$ such that 
\be
     C_i^{-1} e^{2 \pi i J} C_i = M_i,~~~~~i=1,2,3
\label{conn2}
\ee
and we look for a matrix valued meromorphic
function $Y(z;u)$  such that 
\be
   Y(z;u)= \left\{ \matrix{ 
                        G_{\infty} \bigl(I+O({1\over z})\bigr) ~z^{-\hat{\mu}}
                        z^{R}~C_{\infty},~~~~z\to \infty 
\cr 
\cr
                        G_i\bigl(I+O(z-u_i)\bigr)~(z-u_i)^J~C_i,~~~~z\to
                        u_i,~~~i=1,2,3 \cr
                                     } \right.
\label{solrh}
\ee
$G_{\infty}$ and $G_i$ are  invertible matrices depending on $u$.  The
coefficient of the fuchsian system are then given by $
 A(z;u_1,u_2,u_3):= {d Y(z;u) \over dz} Y(z;u)^{-1}
$.

\vskip 0.2 cm
A $2 \times 2$ R.H. is always solvable at a fixed $u$ \cite{AB}. 
As a 
function of $u=(u_1,u_2,u_3)$, the solution $A(z;u_1,u_2,u_3)$ 
extends to a meromorphic function on the universal covering of 
${\bf C}^3\backslash \cup_{i\neq j}\{u~|~u_i=u_j \}$. The monodromy
matrices are considered up to the conjugation
\be
 M_i \mapsto M_i^{\prime}=B^{-1} M_i B, ~~~~~~~\det{B}\neq 0,~~~~i=1,2,3, \infty
\label{conjj}
\ee
and the coefficients of the fuchsian system itself are considered up
to conjugation $A_i \mapsto F^{-1} A_i F$ $(i=1,2,3)$, by an invertible matrix
$F$.  Actually, two conjugated
fuchsian systems admit fundamental matrix solutions with the same
monodromy, and a given fuchsian system defines the monodromy up to
conjugation.  

On the other hand, a triple of monodromy matrices $M_1$, $M_2$, $M_3$
may be realized by two fuchsian systems which are not conjugated.    
This corresponds to the fact that the solutions $C_{\infty}$, $C_{i}$
of (\ref{conn1}), (\ref{conn2})
are not unique, and the choice of different particular solutions may
give rise to   fuchsian systems which are not conjugated. If this is
the case, there is no one-to-one correspondence between monodromy matrices (up
to conjugation) and solutions of $PVI_{\mu}$. It is proved in \cite{M} that:

\vskip 0.2 cm 
{\it The
R.H. has a unique solution, up to conjugation, for $2 \mu \not \in {\bf
Z}$ or for $2 \mu \in {\bf Z}$ and $R\neq 0$}. $^2$


\vskip 0.3 cm 

Once the R.H. is solved, the sum of the matrix coefficients $A_i(u)$ of the 
solution $A(z;u_1,u_2,u_3)=\sum_{i=1}^3
{A_i(u)\over z-u_i} $ must be diagonalized (to give
$-$ diag$(\mu,-\mu)$). $^3$
%
%
%
  After that, a branch y(x)  of $PVI_{\mu}$ can be 
computed from $[A(q;u_1,u_2,u_3)]_{12}=0$.   
 The fact that the R.H. has a unique solution for the given monodromy data
 (if  $2 \mu \not \in {\bf
Z}$ or  $2 \mu \in {\bf Z}$ and $R\neq 0$) means that there is a one-to-one
 correspondence between the triple  $M_1$, $M_2$, $M_3$  and the branch
 $y(x)$.

\vskip 0.2 cm 
\noindent
 We review some known results \cite{DM} \cite{M}. 
\vskip 0.15 cm 

1) One  $M_i=I$ if and only if 
$q(u)\equiv u_i$. This does not correspond to a solution of
$PVI_{\mu}$. 

\vskip 0.15 cm 

2) If the $M_i$'s, $i=1,2,3$, commute, then $\mu$ is integer (as it follows
from the fact that the $2\times 2$ matrices with 1's on the diagonals 
commute if and only if
they  can be simultaneously put in upper or
lower triangular form). There are solutions of
$PVI_{\mu}$ only for 
  $$ 
      M_1=\pmatrix{ 1 & i \pi a \cr
                    0 & 1       \cr},~~~M_2=\pmatrix{ 1 & i \pi  \cr
                    0 & 1       \cr},~~~M_3=\pmatrix{ 1 & i \pi(1- a) \cr
                    0 & 1       \cr},~~~~~a\neq 0,1
$$
In this case $R=0$ and $M_{\infty}=I$. For $\mu=1$ the solution is 
  $$ y(x)= {a x \over 1-(1-a)x}$$
and for other integers $\mu$ the solution is obtained from $\mu=1$ by a
birational transformation \cite{DM} \cite{M}. 

\vskip 0.15 cm 

3) Non commuting $M_i$'s. 

 The parameters in the space of the
 monodromy representation,  independent of conjugation, are 
$$
  2-x_1^2:=\hbox{ tr}(M_1 M_2), ~~~2-x_2^2:= \hbox{ tr}(M_2 M_3),~~~
 2-x_3^2:=\hbox{ tr}(M_1M_3)
$$
The triple $(x_0, x_1, x_{\infty})$ in the Introduction  is  
$(x_1,x_2,x_3)$.

\vskip 0.13 cm

3.1) If at least two of the $x_j$'s are zero, then one of the
$M_i$'s is $I$, or the matrices commute. We return to the case 1 or 2.
Note that  $(x_1,x_2,x_3)=(0,0,0)$ in case 2.
\vskip 0.13 cm 

3.2) At most one of the $x_j$'s is zero. We say that the triple
$(x_1,x_2,x_3)$ is {\it admissible}. In this case it is possible to
fully parameterize the monodromy using the triple $(x_1,x_2,x_3)$. 
Namely, there exists a fundamental matrix solution 
such that: 
$$ 
  M_1= \pmatrix{ 1 & -x_1 \cr 
                 0 & 1     \cr} ,~~~~M_2= \pmatrix{ 1 & 0 \cr 
                 x_1 & 1     \cr},~~~~
M_3=\pmatrix{ 1 +{x_2 x_3\over x_1} & -{x_2^2\over x_1} \cr
               {x_3^2\over x_1}     &  1-{x_2 x_3\over x_1} \cr
             },
$$
if $x_1\neq 0$. If $x_1=0$ we just choose a similar parameterization
starting from $x_2$ or $x_3$. 
The relation 
$$ M_3M_2M_1 \hbox{ similar to } e^{-2\pi i \hat{\mu} }e^{2\pi i R}$$ 
implies 
$$ 
    x_1^2+x_2^2+x_3^2-x_1x_2x_3 = 4 \sin^2(\pi \mu)
$$
The conjugation (\ref{conjj}) changes the triple by two signs. Thus
the true parameters for the monodromy data are classes of equivalence
of triples $(x_1,x_2,x_3)$ defined by the change of two signs. 

We have to distinguish three sub-cases of 3.2): 

i) $2\mu \not \in {\bf Z}$. There is a one to one correspondence between (classes of equivalence of) 
 monodromy data $(x_0,x_1,x_{\infty})\equiv (x_1,x_2,x_3)$ 
and the branches of transcendents of $PVI_{\mu}$. 
 The 
 connection problem was solved in \cite{DM} for the class of
transcendents with critical  behavior 
\be
y(x)= a^{(0)} x^{1-\sigma^{(0)}}(1+O(|x|^{\delta})),~~~~x\to 0,
\label{loc1}
 \ee
\be
y(x)= 1-a^{(1)}(1-x)^{1-\sigma^{(1)}} (1+O(|1-x|^{\delta})),~~~~x\to 1,
\label{loc2}
 \ee
\be
y(x)= a^{(\infty)}
 x^{-\sigma^{(\infty)}}(1+O(|x|^{-\delta})),~~~~x\to \infty,
\label{loc3}
\ee
where $a^{(i)}$ and
 $\sigma^{(i)}$ are complex numbers such that $a^{(i)}\neq 0$ and 
 $0\leq \Re \sigma^{(i)}<1$. $\delta$ is a small positive number. 
This behavior is true if $x$ converges
 to the critical points inside a sector in the $x$-plane with vertex
 on the corresponding critical point and finite angular width.  
In \cite{DM} all the 
algebraic solutions are classified and related to the 
finite reflection groups $A_3$, $B_3$, $H_3$.

ii) The case  $\mu$ half integer was studied in \cite{M}. 
There is an infinite set of {\it Picard type
solutions }  in one to one correspondence to triples of monodromy
data not in the equivalence class of  $(2,2,2)$. 
These solutions form a two parameter family, behave asymptotically as 
the solutions of the case
$i)$, and comprise a denumerable subclass of algebraic solutions. In this 
case $R\neq 0$. For any half integer $\mu\neq {1\over 2}$ there is also a one
parameter family of   
{\it Chazy solutions}.  In this case $R=0$ and  the
one to one correspondence with monodromy data is lost. In fact, they
form an infinite family but  any element of the family corresponds to 
 the class of equivalence of the triple $(2,2,2)$. 
 The result of our paper applies to the 
Picard's solutions with $x_i\neq \pm 2$. 

iii) $\mu$ integer. In this case $R\neq 0$. $^4$
There is a one to one correspondence
between monodromy data  
and  branches. To our knowledge,
this  case 
has not  been studied before our paper. 
There are relevant examples of Frobenius
manifolds included in this case, like the case of Quantum 
Cohomology of ${\bf CP}^2$. For this manifold, 
$\mu=-1$, the triple $(x_1,x_2,x_3)=(3,3,3)$ (see  \cite{Dub2} \cite{guz}) 
 and the real part of
$\sigma$ is equal to 1.

\vskip 0.15 cm 
 In this paper  we find the critical behavior and we solve the
 connection problem for any $\mu \neq 0$ and for all the triples 
$(x_1,x_2,x_3)$ except for the points $
    x_i =  \pm 2$ $\Longrightarrow$ $\sigma^{(i)} = 1$, $i=0,1,\infty
$.


\section{ Critical Behavior -- Theorem 1}\label{Local Behaviour -- Theorem 1}

Theorem 1 has been stated in the Introduction and will be proved in section 
\ref{proof of theorem 1}. 
Here we give some comments about the domain $D(\sigma)$. The superscript of $\sigma^{(i)}, 
a^{(i)}$ will be omitted in this section and we 
  concentrate on a small punctured
neighborhood of $x=0$ ($x=1,\infty$ will be treated 
in section \ref{Singular Points x=1, x=infty  (Connection Problem)}).  
The point $x$ can be read as a 
point in the universal covering of ${\bf C}_0:={\bf
C}\backslash \{0\}$ with $0<|x|<\epsilon$ ($\epsilon<1$). Namely, 
 $x=|x| e^{i
\arg(x)}$, where $-\infty< \arg(x) < +\infty$. Let $\sigma $ be 
such that $\sigma\not \in (-\infty,0)\cup 
[1,+\infty)$. In the Introduction we defined the domains $
   D(\epsilon;\sigma;\theta_1,\theta_2,\tilde{\sigma})$, or 
$D(\sigma;\epsilon)$ for real $\sigma$. Theorem 1  
 holds in these domains; the 
small number $\epsilon$  depends on $\tilde{\sigma}$, $\theta_1$ and
$a$. In the following,  
we may sometimes 
 omit $\epsilon$, $\tilde{\sigma}$, $\theta_i$ and write simply
$D(\sigma)$.

\begin{figure}
\epsfxsize=13cm
\centerline{\epsffile{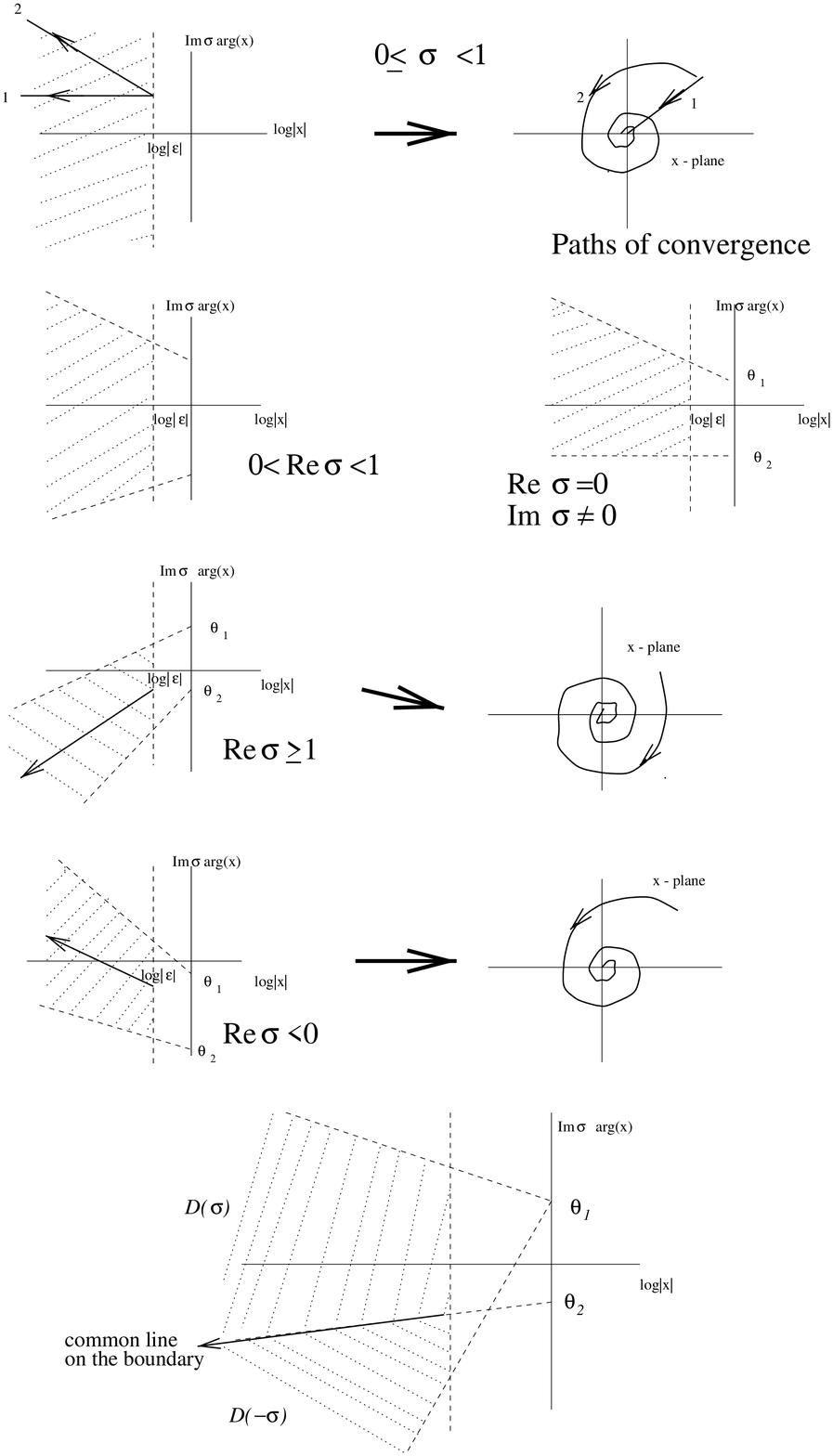}}
\caption{ We represent the domains $D(\epsilon;\sigma;\theta_1,\theta_2)$ in
the $(\ln |x|,\Im\sigma\arg(x))$-plane.
  We also represent some lines along
which $x$ converges to 0. These lines are also represented in the $x$-plane:
they are radial paths or spirals.}
\label{figura1}
\end{figure}

 We observe that  
$
    |{x}^{\sigma}| = |x|^{\sigma^{\prime}(x)}$  where  
 $\sigma^{\prime}(x) := \Re \sigma- {\Im \sigma ~\arg(x) \over \log|x|}
$. 
In particular, for real $\sigma$ we have $\sigma^{\prime}(x)=\sigma$. 
The exponent $\sigma^{\prime}(x)$ satisfies
the restriction $0\leq\sigma^{\prime}(x)<1$ for $x\to 0$, 
 if $x$ lies in the domain, because 
$$
    - {\theta_2 \Im \sigma\over \ln|x|}
\leq \sigma^{\prime}(x) \leq \tilde{\sigma} -{\theta_1 \Im \sigma \over 
\ln|x|},
$$
and  $ - {\theta_2 \Im \sigma\over \ln|x|}\to 0$,  
$ \left(\tilde{\sigma} -{\theta_1 \Im \sigma \over \ln|x|}\right) \to
 \tilde{\sigma}<1$ for $x\to 0$. 
 Figure \ref{figura1} shows the domains in 
the $(\ln|x|,\Im \sigma \arg x)$-plane ( in 
the  $(\ln|x|, \arg x)$-plane if $\Im \sigma=0$).

 In figure \ref{figura1}  
we draw  the  paths along which $x\to 0$. Any regular path is allowed if $\Im
\sigma=0$. If $\Im
\sigma\neq 0$, we considered the family of paths  (\ref{spirale})   connecting a point $x_0\in 
 D(\sigma)$ to $x=0$.  In general, these paths are spirals, represented in figure
\ref{figura1} both in the plane $(\ln|x|, \Im \sigma \arg x)$ and in the 
$x$-plane.   
They are radial paths if $0\leq \Re \sigma<1$ and $  \Sigma= 
\Re \sigma$, because in this case $\arg(x)=\hbox{constant}$.  
But there are only  spiral paths 
 whenever $\Re \sigma<0$ and $\Re \sigma \geq 1$. In particular, the 
 paths 
 $$
\Im \sigma \arg(x) =\Im \sigma \arg(x_0)+ \Re \sigma \log{|x|\over |x_0|} 
$$
  are parallel to one of the  boundary lines of
 $D(\sigma)$ in the plane $(\ln|x|,\Im \sigma~\arg (x))$ and  the 
critical behavior is (\ref{annaalbertostellacecilia}).  
The boundary line is $\Im \sigma \arg)x)= \Re
 \sigma \ln|x|+ \Im \sigma \theta_2$ and it is shared by
$ D(\sigma)$ and   $
 D(-\sigma)$ (with the same $\theta_2$ -- see also Remark 2 of section 
\ref{Parametrization of a branch through Monodromy Data --
Theorem 2}).

\begin{figure}
\epsfxsize=12cm
\centerline{\epsffile{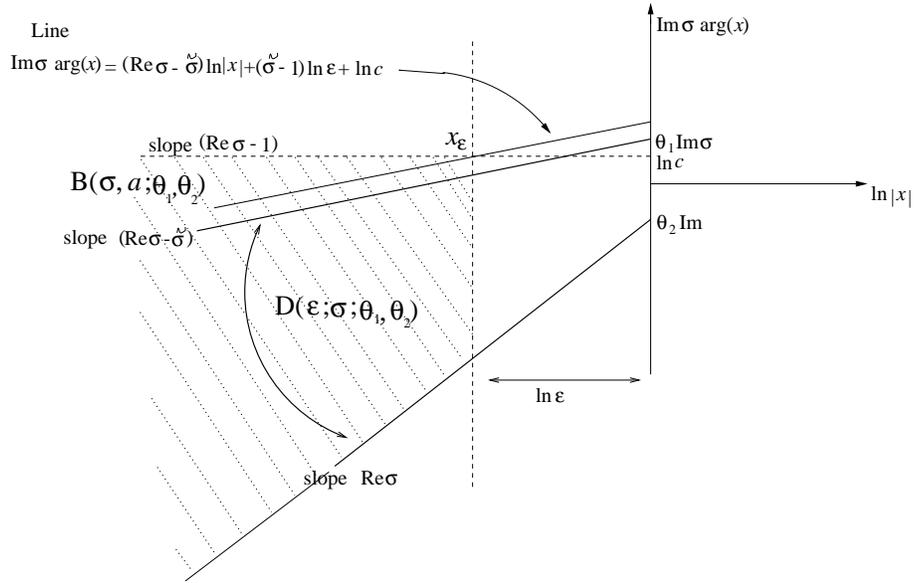}}
\caption{Construction of the domain $B(\sigma,a;\theta_2, \tilde{\sigma})$ for
$\Re \sigma =1$.}
\label{figure13}
\end{figure}

\vskip 0.3 cm 
\noindent
{\bf $\bullet$ Important Remark  on the Domain:}. Consider the domain  
$D(\epsilon;\sigma;\theta_1,\theta_2,\tilde{\sigma})$ for $\Im \sigma \neq 0$. 
In  Theorem 1  we can choose
$\theta_1$ arbitrarily. Apparently, if we increase $\theta_1 \Im \sigma $ the
domain $D(\epsilon;\sigma;\theta_1,\theta_2)$ becomes larger. But $\epsilon$
itself depends on $\theta_1$. In the proof of Theorem 1 (section 
\ref{proof of theorem 1}) we will show that 
$$   
   \epsilon^{1-\tilde{\sigma}}\leq c~e^{-\theta_1 \Im \sigma} 
$$   
where $c$ is a constant, depending on $a$. Equivalently, $\theta_1 \Im \sigma
\leq (\tilde{\sigma}-1)\ln \epsilon + \ln c$. This means that if we increase
$\Im \sigma \theta_1$ we have to decrease $\epsilon$. Therefore, for $x \in   
D(\epsilon;\sigma;\theta_1,\theta_2,\tilde{\sigma})$ we have: 
$$ 
   \Im \sigma \arg (x) \leq (\Re \sigma -\tilde{\sigma}) \ln|x|+ \theta_1 
\Im \sigma \leq (\Re \sigma -\tilde{\sigma}) \ln|x| + (\tilde{\sigma}-1)\ln
\epsilon  +\ln c
$$
We advise the reader to visualize a point $x$ in the plane
$(\ln|x|, \Im\sigma \arg(x))$.  With this visualization in mind, let 
$x_{\epsilon}$ be  the point   $\left\{ \Im \sigma \arg x= 
(\Re \sigma -\tilde{\sigma}) \ln|x|+
(\tilde{\sigma}-1)\ln 
\epsilon  +\ln c \right\} \cap \{|x|=\epsilon\}$ (see figure \ref{figure13}). 
Namely, 
$$ 
\Im \sigma  \arg x_{\epsilon}= (\Re \sigma -1) \ln \epsilon + \ln c
$$ 
This has the following implication. 
 Let $\sigma$, $a$, $\tilde{\sigma}$, $\theta_2$ be fixed. 
 The union  of the domains
$D(\epsilon=\epsilon(\theta_1);\sigma;\theta_1,\theta_2,\tilde{\sigma})$ 
obtained by letting $\theta_1$ vary 
 is
$$ 
  \bigcup_{\theta_1} D\bigl(\epsilon(\theta_1);\sigma;\theta_1,\theta_2,\tilde{\sigma}\bigr)
\subseteq B(\sigma,a;\theta_2, \tilde{\sigma})
$$ 
where 
\be
   B(\sigma,a;\theta_2, \tilde{\sigma}):= \left\{|x|< 1 ~\hbox{ such that }~ \Re \sigma \ln|x|
+\theta_2 \Im \sigma \leq \Im \sigma \arg (x) < (\Re \sigma -1) \ln|x| +\ln
c\right\}
\label{come se non bastasse!!}
\ee
The dependence on $a$ of the domain $B$  is motivated by the fact
that $c $ depends on $a$ (but not on $\theta_1$, $\theta_2$).

 If $0\leq \Re \sigma <1$, the above result is not a limitation on the values
 of $\arg(x)$ in  $D(\epsilon;\sigma;\theta_1,\theta_2,\tilde{\sigma})$, 
 provided that
 $|x|$ is sufficiently small. 

Also in the case $\Re \sigma <0$ there is no
 limitation, because any point $x$, such that $|x|<\epsilon$, can be 
included in  $D(\epsilon;\sigma;\theta_1,\theta_2,\tilde{\sigma})$ 
for a suitable $\theta_2$. In fact, we can always decrease $\Im \sigma \theta_2$ without
 affecting $\epsilon$. 

\begin{figure}
\epsfxsize=12cm
\centerline{\epsffile{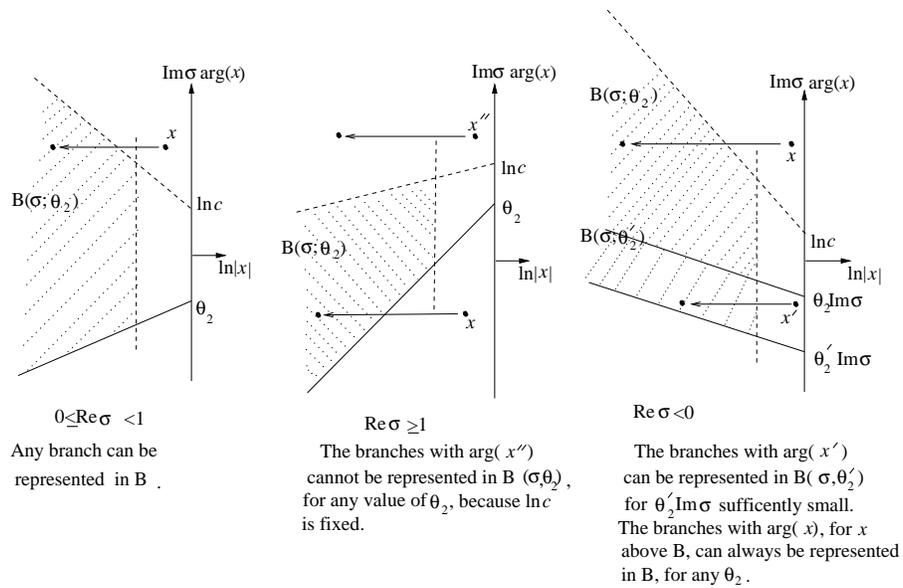}}
\caption{For $\Re \sigma\geq 1$ we can not include all values of $\arg(x)$ in
 $B$}
\label{figure14}
\end{figure}

   But if $\Re
 \sigma \geq 1$, the situation is different. 
Actually, all the  points  $x$ which 
lie above the set  $B(\sigma,a;\theta_2, \tilde{\sigma})$ 
 in the $(\ln|x|, \Im\sigma
 \arg(x))$-plane  can 
 never  be included in any $D(\epsilon;\sigma;\theta_1,\theta_2,
\tilde{\sigma})$. See
 figure \ref{figure14}. This is an important restriction on the domains of 
Theorem 1.

\section{Parameterization of a branch through Monodromy Data --
Theorem 2}\label{Parametrization of a branch through Monodromy Data --
Theorem 2}

The second step in our discussion is to compute the relation between 
the parameters $\sigma$, $a$ of 
Theorem 1, stated for $x=0$,   and the 
monodromy data $(x_0,x_1,x_{\infty})$, to which a unique transcendent  $y(x;x_0,x_1,
x_{\infty})$ is associated. Also in this section, $\sigma^{(0)},a^{(0)}$ are denoted 
$\sigma,a$. The points $x=1,\infty$ are studied in section  
\ref{Singular Points x=1, x=infty  (Connection Problem)}. 

We consider the fuchsian system (\ref{fuchsMJ}) for the
special choice 
$$ 
  u_1=0,~~~u_2=x,~~~u_3=1
$$ 
The labels $i=1,2,3$ will be substituted by the labels $i=0,x,1$, and
the system becomes 
\be 
 {d Y \over dz} = \left[ {A_0(x) \over z} +{A_x(x)\over z-x}+ {A_1(x)
 \over z-1} \right]~Y
\label{seihat}
\ee
Also, the triple $(x_1,x_2,x_3)$ will be denoted by
$(x_0,x_1,x_{\infty})$, as in \cite{DM} and  in the Introduction.  
We consider only   admissible triples 
and  $x_i\neq \pm 2$, $i=0,1,\infty$. 
We recall that an {\it admissible} triple is 
defined in \cite{DM} by the condition that 
only one $x_i$, $i=0,1,\infty$ may be zero. Two admissible triples 
are {\it equivalent}  if their elements  differ just by the change of two
signs and  
\be
x_0^2+x_1^2+x_{\infty}^2-x_0x_1x_{\infty}=4 \sin^2(\pi \mu).
\label{10ii}
\ee

 In the Introduction we called
 $y(x;x_0,x_1,x_{\infty})$ the branch in one to one
 correspondence with an equivalence class of $(x_0,x_1,x_{\infty})$.
  The branch has analytic continuation on the universal covering of 
${\bf P}^1\backslash \{0,1,\infty\}$. We  also  denote this 
continuation by   $y(x;x_0,x_1,x_{\infty})$, where  $x$ is now  a point in 
the universal covering. 

Theorem 2 has been stated in full generality in the Introduction and it will be proved in section \ref{dimth2}.  The result is a generalization 
of the formulae of \cite{DM} to any $\mu \neq 0$ 
(including half-integer $\mu$)  and to 
all $x_i\neq \pm 2$, $i=0,1,\infty$.

   The proof of the theorem 
is valid also for the {\it resonant}
 case $2\mu \in {\bf Z}\backslash \{0\}$.  To read the formulae in this case, 
it is enough to just substitute an integer for
 $2\mu$ in the  formulae i) or I) of the theorem. The cases 
 ii), iii); II), III) do not  occur when $2\mu \in {\bf Z}\backslash \{0\}$.

Note that for $\mu$ integer the case ii), II) degenerates to
 $(x_0, x_1,x_{\infty})=(0,0,0)$
and $a$ arbitrary. This is the case in which the triple is 
 not a good parameterization for the monodromy (not admissible triple).
 Anyway,
we know that in this case there is a one-parameter family of rational
solutions \cite{M}, which are all obtained by a birational
transformation from the family 
$$ 
    y(x)={a x \over 1 -(1-a) x} , ~~~~~\mu =1
$$
At  $x=0$ the the behavior is  
$y(x)= a x (1+O(x))$, and then the limit of Theorem 2 for $\mu \to n \in
              {\bf Z} \backslash \{0\}$ and $\sigma=0$ yields the
above one-parameter family. Recall that $R=0$ in this case.

\vskip 0.2 cm
\noindent
{\it Remark 1: } We repeat the remark to Theorem 2 we 
 made in the Introduction; namely, 
 the  equation  (\ref{LERBAVOGLIO}) does not 
determine $\sigma$ uniquely. We decided to choose $\sigma$ such that $
0\leq \Re \sigma \leq 1
$, so that all all the solutions of (\ref{LERBAVOGLIO})
 are $
\pm \sigma +2 n,$ $n\in {\bf Z}$.  
If $\sigma$ is real, we can only 
choose  $ 0\leq \sigma<1$. 
With this convention, there is a one to one correspondence between 
$(\sigma,a)$ and 
(a class of equivalence of) an admissible triple $(x_0,x_1,x_{\infty})$. 

We observed that  $a= a(\sigma;x_0,
x_1,x_{\infty})$ is affected by  $
\pm \sigma +2 n,$ $n\in {\bf Z}$.  Hence, $y(x;x_0,x_1,x_{\infty})$ 
has different critical behaviors in different domains 
$D( \pm \sigma +2 n)$ . Outside their union, we expect movable poles.

\vskip 0.2 cm
\noindent
{\it Remark 2:}
The domains $D(\sigma)$ and $D(-\sigma)$, with the same $\theta_2$,
intersect along the common boundary  $  \Im \sigma \arg(x)= 
 \Re \sigma \log|x| +
 \theta_2 \Im \sigma$ (see figure \ref{figura1}). 
 The critical behavior of  $y(x;
 x_0,x_1,x_{\infty})$  along the common boundary  is given  in
 terms of $(\sigma(x_0), a(\sigma;x_0,x_1,x_{\infty}))$ and 
$(-\sigma(x_0), a(-\sigma;x_0,x_1,x_{\infty}))$  respectively. 
According to Theorem 1, the critical behaviors in $D(\sigma)$ and 
$D(-\sigma)$ are different, but they become equal on the common boundary. 
 Actually,  along the boundary of $D(\sigma)$ the behavior 
is given by (\ref{annaalbertostellacecilia}), which we rewrite as:  
 $$
 y(x) =
 A(x;\sigma,a(\sigma)) ~x~\left(1+O(|x|^{\delta})\right)$$
 where $\delta$ is a small number between 0 and 1 and 
$$
 A(x;\sigma, a(\sigma))= a (C e^{i\alpha(x;\sigma)})^{-1}+{1\over 2} +{1\over
 16 a} 
 Ce^{i\alpha(x;\sigma)} 
$$    
$$x^{\sigma}=
 Ce^{i\alpha(x;\sigma)},~~~~~C=e^{-\theta_2\Im\sigma},~~~~
~\alpha(x;\sigma)=\Re\sigma\arg(x) 
+\Im\sigma\ln|x|\bigl|_{\Im \sigma \arg(x) =
 \Re \sigma \log|x| + \theta_2 \Im \sigma }
$$   
We observe that  $\alpha(x;-\sigma)=-\alpha(x;\sigma)$. At the end of section 
\ref{dimth2}  we prove that 
$a(\sigma) = {1 \over 16 a(-\sigma)}$. This  implies
that 
     $$ A(x;- \sigma, a(-\sigma))= A(x;\sigma,a(\sigma))$$
Therefore, the critical  behavior, as prescribed
by Theorem 1 in $D(\sigma)$ and $D(-\sigma)$, is the same along the
common boundary of the two domains. 

\vskip 0.3 cm
We end the section with the following

\vskip 0.2 cm
\noindent
{\bf Proposition} [unicity]: {\it Let $\sigma \not \in 
 (-\infty,0)\cup [1,+\infty)$ and $a\neq 0$.  Let $y(x)$ be a solution 
of $PVI_{\mu}$ such that 
 $y(x)= a x^{1-\sigma}$(1+higher order terms)  as $x\to 0 $ in the 
 domain $D(\epsilon;\sigma)$. Suppose that the triple $(x_0,x_1,x_{\infty})$ 
computed by the formulae of Theorem 2 in terms of $\sigma$ and $a$ 
is admissible.  Then, $y(x)$ coincides with $y(x;\sigma,a)$ of
Theorem 1} 
\vskip 0.2 cm
\noindent
{\it  Proof:} see section \ref{dimth2}.



 \section{ Other Representations of the Transcendents -- 
Theorem 3}\label{beyond}

We need to further investigate the critical behavior close to $x=0$, 
in order to extend the results of Theorem 1 for $x\to 0$ along paths not 
allowed by the theorem. In this section we discuss the critical behavior of 
the elliptic representation of Painlev\'e transcendents. 
 According to Remark 1 of section 
\ref{Parametrization of a branch through Monodromy Data --
Theorem 2}   we
 restrict to  $   0 \leq \Re \sigma \leq 1 $ for $  \Im \sigma\neq 0$, 
  or $   0\leq \sigma <1 $ for $ \sigma$ real. 

\begin{figure}
\epsfxsize=15cm
\centerline{\epsffile{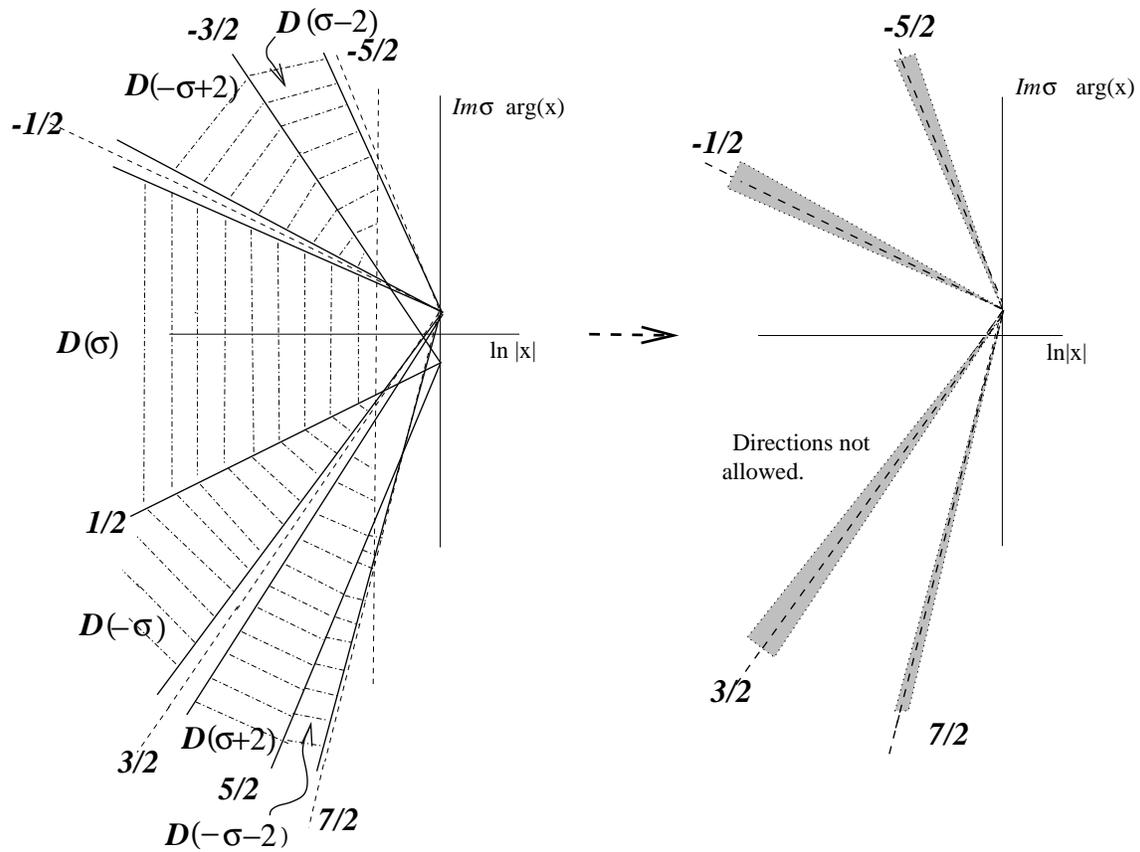}}
\caption{ Domains for $\sigma={1\over 2} +i \Im \sigma$. The numbers close to
the lines are their slopes. 
The small sectors 
around the dotted lines 
represented in the right figure are not contained in the union of the
domains. If $x\to 0$ along a direction which ends in one of these sectors, we
do not know the behavior of the transcendent.}
\label{figure8}
\end{figure}

 In figure \ref{figure8} (left) we draw domains  $D(\sigma)$, $D(-\sigma)$,
$D(-\sigma+2)$, $D(2-\sigma)$, etc, where $y(x;x_0,x_1,x_{\infty})$ has 
different critical behaviors.  
Some small
sectors remain uncovered by the union of the domains (figure \ref{figure8}
(right)). If $x\to 0$ inside these sectors, we do not know the behavior of
the transcendent. For example, if $\Re \sigma =1$, a 
 radial path converging to $x=0$ will end up in a forbidden small 
sector (see also figure \ref{figure12} for the case $\Re \sigma=1$).

\begin{figure}
\epsfxsize=12cm
\centerline{\epsffile{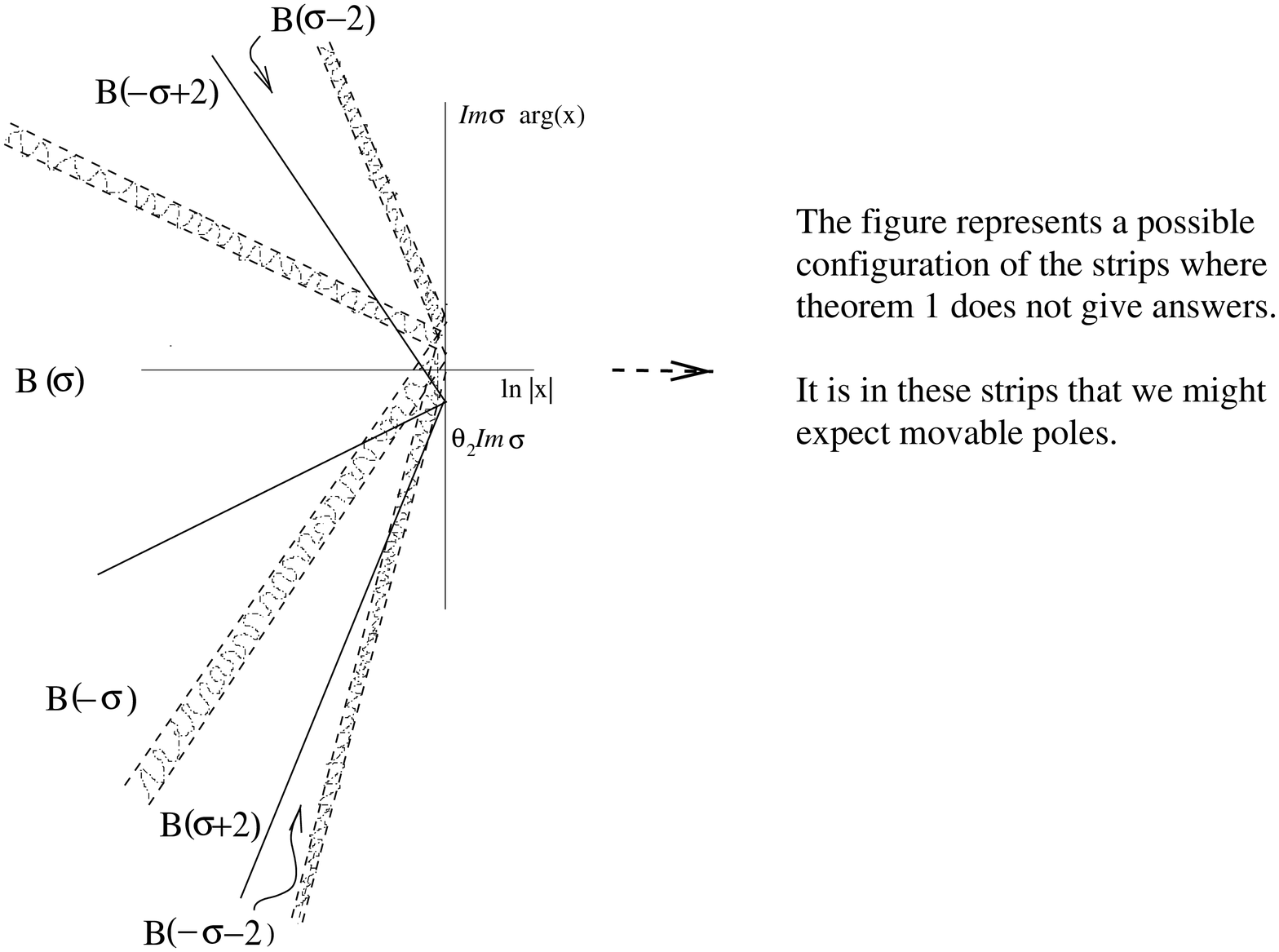}}
\caption{ }
\label{figura18}
\end{figure}

 If we draw, for the same $\theta_2$,  
 the domains $B(\sigma)$, $B(-\sigma)$, $B(-\sigma+2)$, etc, 
 defined in 
 (\ref{come se non bastasse!!}) we obtain strips in the $(\ln|x|,
 \Im\sigma\arg(x) )$-plane which are {\it certainly forbidden} to Theorem 1
 (see figure \ref{figura18}).
 In the strips we know nothing about the transcendent. We guess
 that there might be poles there, as we verify in one example later.

 \vskip 0.2 cm 

What is the behavior along the directions not allowed by Theorem 1? 
 In the very particular case $(x_0,x_1,x_{\infty}) \in \{ (2,2,2), (2,-2,-2),
(-2,-2,2), (-2,2,-2)\}$ it is known that 
$PVI_{\mu=-1/2}$ has a 1-parameter family of {\it classical} 
solutions. The  critical 
 behavior of a branch for {\it radial} 
convergence to the critical points 0, 1 ,$\infty$ was computed in \cite{M}:  
$$ 
   y(x)=\left\{\matrix{          
                      -\ln(x)^{-2}(1+O(\ln(x)^{-1})),~~~~~x\to 0 \cr\cr
                  1+ \ln(1-x)^{-2}(1+O( \ln(1-x)^{-1})),~~~~~x\to 1 \cr\cr
                   -x\ln(1/x)^{-2}(1+O(\ln(1/x)^{-1}),~~~~~x\to \infty \cr
     } \right.
$$ 
The branch is specified by $|\arg(x)|<\pi$, $|\arg(1-x)|<\pi$.  
This behavior is completely different from $\sim a(x) x^{1-\sigma}$ 
 as $x\to 0$. Intuitively, as $x_0$ approaches the value 2,  
$1-\sigma$ approaches 0 and the decay of $y(x)\sim a x^{1-\sigma}$ becomes 
logarithmic. These solutions were called {\it Chazy solutions} in \cite{M}, 
because they can be computed as functions of solutions of the 
 Chazy equation. 

 \vskip 0.3 cm

 This section is devoted to the investigation of
the critical  behavior at $x=0$ 
in the regions not allowed in  theorem 1.


\subsection{Elliptic Representation}\label{Elliptic Representation}

 The transcendents of  $PVI_{\mu}$ can be represented in the elliptic form
 \cite{fuchs} 
$$ 
  y(x)= \wp \left({u(x)\over 2};\omega_1(x),\omega_2(x)\right)+{1+x\over 3}
$$
where $\wp(z;\omega_1,\omega_2)$ is the Weierstrass elliptic function of
 half-periods $\omega_1$, $\omega_2$. $ u(x)$ 
 solves the non-linear differential
 equation 
\be
  {\cal L}(u)= {\alpha\over x(1-x)} {\partial \over \partial u} \left[
    \wp\left( {u\over
    2};\omega_1(x),\omega_2(x)\right)\right],~~~\alpha={(2\mu-1)^2\over 2}
\label{difficilissima}
\ee
where the differential linear  operator ${\cal L}$ applied to $u$ is 
$$ 
   {\cal L}(u):= 
x(1-x)~{d^2 u\over dx^2}+(1-2x)~{du\over dx}-{1\over 4}~u.  
$$
 The half-periods are two independent solutions of the hyper-geometric 
equation ${\cal L}(u)=0$:  
$$
\omega_1(x):= {\pi \over 2} ~F\left( x\right),~~~~
\omega_2(x):= -{i\over 2} [F(x) \ln x +F_1(x)] 
$$
where $F(x)$ is the hyper-geometric function 
$$
F(x):=F\left({1\over 2},{1\over 2},1;x\right)=  
\sum_{n=0}^{\infty}{ \left[\left({1 \over 2}\right)_n\right]^2
  \over (n!)^2 } x^n,
$$
and 
$$
F_1(x):=  
\sum_{n=0}^{\infty}{ \left[\left({1 \over 2}\right)_n\right]^2
  \over (n!)^2 } 2\left[ \psi(n+{1\over 2}) - \psi(n+1)\right]
x^n
$$
$$\psi(z) = 
{d \over dz}\ln \Gamma(z),~~~ \psi\left({1\over 2}\right) = -\gamma -2 \ln
2,~~~ \psi(1)=-\gamma,~~~\psi(a+n)=\psi(a)+\sum_{l=0}^{n-1} {1\over a+l}.
$$

The solutions $u$ of (\ref{difficilissima}) were  not  studied in the 
 literature,  
so we did that and we 
 proved a general result in Theorem 3. But first,  
 we give a special example, already known to Picard.

 \vskip 0.3 cm
\noindent
{\it Example:} The equation $PVI_{\mu=1/2}$ has a two parameter family of 
solutions discovered  
 by Picard \cite{Picard} \cite{Okamoto}  \cite{M}. It is easily obtained
from (\ref{difficilissima}). Since  $\alpha=0$, $u$ solves the  hyper-geometric
equation ${\cal L}(u)=0$ and has the general form 
$$
{u(x)\over 2}:= \nu_1 \omega_1(x)+\nu_2 \omega_2(x),~~~ \nu_i \in {\bf C},
~~~0\leq \Re\nu_i <2,~~~(\nu_1,\nu_2)\neq (0,0),
$$
 A branch of $y(x)$ is specified by a 
branch of $\ln x$ in $\omega_2(x)$.  
The monodromy data computed in \cite{M} are 
$$ x_0=-2\cos \pi r_1,~~~x_1= -2\cos \pi r_2,~~~x_{\infty}=-2\cos \pi r_3,$$
$$
   r_1={\nu_2\over 2},~~~r_2=1-{\nu_1\over 2},~~~r_3={\nu_1-\nu_2\over 2},~~~~
\hbox{ for } \Re\nu_1>\Re\nu_2$$
$$
   r_1=1-{\nu_2\over 2},~~~r_2={\nu_1\over 2},~~~
r_3={\nu_2-\nu_1\over 2},~~~~
\hbox{ for } \Re\nu_1<\Re\nu_2.
$$
 The modular parameter is now a function of $x$: 
$$ 
      \tau(x) = {\omega_2(x) \over \omega_1(x)}=  {1 \over \pi} 
                    (\arg x - i \ln|x|) +{4i\over \pi}\ln 2+O(x),~~~x \to 0.
$$
We see that $\Im \tau >0$ as $x \to 0$.  Now, if 
\be
\left|\Im {u(x)\over 4 \omega_1}\right|< \Im \tau,
\label{tau}
\ee
 we can expand the 
Weierstrass function in Fourier series.  Condition (\ref{tau}) becomes
$$
     {1\over 2} \left| \Im \nu_1 +{\Im \nu_2 \over \pi} \arg(x)-
        {\Re \nu_2 \over \pi} \ln|x|+{4\ln 2\over \pi} \Re\nu_2\right|<
        -{\ln|x| \over  
                             \pi} +{4\ln 2\over \pi} +O(x),~~~~\hbox{ as }
x\to 0
$$
For $\Im \nu_2\neq 0$ this can be written as follows: 
\be 
   (\Re \nu_2 +2) \ln|x| +c_1 < \Im \nu_2 \arg(x) 
< (\Re \nu_2 -2) \ln|x| +c_2.
\label{dominio}
\ee
$$
c_1:= - \pi \Im \nu_1-4\ln 2 ~(\Re\nu_2+2),~~~~~
c_2:=  - \pi \Im \nu_1-4\ln 2 ~(\Re \nu_2 -2)
$$
On the other hands, if $\Im \nu_2=0$ any value of $\arg x$ is allowed.  
The Fourier expansion is 
$$
y(x) = {x+1\over 3} +{1 \over F(x)^2 }\left[  {1 \over \sin^2\left( 
-{1\over 2} [i \nu_2(\ln(x) +{F_1(x) \over F(x)})-\pi \nu_1] \right)}
-{1\over 3} + \right.
$$
$$\left.
 +8 \sum_{n=1}^{\infty} { x^{2n} \over e^{-2 n {F_1(x)\over 
F(x)}}
- x^{2n}} \sin^2\left( 
-{n\over 2} [i \nu_2(\ln(x) +{F_1(x) \over F(x)})-\pi \nu_1] \right) \right]
$$
$$
= {x \over 2} +(1 -{x\over 2} +O(x^2)) \left[ {1\over \sin^2\left( 
-{1\over 2} [i \nu_2(\ln(x) +{F_1(x) \over F(x)})-\pi \nu_1] \right)}+
\right.
$$
$$
\left.
 -{1 \over 4} \left[{e^{i\pi \nu_1} \over 
16^{\nu_2-1}}\right]^{-1} x^{2-\nu_2} +
O(x^2+ x^{3-\nu_2}+ x^{4 - \nu_2})  \right], ~~
~~~x\to 0 \hbox{ in the domain } 
(\ref{dominio})
$$
 As far as { \it radial} convergence is concerned, we have:

\vskip 0.2 cm
a) $0<\Re \nu_2<2 $,
$$
     {1 \over \sin^2(...)}= -{1 \over 4} \left[{e^{i\pi \nu_1} \over 
16^{\nu_2-1}} \right]~x^{\nu_2} ~\bigl( 1 +O(|x^{\nu_2}|)\bigr),$$
and so 
\be
y(x)=  \left\{-{1 \over 4} \left[{e^{i\pi \nu_1} \over 
16^{\nu_2-1}} \right]~x^{\nu_2}+ {1\over 2} x 
-{1 \over 4} \left[{e^{i\pi \nu_1} \over 
16^{\nu_2-1}}\right]^{-1} x^{2-\nu_2}\right\} ~\left(1 +O(x^{\delta}) \right),
~~~~\delta>0. \label{behaviour}
\ee
This is the same critical behavior of Theorem 1. By virtue of the 
Proposition of section 
\ref{Parametrization of a branch through Monodromy Data --
Theorem 2}, 
 the transcendent here coincides with $y(x;\sigma,a)$  of Theorem 1 if we 
 identify $1- \sigma $ with $\nu_2$ for $0<\Re \nu_2<1$, or  with 
$2-\nu_2$ for 
$1<\Re \nu_2<2$. In the case $\Re \nu_2 =1$ the three terms $x^{\nu_2}$, 
$x$, $x^{2-\nu_2}$ have the same order and we find again the behavior 
(\ref{cecilia}) of 
Theorem 1 (oscillatory case):
$$
y(x) = \left\{a x^{\nu_2} +{x \over 2} +{1\over 16 a} x^{2-\nu_2}
\right\}(1+O(x^{\delta}))= a x^{\nu_2} \left\{1
+{1\over 2a} x^{-i \Im \nu_2}+ {1\over 16 a^2}
 x^{-2i \Im \nu_2}\right\}(1+O(x^{\delta})),
$$
where $a= -{1\over 4} \left[{e^{i\pi \nu_1} \over 
16^{\nu_2-1}} \right]$.

\vskip 0.2 cm
b) $\Re \nu_2 =0$. Put $ \nu_2 = i \nu$ ( namely, $\sigma =1 -i \nu$ ). The 
domain (\ref{dominio}) is now (for sufficiently small $|x|$):
$$
  2 \ln|x|-\pi \Im \nu_1-8\ln 2<\Im \nu_2 \arg(x) < -2 \ln|x| - \pi \Im
\nu_1+8\ln 2,
$$
or 
 \be
 2 \ln|x|+\pi \Im \nu_1-8\ln 2<\Im \sigma \arg(x) < -2 \ln|x| + \pi \Im
\nu_1+8\ln 2.
\label{dominio1}
\ee
For radial convergence we have 
$$
  y(x)= {1 +O(x) \over \sin^2( {\nu\over 2} \ln(x) + {\nu \over 2} {F_1(x) 
\over F(x)} + {\pi \nu_1\over 2})}  + O(x). 
$$ 
This is an oscillating functions, and it may have poles. 
 Suppose for example that $\nu_1$ is real.   
Since $F_1(x)/F(x)$ is a convergent power series ($|x|<1$) with real 
coefficients and defines a  bounded function,  then $y(x)$ has a 
sequence of poles on the positive real axis, converging to $x=0$.
 
In the domain (\ref{dominio1})  spiral convergence of $x$ to zero is also 
allowed and the critical behavior is 
(\ref{behaviour}) because $\arg x$ is not constant. 
 
 Finally, if     
$\nu=0$, namely $\nu_2=0$ (and then $x_0=2$) we have 
$$ 
y(x) = {1\over \sin^2(\pi \nu_1)} (1+O(|x|)).
$$

\vskip 0.2 cm 

\begin{figure}
\epsfxsize=15cm
\centerline{\epsffile{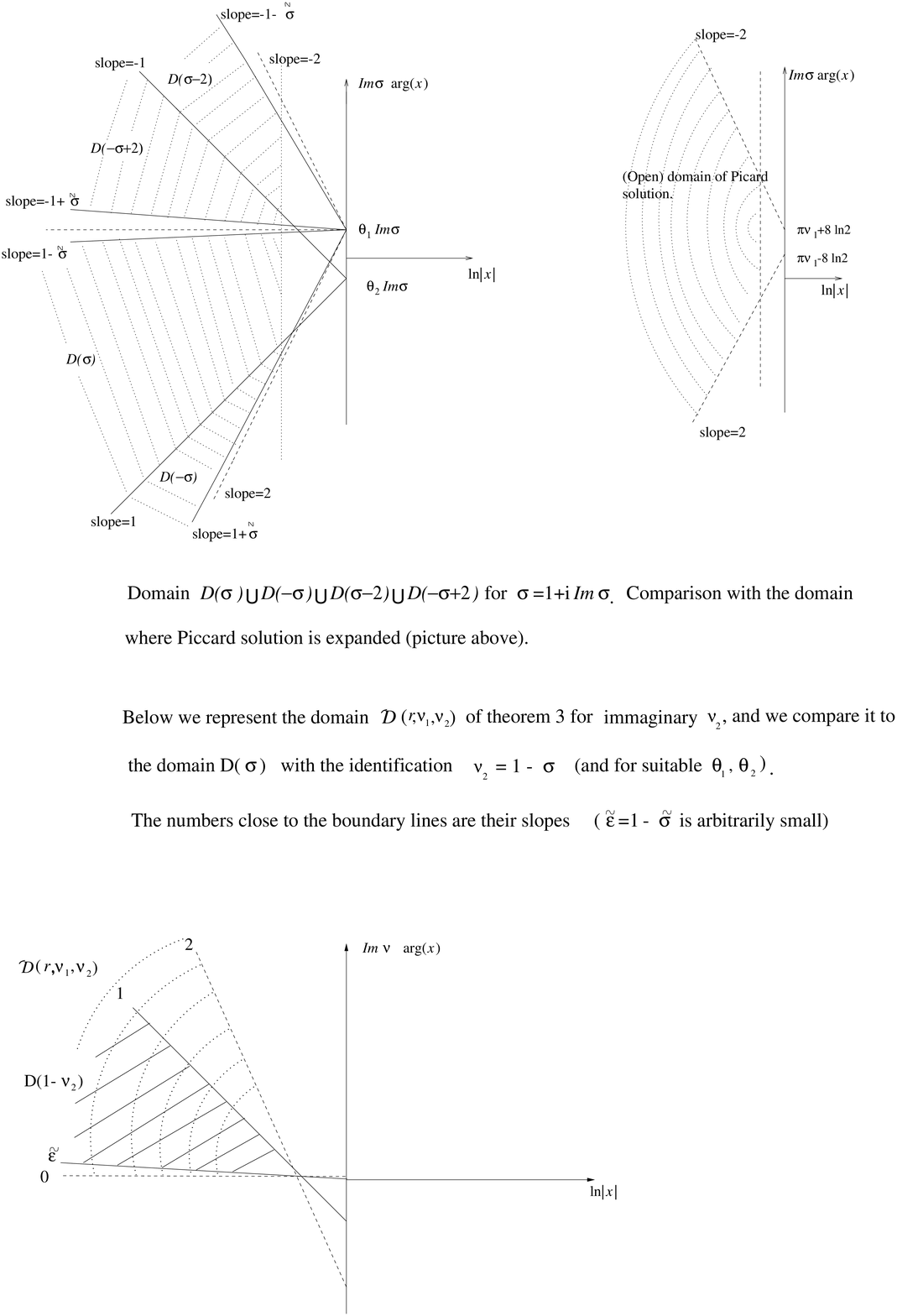}}
\caption{ }
\label{figure12}
\end{figure}

 The case b) in the  above example is good to 
 understand the limitation of Theorem 1 in giving a complete description of 
the behavior of  Painlev\'e transcendents. Actually, Theorem 1  yields  
the behavior (\ref{behaviour}) 
in the domain $D(\sigma)\cup D(-\sigma)$ ($\Re \sigma=1$):
$$
 (1+\tilde{\sigma})\ln|x| + \theta_1 \Im \sigma \leq \Im \sigma \arg x 
\leq (1-\tilde{\sigma}) \ln|x| +\theta_1 \Im \sigma ,
$$ 
 where 
radial convergence to $x=0$ is not allowed.  
On the other hand, the transformations $ \sigma \to\pm( \sigma -2)$, 
gives a further domain $ D(\sigma -2 ) \cup D(-\sigma +2)$:  
$$
  (-1+ \tilde{\sigma})\ln|x| + \theta_1 \Im \sigma \leq \Im \sigma \arg x
\leq -(1+\tilde{ \sigma}) \ln|x| +\theta_1 \Im \sigma, 
$$
but again it is not possible for $x$ to converge to $x=0$ along a radial
path. 
Figure \ref{figure12} 
shows $D(\sigma)\cup D(-\sigma)\cup D(2-\sigma)\cup D(\sigma -2)$. 
Note that a radial  path 
would be allowed if it were possible  to make 
$\tilde{\sigma} \to 1$. 
The interior of the set obtained as the limit 
for 
$\tilde{\sigma}\to 1$  of  
 $D(\sigma)\cup D(-\sigma)\cup D(2-\sigma)\cup D(\sigma -2)$  is like  
 (\ref{dominio1}). Actually, the intersection of (\ref{dominio1}) and 
 $D(\sigma)\cup D(-\sigma)\cup D(2-\sigma)\cup D(\sigma -2)$   is never
empty. On (\ref{dominio1}) the 
 elliptic representation predicts an oscillating
behavior and poles.    
So it is definitely clear that the 
``limit'' of theorem 1 for $\tilde{\sigma} \to 1$ is not trivial.

 \vskip 0.2 cm

\noindent
{\it Remark on the example:} For $\mu$ half integer all the possible 
values of $(x_0,x_1,x_{\infty})$ such that
$x_0^2+x_1^2+x_{\infty}^2-x_0x_1x_{\infty} =4$ 
are covered by Chazy and Picard's solutions, 
with the warning that for $\mu={1 \over 2}$ the image (through birational
transformations) of Chazy solutions  
is $y=\infty$. See \cite{M}.

\vskip 0.3 cm

 We turn to the general case.  The elliptic 
representation has been studied from the point of view of algebraic geometry 
in \cite{Manin2}, but to our knowledge   Theorem 3 and its 
 Corollary, both  stated in the Introduction,  
are  the first  general 
result about its critical behavior appearing  in  the literature. 

 We prove Theorem 3 in 
 section \ref{provaellittica}. 
Here we prove the Corollary. 
The critical behavior is obtained expanding  $y(x)$ in Fourier series: 
\be
   \wp\left({u\over 2};\omega_1,\omega_2\right)= -{\pi^2 \over 12
   \omega_1^2} +{2\pi^2\over \omega_1^2} \sum_{n=1}^{\infty} {n e^{2\pi i n
   \tau} \over 1- e^{2\pi i n \tau}} \left(1- \cos \left(n {\pi u \over 2
   \omega_1} \right)\right) +{\pi^2\over 4 \omega_1^2} ~{1\over
   \sin^2\left( {\pi u \over 4 \omega_1}\right)} 
\label{rondine}
\ee
The expansion can be performed if $\Im \tau(x)>0$ and $\left|\Im
\left({u(x)\over 
\omega_1(x)}\right)\right|<\Im \tau$; these conditions are satisfied in ${\cal
D}(r;\nu_1,\nu_2)$. Let's put $F_1/F=-4\ln 2 +g(x)$, where 
 $g(x)=O(x)$.  
 Taking into account (\ref{rondine}) and  Theorem 3, 
the expansion of $y(x)$ for $x\to 0$ in $  {\cal D}(r;\nu_1,\nu_2)$   is
$$
  y(x)= \left[{1+x\over 3} - {\pi^2\over 12 \omega_1(x)^2}\right]
 +{\pi^2\over
  \omega_1(x)^2 } \sum_{n=1}^{\infty} {n\over 1- \left({e^{g(x)}\over
  16}\right)^{2n} x^{2n}}\left\{
2\left({e^{g(x)}\over
  16}\right)^{2n} x^{2n} +\right.
$$
$$\left. - e^{n(\nu_2+2)g(x)}\left[ {e^{i\pi \nu_1}\over
  16^{2+\nu_2}} x^{2+\nu_2} \right]^n e^{in\pi {v(x)\over \omega_1(x)}} - 
e^{n(2-\nu_2)g(x)} \left[{e^{-i\pi \nu_1}\over
  16^{2-\nu_2}}x^{2-\nu_2}\right]^n e^{-i\pi n {v(x)\over \omega_1(x)} } 
\right\}
+
$$
$$
+ {\pi^2\over 4\omega_1(x)^2}{1\over   \sin^2\left(-i{\nu_2\over 2}
  \ln x +i{\nu_2\over 2} \ln 16 +{\pi \nu_1\over 2} - i{\nu_2 \over 2} g(x)
  +{\pi v(x) \over 2 \omega_1(x)}\right) }
$$ 
We  observe that  $\omega_1(x) \equiv {\pi \over 2} F(x) =
{\pi \over 2} ( 1 +{1\over 4} x +O(x^2))$, 
$ 
{1+x\over 3} - {\pi^2\over 12 \omega_1(x)^2}\equiv {1+x\over 3} -{1\over 3
F(x)} =  {1\over 2} x \bigl(1+O(x)\bigr)
$, $ e^{g(x)}= 1 +O(x)$ and 
$$
e^{\pm i \pi {v(x)\over \omega_1(x)}}= 1+ O\left(|x|+\left|{e^{-i\pi
\nu_1}\over 16^{2-\nu_2}} x^{2-\nu_2} \right|+\left|  
{e^{i\pi \nu_1} \over 16^{\nu_2}} x^{\nu_2}\right| \right)
$$

In order to single out the leading terms, we observe that we are dealing with
the powers $x$, $x^{2-\nu_2}$, 	$x^{\nu_2}$ in ${\cal D}(r;\nu_1,\nu_2)$. If
$0<\nu_2<2$ (the only allowed real values of $\nu_2$) $|x^{\nu_2}|$ is leading
if $0<\nu_2<1$ and $|x^{2-\nu_2}|$ is leading if $1<\nu_2<2$. We have 

$$
 {1\over \sin^2(...)} = -4 {e^{i\pi \nu_1}\over 16^{\nu_2} } x^{\nu_2}
 \left[1+O(|x|+|x^{\nu_2}|+|x^{2-\nu_2}|)\right] 
$$ 
Thus, there exists $0<\delta<1$ (explicitly computable in terms of $\nu_2$)  
such that
$$ 
  y(x)= \left[{1\over 2} x - {1\over 4} \left[{e^{i\pi \nu_1} \over
  16^{\nu_2-1} }  \right] x^{\nu_2} - {1\over 4} \left[{e^{i\pi \nu_1} \over
  16^{\nu_2-1} }  \right]^{-1} x^{2-\nu_2}
  \right](1+O(x^{\delta}))
$$
$$
= \left\{ \matrix{ \sin^2\left({\pi \nu_1\over 2} \right)
 ~x~(1+O(x^{\delta})), ~~~\hbox{ if } \nu_2=1 \cr\cr
 - {1\over 4} \left[{e^{i\pi \nu_1} \over
  16^{\nu_2-1} }  \right] x^{\nu_2}(1+O(x^{\delta}))~~~\hbox{ if } 0<\nu_2<1
  \cr\cr 
- {1\over 4} \left[{e^{i\pi \nu_1} \over
  16^{\nu_2-1} }  \right]^{-1} x^{2-\nu_2}(1+O(x^{\delta}))~~~\hbox{ if }
  1<\nu_2<2 \cr
}\right.
$$ 
This behavior coincides with that  
 of Theorem 1 for $\sigma=0$ in the first case, $\sigma=
1-\nu_2$ in the second, $\sigma=\nu_2-1$ in the third. 

\vskip 0.2 cm

We turn to the case $\Im \nu_2\neq 0$. We consider a path contained in ${\cal
D}(r;\nu_1,\nu_2)$ of equation 
\be
  \Im \nu_2 \arg(x)= (\Re \nu_2 - {\cal V}) \ln|x|+b,~~~~0\leq {\cal V}\leq 2
\label{che il cielo ce la mandi buona}
\ee
with a suitable constant $b$ (the path connects some $x_0\in {\cal
D}(r;\nu_1,\nu_2)$ to $x=0$, therefore $b=\Im \nu_2 \arg x_0 - (\Re \nu_2 - 
{\cal V}) \ln |x_0|$). 
 We have$|x^{2-\nu_2}|=|x|^{2-{\cal V}} e^b$,
$|x^{\nu_2}|= |x|^{\cal V} e^{-b}$ and so  
$$ 
   |x^{\nu_2}|~~\hbox{ is leading for } ~ 0\leq {\cal V}<1,
$$
$$
 |x^{\nu_2}|,~|x|,~|x^{2-\nu_2}|~~\hbox{ have the same order for } ~ 
{\cal V}=1,
$$
$$
|x^{2-\nu_2}|~~\hbox{ is leading for } ~ 1< {\cal V}\leq 2.
$$
If ${\cal V}=0$,  
$$\left| 
{e^{i\pi \nu_1} \over 16^{\nu_2}} x^{\nu_2}\right|<r,~~~~~~\hbox{ but }
~x^{\nu_2}\not\to 0~\hbox{ as } ~x\to 0.
$$
If  ${\cal V}=2$, 
  $$
 \left|{e^{-i\pi \nu_1}\over 16^{2-\nu_2}} x^{2-\nu_2} \right|<r~~~~~~\hbox
{ but } 
~x^{2-\nu_2}\not\to 0~\hbox{ as } ~x\to 0. 
$$ 
This also implies that $v(x)\not \to 0$ as $x\to 0$ 
along the paths with ${\cal V}=0$ or ${\cal V}=2$, while $v(x)\to 0$ for all 
other values $0 <{\cal V}<2$.   We conclude that:

\vskip 0.2 cm 
a) If  $x\to 0$ in ${\cal D}(r;\nu_1,\nu_2)$ along 
(\ref{che il cielo ce la mandi buona}) for ${\cal V}\neq 0,2$, then  
$$ 
 y(x)= \left[{1\over 2} x - {1\over 4} \left[{e^{i\pi \nu_1} \over
  16^{\nu_2-1} }  \right] x^{\nu_2} - {1\over 4} \left[{e^{i\pi \nu_1} \over
  16^{\nu_2-1} }  \right]^{-1} x^{2-\nu_2}
  \right](1+O(x^{\delta})),~~~~0<\delta<1. 
$$ 
The three leading terms have the same order if the convergence is along a path
asymptotic to (\ref{che il cielo ce la mandi buona}) with
${\cal V}=1$. Namely
$$ 
y(x)= 
x~\sin^2\left(i{1-\nu_2\over 2} \ln x +{\pi \nu_1\over 2} +2i(\nu_2-1)\ln 2
\right)~\bigl(1+O(x)\bigr)~~\hbox{ for } {\cal V}=1.
$$
 Otherwise 
$$ 
  y(x) =  - {1\over 4} \left[{e^{i\pi \nu_1} \over
  16^{\nu_2-1} }  \right] x^{\nu_2}(1+O(x^{\delta}))~~\hbox{ for } 
0<{\cal V}<1,
$$
or 
$$ 
 y(x)= - {1\over 4} \left[{e^{i\pi \nu_1} \over
  16^{\nu_2-1} }  \right]^{-1} x^{2-\nu_2} (1+O(x^{\delta}))~~\hbox{ for } 
1<{\cal V}<2.
$$ 
 This is the behavior of 
Theorem 1 with $1-\sigma=\nu_2$
or $2-\nu_2$.

\vskip 0.2 cm 
\noindent
{\bf Important Observation:} 
Let $\nu_2=1-\sigma$ and consider the intersection ${\cal
D}(r;\nu_1,\nu_2)\cap \cal D(\sigma)$ in the $(\ln|x|,\Im \nu_2
\arg(x))$-plane. It is never empty.  See figure \ref{figure12}. We 
 choose $\nu_1$ such that
 $a= - {1\over 4} \left[{e^{i\pi \nu_1} \over
  16^{\nu_2-1} }  \right]$.  According to 
the Proposition in  section \ref{Parametrization of a branch through
Monodromy Data -- Theorem 2},  the 
transcendent of the elliptic representation and $y(x;\sigma,a)$ of 
Theorem 1 coincide on the intersection.
 Equivalently, we can choose the identification $1-\sigma=2-\nu_2$ and repeat
the argument.

\vskip 0.2 cm
b) If ${\cal V}=0$  the term 
$$
   {1\over \sin^2\left( -i{\nu_2\over 2} \ln x +\left[i {\nu_2\over 2} \ln 16
   +{\pi 
   \nu_1\over 2}\right] - i{\nu_2 \over 2} g(x) + {\pi v(x)\over 2\omega_1(x)} 
\right)} 
$$
is {\it oscillatory  as $x\to 0$} 
and does not vanish. Note that there are no poles because 
the denominator does not vanish in ${\cal D}(r;\nu_1,\nu_2)$ since $\left| 
{e^{i\pi \nu_1} \over 16^{\nu_2}} x^{\nu_2}\right|<r <1 $. Then 
$$ 
y(x)= O(x) + {1\over F(x)^2} 
   {1\over \sin^2\left( -i{\nu_2\over 2} \ln x +\left[i {\nu_2\over 2} \ln 16
   +{\pi 
   \nu_1\over 2}\right] - i{\nu_2 \over 2} g(x) + { v(x)\over F(x)} 
\right)} 
$$
$$ 
  = {1+O(x)\over \sin^2\left( -i{\nu_2\over 2} \ln x +\left[i {\nu_2\over 2} \ln 16
   +{\pi 
   \nu_1\over 2}\right] + \sum_{m=1}^{\infty} c_{0m}(\nu_2) \left[{e^{i\pi \nu_1}\over 16^{\nu_2}}x^{\nu_2}\right]^m 
\right)}+O(x) 
$$
The last step is obtained taking into account the non vanishing term 
in (\ref{vdix}) and  $ {\pi v(x)\over 2\omega_1(x)}= {v(x) \over F(x)} =
v(x)(1+O(x))$.  

\vskip 0.2 cm 
c)   If ${\cal V}=2$  the  series   
$$-\sum_{n=1}^{\infty}  {n\over 1- \left({e^{g(x)}\over
  16}\right)^{2n} x^{2n}}~
e^{n(2-\nu_2)g(x)} \left[{e^{-i\pi \nu_1}\over
  16^{2-\nu_2}}x^{2-\nu_2}\right]^n e^{-i\pi n {v(x)\over \omega_1(x)} } 
$$
which appears in $y(x)$ is oscillating. Simplifying we obtain:
$$
  y(x)= 
  O(x) -4(1 +O(x))~ \sum_{n=1}^{\infty} n \left[ 
{e^{-i\pi \nu_1} \over 16^{2-\nu_2}} x^{2-\nu_2}
\right]^n ~e^{-i\pi n {v(x) \over \omega_1(x)}}
$$
$$
=  {1 +O(x) \over  \sin^2\left( i{2-\nu_2\over 2} \ln x +
\left[i {\nu_2-2\over 2} \ln 16
   +{\pi 
   \nu_1\over 2}\right] + \sum_{m=1}^{\infty} 
b_{0m}(\nu_2) \left[{e^{-i\pi \nu_1}\over 16^{2-\nu_2}}x^{\nu_2}\right]^m 
\right)}+
O(x)
$$

\vskip 0.2 cm 

The observation at the end of point a) makes it
possible to 
investigate the behavior of the transcendents of Theorem 1 along a 
path (\ref{spirale}) with $\Sigma=1$. The path (\ref{spirale}) coincides with 
 (\ref{che il cielo ce la mandi buona}) for
${\cal V}=0$ if we define $1-\sigma:=\nu_2$,  for ${\cal V}=2$ if we define
$1-\sigma=2-\nu_2$. 

 In particular, we can analyze the radial convergence when $\Re
\sigma=1$. We identify $\nu_2=1-\sigma$ and choose $\nu_2=i\nu$, $\nu \neq 0$
real. Namely, $\sigma=1-i\nu$. 
 Let $x\to 0$ in ${\cal D}( r; \nu_1,i\nu)$ along the line $\arg(x)$=
constant (it is the line with ${\cal V}=0$). We have
$$
   y(x)=   {1\over F(x)^2} {1\over \sin^2\left({\nu\over 2} \ln x - \nu
   \ln 16 +{\pi \over 2} \nu_1 +{\nu\over 2} g(x) +{\pi v(x) \over 2
   \omega_1(x)} \right)} +O(x)
$$
$$
={ 1+O(x)\over \sin^2\left({\nu\over 2} \ln x - \nu \ln 16 +{\pi
\nu_1\over 2} + \sum_{m=1}^{\infty} c_{0m}(\nu) \left[\left({e^{i\pi \nu_1}
\over 16^{i\nu}}\right)x^{i\nu}\right]^m+O(x)\right)}+O(x)
$$
$$
= 
{ 1+O(x)\over \sin^2\left({\nu\over 2} \ln x - \nu \ln 16 +{\pi
\nu_1\over 2} + \sum_{m=1}^{\infty} c_{0m}(\nu) \left[\left({e^{i\pi \nu_1}
\over 16^{i\nu}}\right)x^{i\nu}\right]^m\right)}
$$
The last step is possible because $\sin(f(x)+O(x))=\sin(f(x))+O(x)=
\sin(f(x))\bigl(1+O(x)\bigr)$ if $f(x)\not  \to 0$ as $x\to 0$; this is our
case for $f(x)={\nu\over 2} \ln x - \nu \ln 16 +{\pi
\nu_1\over 2} + \sum_{m=1}^{\infty} c_{om}(\nu) \left[\left({e^{i\pi \nu_1}
\over 16^{i\nu}}\right)x^{i\nu}\right]^m $ in  ${\cal D}$. 

 \vskip 0.2 cm 
We observe that for 
  $\Re \sigma=1$ we  have a 
limitation on $\arg (x)$ in   ${\cal D}( r; \nu_1,i\nu)$,  
namely 
\be 
  - \pi \Im \nu_1 - \ln r < \nu \arg (x)
\label{doppiaonda}
\ee
This is the analogous of the limitation imposed by
$B(\sigma,a;\theta_2, \tilde{\sigma})$  of (\ref{come se non bastasse!!}).

\vskip 0.2 cm
\noindent
{\it Remark:}  If $\Im \nu_2 \neq 0$, the freedom $\nu_2 \mapsto \nu_2 +2N$, 
$N\in {\bf Z}$, is the analogous of the freedom $\sigma \mapsto 
\pm \sigma +2n$. Moreover, Theorem 3 yields different critical behaviors for 
the same transcendent on 
the different domains corresponding to $\nu_2 +2N$.

\vskip 0.2 cm
 As a last remark we observe that the coefficients in the expansion of $v(x)$
 can be computed by direct substitution of $v$ into the elliptic form of
 $PVI_{\mu}$, the right hand-side being expanded in Fourier series.


\subsection{Shimomura's Representation}

In \cite{Sh} and \cite{IKSY} S. Shimomura proved the following statement for 
 the Painlev\'e VI 
equation with any 
value of 
the parameters $\alpha, \beta,\gamma,\delta$. 
\vskip 0.2 cm
{
\it For any complex number $k$ and for any $ \sigma \not \in   (-\infty,0]
\cup[1,+\infty)$  there is a 
sufficiently small $r$  such that 
the Painlev\'e VI equation for given $\alpha,\beta,\gamma,\delta$  has a 
holomorphic solution in the domain 
$$
  {\cal D}_s(r;\sigma,k)= \{{x} \in \tilde{C_0} ~|~ |x|<r, 
~|e^{-k}{x}^{1-\sigma}|<r,~|e^k {x}^{\sigma}|<r \}$$ 
with the following representation:
$$
   y({x};\sigma,k)= {1 \over \cosh^2({\sigma-1\over 2}\ln {x}
    +{k\over 2} +{v({x})\over 2})},
$$
where 
$$ 
  v({x})= \sum_{n\geq 1} a_n(\sigma) {x}^n+ \sum_{n\geq 0, 
~m\geq 1} b_{nm}(\sigma) {x}^n (e^{-k}{x}^{1-\sigma})^m 
+ \sum_{n\geq 0,~m\geq 1} c_{nm}(\sigma) {x}^n 
(e^{k}{x}^{\sigma})^m,
$$
 $a_n(\sigma), b_{nm}(\sigma), c_{nm}(\sigma)$ are rational functions of 
$\sigma$  and the series defining $v({x})$ is convergent 
(and holomorphic) in ${\cal D}(r;\sigma,k)$. Moreover, there exists a constant 
$M=M(\sigma)$ such that
\be
  |v(x)|\leq M(\sigma) ~\left(|x|+|e^{-k}{x}^{1-\sigma}|+
                     |e^{k}{x}^{\sigma}| \right).
\label{MMM}
\ee
}
\vskip 0.2 cm

\begin{figure}
\epsfxsize=15cm
\centerline{\epsffile{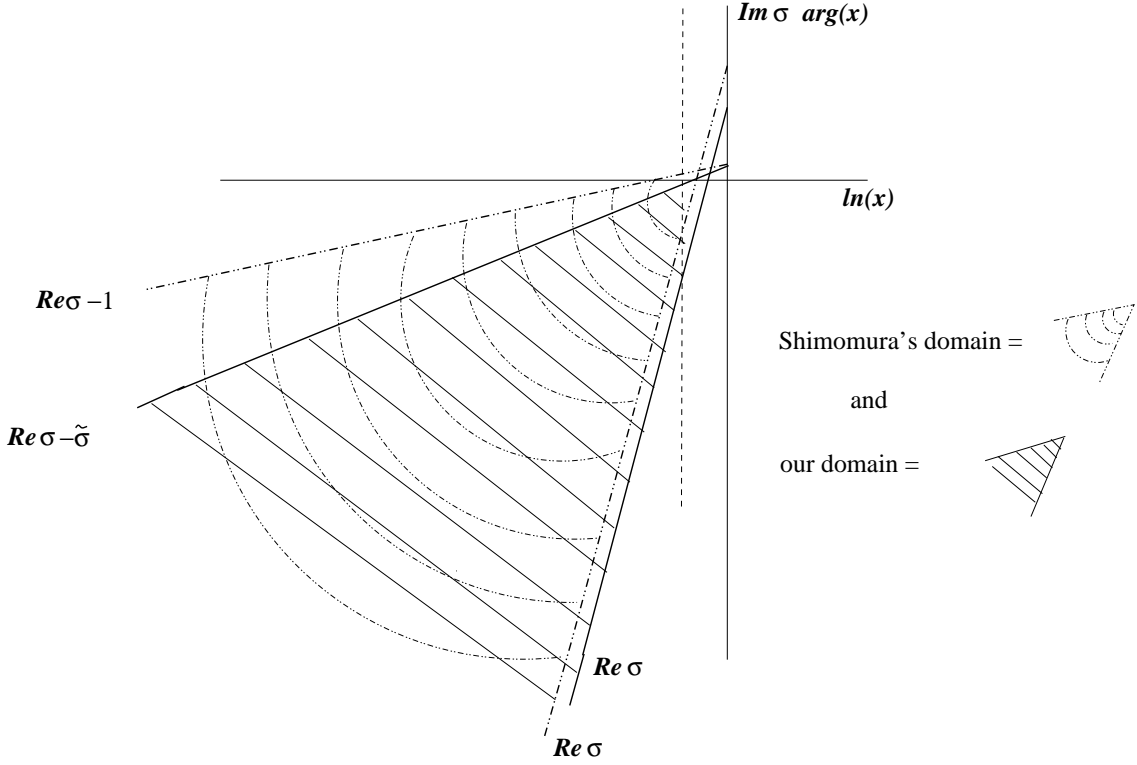}}
\caption{The domains  ${\cal
D}_s(r;\sigma,k)$ and $D(\epsilon;\sigma;\theta_1,\theta_2,\tilde{\sigma})$}
\label{figure10}
\end{figure}

The domain ${\cal D}(r;\sigma,k)$ is specified by the conditions:
\be
 |x|<r,~~~ \Re\sigma \ln|x| + [\Re k - \ln r] <
                                                    \Im \sigma \arg(x)
           < (\Re \sigma -1 )\ln |x| + [\Re k + \ln r].
\label{UFFAUFFA}
\ee
 This is an open domain in the plane $(\ln|x|,\arg(x))$. It can  be compared 
with 
the domain $D(\epsilon;\sigma,\theta_1,\theta_2)$ of Theorem 1 (figure
\ref{figure10}). 
Note that (\ref{UFFAUFFA})
imposes a  limitation on $\arg(x)$. 
 For example, if $\Re \sigma=1$ we have 
$$ 
  \Im \sigma \arg(x)< [\Re k + \ln r],~~~~~(\ln r<0)
$$
 This is similar to (\ref{doppiaonda}). 
 We will show  that  Shimomura's transcendents coincide  with those of
 Theorem 1 (see point a.1) below). So,  the above limitation turns out to be 
 the analogous of
 the limitation imposed to $D(\epsilon;\sigma;\theta_1,\theta_2)$  by
$B(\sigma,a;\theta_2, \tilde{\sigma})$  of (\ref{come se non bastasse!!}).

\vskip 0.2 cm 

 Like the elliptic representation,  Shimomura's 
 allows us to investigate what 
happens when ${x}\to 0$ along a path (\ref{spirale}) with $\Sigma=1$, 
contained in ${\cal D}_s(r;\sigma,k)$. It is a  radial path if $\Re
\sigma=1$.    Along  (\ref{spirale}) 
we have $|x^{\sigma}|=  |x|^{\Sigma} e^{-b}$.  
We suppose $\Im \sigma\neq 0$. 
\vskip 0.2 cm 

  a) $0\leq \Sigma <1$. 
 We observe that $|x^{1-\sigma} e^{-k}|\to 0$ as $x\to 0$ 
along the line.  Then:
$$
    y({x};\sigma,k)= {1 \over \cosh^2({\sigma-1\over 2}\ln x
    +{k\over 2} +{v({x})\over 2})}= {4 \over 
    x^{\sigma-1}e^ke^{v(x)}+ x^{1-\sigma}e^{-k} e^{-v(x)} +2}
$$
$$
= 4 e^{-k} e^{-v(x)} x^{1-\sigma} ~{1\over 
(1+e^{-k} e^{-v(x)} x^{1-\sigma})^2}= 4 e^{-k} e^{-v(x)} x^{1-\sigma} ~
\left( 1+e^{-v(x)} O(|e^{-k}x^{1-\sigma}|)\right) .
$$

\vskip 0.2 cm
 Two sub-cases:

\vskip 0.2 cm
a.1) $\Sigma \neq 0$. Then $|x^{\sigma} e^k| \to 0$ and $v(x) \to 0$ 
(see (\ref{MMM})). Thus
$$
     y({x};\sigma,k)= 4e^{-k} x^{1-\sigma} \left(1 + O(|x| + 
| e^k x^{\sigma}| + |e^{-k}x^{1-\sigma}|)\right)
$$
By the  Proposition in section  \ref{Parametrization 
of a branch through Monodromy Data --
Theorem 2},  
$y(x;\sigma,k)$ and $y(x;\sigma,a)$ coincide, for $a=4 e^{-k}$,  
 in  $D_s(r;\sigma,k)\cap D(\epsilon;\sigma;\theta_1,\theta_2)$. The 
intersection is not
empty for any $\theta_1$, $\theta_2$.  See figure \ref{figure10}.

\begin{figure}
\epsfxsize=15cm
\centerline{\epsffile{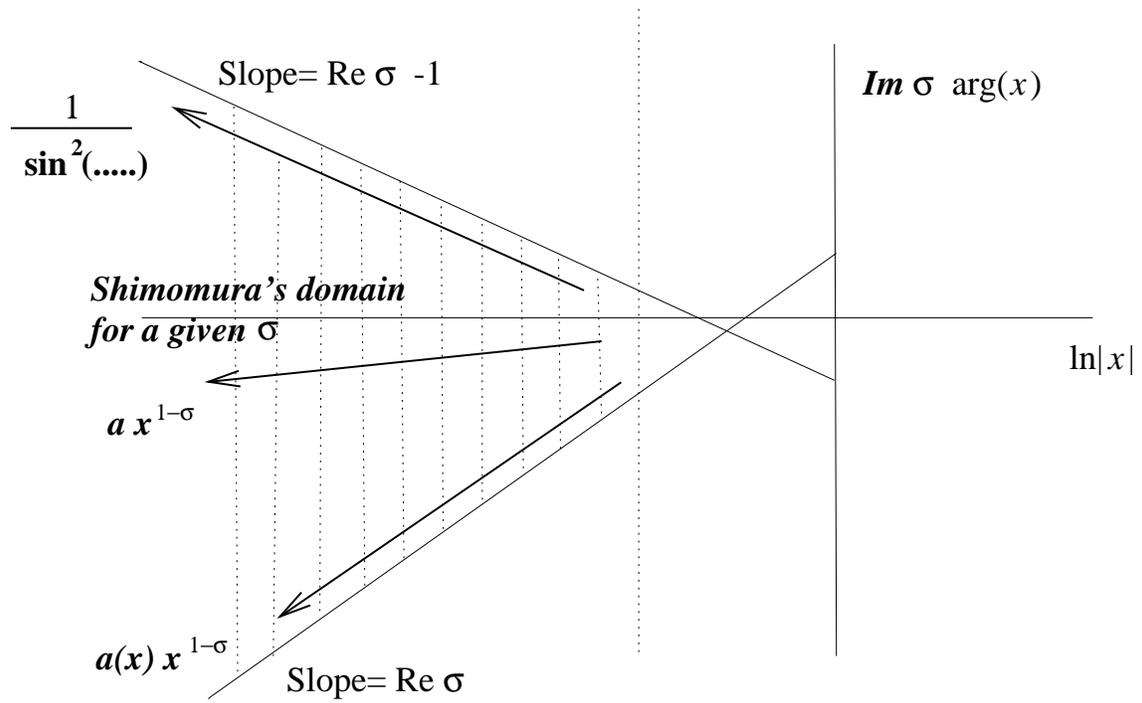}}
\caption{Critical behavior of $y(x;\sigma,k)$ along different lines in ${\cal
D}_s(r;\sigma,k)$}
\label{figure11}
\end{figure}

\vskip 0.2 cm

a.2) $\Sigma =0$.  $|x^{\sigma} e^k| \to $ constant$<r$, so $|v(x)|$ does not
vanish.  Then 
$$
y(x)=
 a(x) x^{1-\sigma}~\left( 1+O(|e^{-k}x^{1-\sigma}|)\right),~~~~
  a(x)= 4 e^{-k} e^{-v(x)},
$$ 
and  $a(x) $ must coincide with (\ref{cecilia})  of Theorem 1:
 
\vskip 0.2 cm

b) $\Sigma =1$. In this case Theorem 1 fails. 
 Now   
 $|x^{1-\sigma} e^{-k}| \to $ (constant$\neq 0)<r$. Therefore $y(x)$ does not
 vanish as $x\to 0$. We keep the representation 
$$
    y({x};\sigma,k)= {1 \over \cosh^2({\sigma-1\over 2}\ln x
    +{k\over 2} +{v({x})\over 2})} \equiv 
{1 \over \sin^2(i{\sigma-1\over 2}\ln x
    +i{k\over 2} +i{v({x})\over 2}-{\pi \over 2})}
$$
 $v(x)$ does not vanish and $y({x})$ is 
oscillating 
as ${x}\to 0$, with no limit. 
We remark that like in the elliptic representation, $\cosh^2(...)$ does not
vanish in ${\cal D}_s(r;\sigma,k)$, so we do not have poles. Figure
\ref{figure11} synthesizes points a.1), a.2), b). 

\vskip 0.3 cm
\noindent
 As an application, we consider the case $\Re \sigma =1$, namely $
  \sigma = 1 -i \nu,$ $ \nu \in \bf{R}\backslash\{0\}$. 
Then, the path corresponding to $\Sigma =1$ is a {\it  radial}  path 
 in the $x$-plane  and 
$$ 
y({x}; 1 -i \nu,k) 
= 
 {1+O(x) \over \sin^2 \left({\nu \over 2} \ln(x) + 
{i k \over 2} - {\pi \over 2}+{i\over 2} \sum_{ 
~m\geq 1} b_{0m}(\sigma) (e^{-k} x^{1-\sigma})^m\right)} 
$$



\section{Analytic Continuation of a Branch}\label{Analytic Continuation of a
Branch}

\begin{figure}
\epsfxsize= 14cm
\centerline{\epsffile{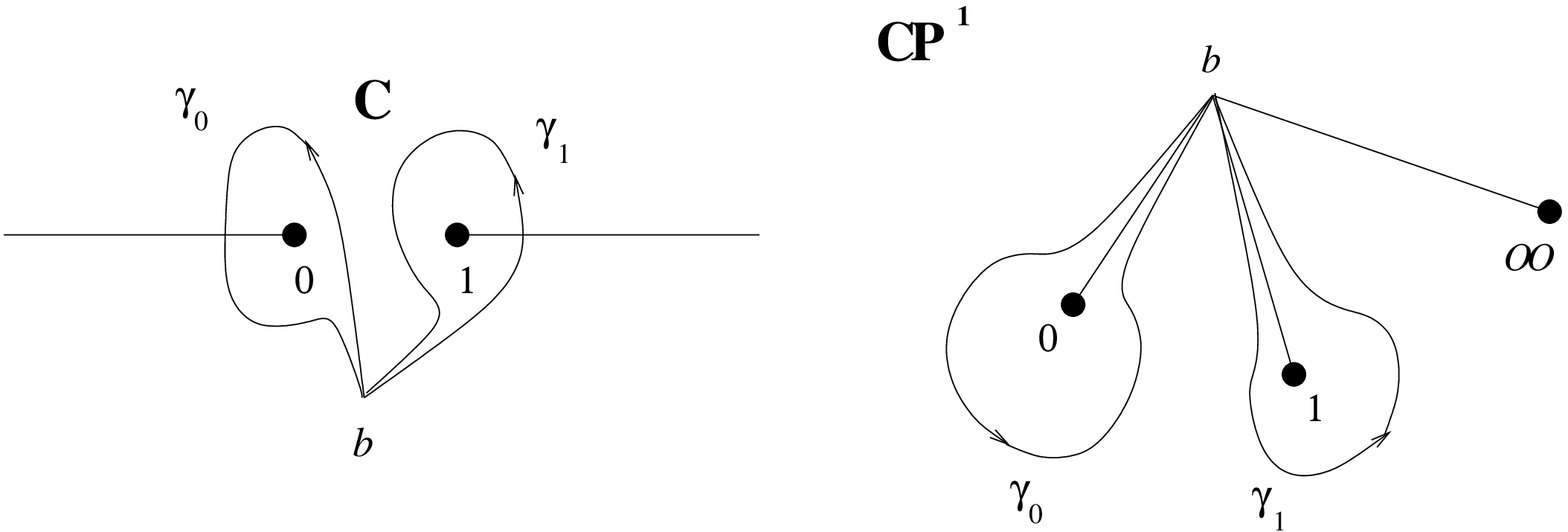}}
\caption{Base point and loops in ${\bf C}\backslash \{0,1\}$ and in ${\bf
  CP}^1\backslash \{0,1,\infty\}$. }
\label{figura2}
\end{figure}

We describe the analytic continuation of the transcendent 
  $y(x;\sigma,a)$.  
We choose
 a basis $\gamma_0$, $\gamma_1$ of two loops around 0 and 1 
respectively in the fundamental group $\pi({\bf P}^1\backslash \{0,1,\infty\}
 , b)$, where $b$ is the base-point (figure \ref{figura2}).  
   The analytic
 continuation of a branch 
 $y(x;x_0,x_1,x_{\infty})$  along paths encircling $x=0$ 
and $x=1$ (a loop around $x=\infty$ is homotopic  to the  product of
 $\gamma_0$, $\gamma_1$) is given by the action of the group of the 
pure braids on 
the monodromy data. This action is computed in \cite{DM}, to which we refer. 
 For a counter-clockwise loop around 0 we
 have 
to transform 
($x_0,x_1, x_{\infty}$)  by the action of the braid $\beta_1^2$, where 
$$
  \beta_1: ~~(x_0,x_1,x_{\infty}) \mapsto(-x_0,x_{\infty}-x_0 x_1,x_1)$$
$$ 
    \beta_1^2:~~ 
 ~~(x_0,x_1,x_{\infty}) \mapsto(x_0,~x_1+x_0x_{\infty}-x_1 x_0^2,~x_{\infty}-x_0 
x_1)
$$
The analytic continuation of the branch $y(x;x_0,x_1,x_{\infty})$ is 
the new branch $y(x; x_0,~x_1+x_0x_{\infty}-x_1 x_0^2,~x_{\infty}-x_0 
x_1)$.   
For a counter-clockwise loop around 1 we need the braid $\beta_2^2$, given by 
$$
 \beta_2: ~~(x_0,x_1,x_{\infty}) \mapsto(x_{\infty},-x_1,x_0-x_1 x_{\infty})
$$
$$
  \beta_2^2:~~(x_0,x_1,x_{\infty}) \mapsto(x_0-x_1x_{\infty},~x_1,~x_{\infty}
+x_0 x_1-x_{\infty}x_1^2)
$$
The analytic continuation of $y(x;x_0,x_1,x_{\infty})$ is 
the new branch $y(x; x_0-x_1x_{\infty},~x_1,~x_{\infty}
+x_0 x_1-x_{\infty}x_1^2)$. 
\vskip 0.2 cm
 A generic loop ${\bf P}^1\backslash \{0,1,\infty\}$ is represented by a braid
  $\beta$, which is a product of factors $\beta_1$ and $\beta_2$. 
  The braid $\beta$
 acts on $(x_0,x_1,x_{\infty})$ and
gives  a new triple $(x_0^{\beta},x_1^{\beta},x_{\infty}^{\beta})$ and a 
 new {\it branch} $y(x;x_0^{\beta},x_1^{\beta},x_{\infty}^{\beta})$.  

On the other hand,  
$y(x;x_0,x_1,x_{\infty})$ is the branch of a transcendent which 
has analytic continuation on the universal covering   of 
${\bf P}^1\backslash \{0,1,\infty\}$. We still denote  this  
transcendent by $y(x; x_0,x_1,x_{\infty})$, 
where $x$ is now regarded as a point in the universal covering. 
 A loop  transforms $x$  to a new point $x^{\prime}$ in the covering. The  
transcendent at $x^{\prime}$ is  $y(x^{\prime};x_0,x_1,x_{\infty})$.  
Let $\beta$ be the corresponding braid. We have:
\be
  y(x;x_0^{\beta},x_1^{\beta},x_{\infty}^{\beta})= 
y(x^{\prime};x_0,x_1,x_{\infty})
\label{SARAVERO?}
\ee

 Let $\sigma$, $a$ be associated to $(x_0,x_1,x_{\infty})$ according to 
 Theorem 2.  Let $x\in D(\sigma)$. At $x$ we have $y(x;x_0,x_1,x_{\infty})= 
y(x;\sigma,a)$. Let $\sigma^{\beta}$, 
$a^{\beta}=a(\sigma^{\beta}; x_0^{\beta},x_1^{\beta},x_{\infty}^{\beta})$ 
be associated to 
 $(x_0^{\beta},x_1^{\beta},x_{\infty}^{\beta})$. If 
 $ D(\sigma)\cap D(\sigma^{\beta})$ is 
not empty and $x$  also belongs to 
 $D(\sigma^{\beta})$, then $y(x;x_0^{\beta},x_1^{\beta},x_{\infty}^{\beta})
= y(x; \sigma^{\beta}, a^{\beta})$ at $x$.  
If $x\not \in D(\sigma^{\beta})$, it belongs to one and only one of the 
domains 
 $D(\pm \sigma^{\beta}+2n)$  and  $y(x;x_0^{\beta},x_1^{\beta},
x_{\infty}^{\beta})
= y(x; \pm \sigma^{\beta}+2n, \tilde{a}^{\beta})$ at $x$, where 
$\tilde{a}^{\beta} = a( \pm \sigma^{\beta}+2n;x_0^{\beta},x_1^{\beta},
x_{\infty}^{\beta})$.  We note however that if $\Re \sigma^{\beta}=1$, it may 
happen that $x$ lies in the  strip between $B(\sigma^{\beta})$ and 
$B(2-\sigma^{\beta})$,  where there may be poles (see the beginning of section 
\ref{beyond}).  In this case, we are not able to describe the analytic 
continuation (actually, the new branch may have a pole in $x$). But in 
this case, we can slightly change $\arg x$ in such a way that $x$ falls 
in a domain $D(\pm \sigma^{\beta} +2n)$.

 As an example, let us start at
$x\in D(\sigma)$; we  perform the loop $\gamma_1$ around $1$ and we go back 
to $x$. 
If $x$ also belongs to $ 
D(\sigma^{\beta_2^2})$ the transformation is
$$
   \gamma_1:~y(x;\sigma,a) \longrightarrow
   y(x;\sigma^{\beta_2^2},a^{\beta_2^2}).
$$
If $x \not \in D(\sigma^{\beta_2^2})$ but $x$ belong to one of the $D(\pm 
\sigma^{\beta_2^2} +2n)$ we have 
$$
   y(x;\sigma,a) \longrightarrow
   y(x;\pm\sigma^{\beta_2^2}+2n,\tilde{a}^{\beta_2^2}). 
$$

Again,  let us  start at
$x\in D(\sigma)$; we  perform the loop $\gamma_0$ around $0$ and we go 
back to $x$. 
The transformation of $(\sigma,a)$ according to the braid
$\beta_1$ is 
\be
(\sigma^{\beta_1^2},a^{\beta_1^2})=(\sigma,a e^{-2\pi i \sigma})
\label{analy}
\ee
as it follows from the fact that 
 $x_0$ is not affected by $\beta_1^2$, then $\sigma$ does not change, 
and from the explicit 
computation of $a(\sigma,
x_0^{\beta_1^2},x_1^{\beta_1^2},x_{\infty}^{\beta_1^2})$
 through Theorem 2  
(we will do  it at the end of  section \ref{dimth2}). Therefore, the effect 
of $\gamma_0$ is 
$$
\gamma_0:~y(x;\sigma,a) \longrightarrow
   y(x;\sigma^{\beta_1^2},a^{\beta_1^2})= y(x;\sigma,a e^{-2\pi i \sigma})
$$
 Since we are considering a loop around $0$,  it makes sense to 
consider is as a loop in ${\bf C}\backslash \{0\}\cap \{|x|<\epsilon\}$. 
The loop is $x \mapsto x^{\prime} = e^{2\pi i} x$. Suppose that also 
$x^{\prime}\in D(\sigma)$. Then, we can represent the analytic 
continuation on the universal covering as 
$$ 
  y(x;\sigma,a)\longrightarrow 
y(x^{\prime};\sigma,a)
$$
On the other hand, according to (\ref{SARAVERO?}), we must have $y(x^{\prime}; 
\sigma,a)=y(x; \sigma^{\beta_1^2},a^{\beta_1^2})$. This is immediately verified because:  
$$
y(x^{\prime};\sigma,a)
= a [x^{\prime}]^{1-\sigma}
  \left(1+O\left(\left|x^{\prime}\right|^{\delta}\right)\right) 
$$
$$
  = a e^{-2\pi i \sigma} x^{1-\sigma} (1+O(|x|^{\delta}))\equiv y(x;\sigma, a
  e^{-2\pi i \sigma} )
$$
Thus  Theorem 1 is in accordance with the
analytic continuation obtained by the action of the braid group.

\begin{figure}
\epsfxsize=15cm
\centerline{\epsffile{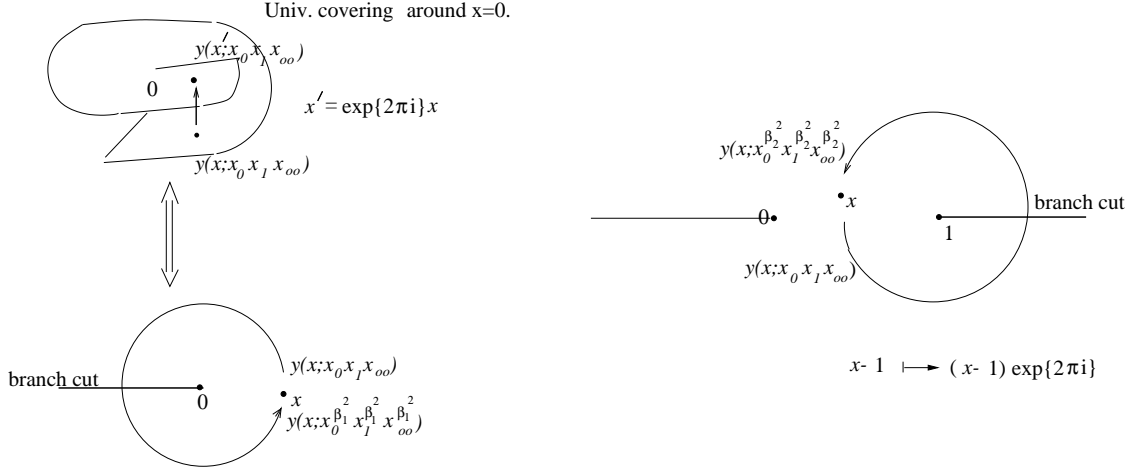}}
\caption{Analytic continuation of a branch for a loop around
$x=0$ and a loop around $x=1$. We also draw the analytic continuation on the
universal covering}
\label{figure16}
\end{figure}



\section{Singular Points $x=1$, $x=\infty$  (Connection
Problem)}\label{Singular Points x=1, x=infty  (Connection Problem)}

 In this section we restore the notation $\sigma^{(0)}$ and $a^{(0)}$
 to denote the parameters of Theorem 1 near the critical point
 $x=0$. We describe now the analogous of Theorem 1 near $x=1$ and
 $x=\infty$. 
The three critical points $0$, $1$, $\infty$ are equivalent thanks
  to the symmetries  discussed in 
\cite{Okamoto} and \cite{DM}. 

\begin{figure}
\epsfxsize=15cm
\centerline{\epsffile{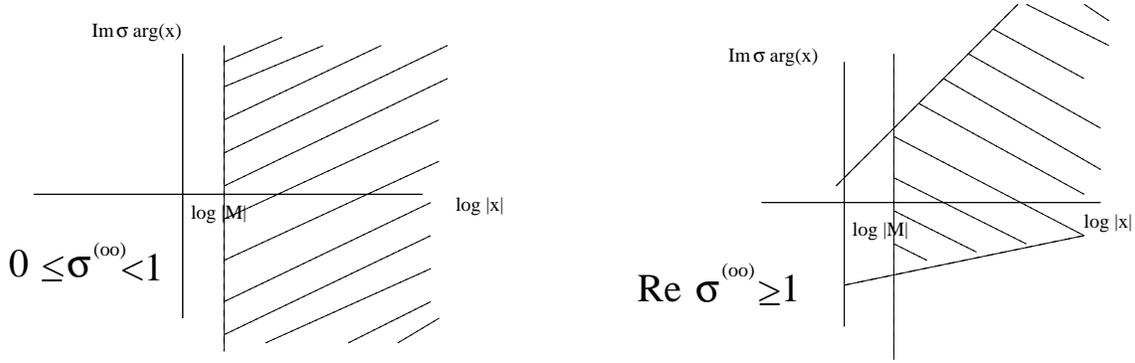}}
\caption{Some examples of the domain
$D(M;\sigma;\theta_1,\theta_2,\tilde{\sigma})$ } 
\label{figura3}
\end{figure}

\vskip 0.15 cm 
 a) Let 
\be
x={1\over t} ~~~~y(x):={1\over t}~\hat{y}(t)
\label{sim1}
\ee
 Then $y(x)$
 is a solution of $PVI_{\mu}$ (variable $x$)  if and only if
 $\hat{y}(t)$ is a solution of $PVI_{\mu}$ (variable $t$). The
 singularities $0$ and $\infty$ are exchanged. Theorem 1 holds for 
$\hat{y}(t)$ at $t=0$ with some parameters $\sigma$, $a$ that we call now 
$\sigma^{(\infty)}$, $a^{(\infty)}$. Then,  
 we go back to $y(x)$ and find a transcendent $y(x;\sigma^{(\infty)},
 a^{(\infty)})$ with behavior  
\be
    y({x}; \sigma^{(\infty)},
 a^{(\infty)}) = a^{(\infty)} {{x}}^{ \sigma^{(\infty)} }\left(
 1+O({1\over |x|^{\delta}}) \right)~~~~~{x} \to \infty
\label{asy2}
\ee
in
$$ 
    D(M; \sigma^{(\infty)};\theta_1,\theta_2,\tilde{\sigma}):=\{ 
 {x}\in 
\widetilde{ {\bf C} \backslash \{\infty\} } 
\hbox{ s.t. } 
|x|>M,~ 
e^{ -\theta_1\Im\sigma^{(\infty)} }
 |x|^{-\tilde{\sigma}}
\leq |{x}^{ -\sigma^{(\infty)}}|
\leq  e^{-\theta_2\Im\sigma^{(\infty)}}
 $$
\be
0<
   \tilde{\sigma} < 1
\}
\label{ANCORAACASAinf}
\ee
where $M>0$ is sufficiently big and $0<\delta<1$ is small (figure
\ref{figura3}).

\vskip 0.15 cm 
b) Let 
\be       
x=1-t,~~~~y(x)=1-\hat{y}(t)
\label{sim2}
\ee
$y(x)$ satisfies $PVI_{\mu}$ if and only if $\hat{y}(t)$ satisfies
$PVI_{\mu}$. Theorem 1 holds for $\hat{y}(t)$ at $t=0$ with some 
parameters $\sigma$, $a$ that we now call 
$\sigma^{(1)}$ and $a^{(1)}$.  Going back to $y(x)$ we obtain a transcendent  
 $y(x;\sigma^{(1)},a^{(1)})$ such that 
\be
y({x},\sigma^{(1)},a^{(1)})= 1-a^{(1)} (1-{x})^{1-\sigma^{(1)}}
(1+O(|1-x|^{\delta }))~~~~{x}\to 1 
\label{asy1}
\ee
in 
$$
  D(\epsilon;\sigma^{(1)};\theta_1,\theta_2,\tilde{\sigma}):= \{ {x}\in
   \widetilde{{\bf C}\backslash \{1\}} 
   \hbox{ s.t. } |1-x|<\epsilon ,~~ e^{-\theta_1 \Im \sigma}
   |1-x|^{\tilde{\sigma}} \leq |(1-{x})^{\sigma^{(1)}}| 
 \leq e^{-\theta_2 \Im
   \sigma} ,$$ 
\be
0 < \tilde{\sigma} < 1 \}
\label{ANCORAACASAuno}
\ee  

\vskip 0.2 cm 
 Consider a branch $y(x;x_0,x_1,x_{\infty})$. The symmetries 
in a) and b) affect the monodromy data, according to the following 
 formulae proved 
in \cite{DM}: 
\be 
   y(x;x_0,x_1,x_{\infty})= {1\over t} \hat{y}(t;x_{\infty},
-x_1,x_0-x_1x_{\infty}),~~~~x={1\over t} 
\label{favero1}
\ee
\be
   y(x;x_0,x_1,x_{\infty})= 1- \hat{y}(t; x_1,x_0,x_0
   x_1-x_{\infty}),~~~~x=1-t
\label{favero2}
\ee

\vskip 0.2 cm
 We are ready to solve the connection problem for the transcendents of Theorem 1, so extending the result of \cite{DM}. We recall that 
we always assume that $0\leq \Re \sigma^{(i)} \leq 1$, $i=0,1,\infty$;
 otherwise we 
write $\pm \sigma^{(i)} +2n$, $n\in {\bf Z}$. 

 We consider a transcendent 
$y(x;\sigma^{(0)}, a^{(0)})$. We choose a point $x\in D(\sigma^{(0)})$.
 At $x$  
there exists  a unique branch 
 $y(x;x_0,x_1,x_{\infty})$ whose analytic continuation in $D(\sigma^{(0)})$ is 
precisely 
 $y\bigl(x;\sigma^{(0)}, a(\sigma^{(0)})\bigr)$, where the 
triple of  monodromy
 data $(x_0,x_1,x_{\infty})$  corresponds 
 to $\sigma^{(0)}$, $a^{(0)}$ according to  Theorem 2.

 If we increase the absolute value of the point and 
we keep $\arg x$ constant, we obtain a new  point $X= |X| 
\exp\{ i \arg x\}$, where $|X|$ is big.
The branch $y(x;x_0,x_1,x_{\infty})$   is also defined in $X$, 
because we have not change 
$\arg x$.  According to (\ref{favero1}), we compute $\sigma^{(\infty)}$, 
$a^{(\infty)}$ from  the data $(x_{\infty},
-x_1,x_0-x_1x_{\infty})$ by the formulae of Theorem 2. 
Therefore, if   $X\in D(M;\sigma^{(\infty)})$, 
the analytic continuation of $y(x;x_0,
x_1, x_{\infty})=y(x;\sigma^{(0)},a^{(0)})$ at $X$ is
 $y(X;\sigma^{(\infty)},a^{(\infty)})$. 
 
\begin{figure}
\epsfxsize=15cm
\centerline{\epsffile{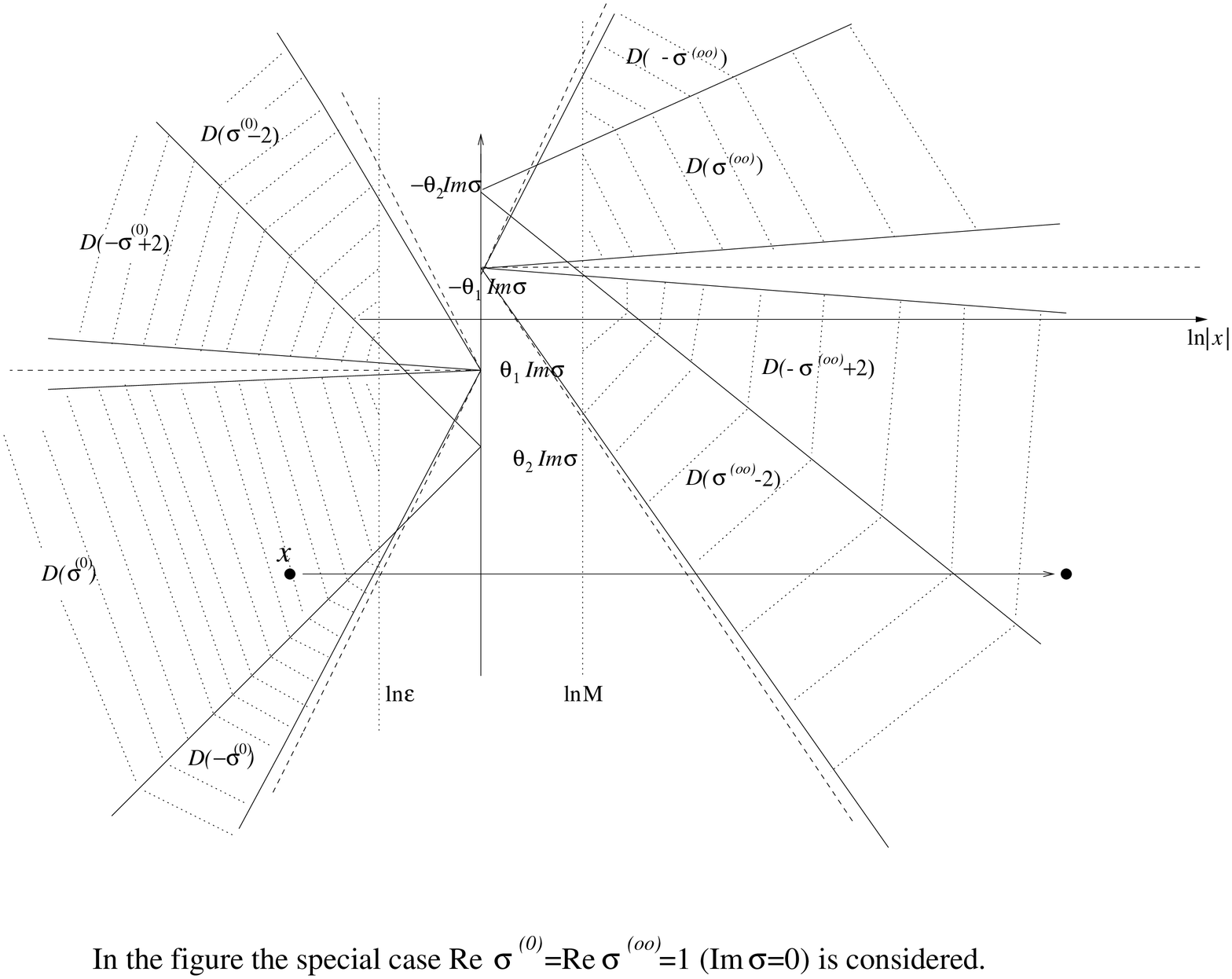}}
\caption{Connection problem for the points $x=0$, $x=\infty$}
\label{figure15}
\end{figure}

We observe that if $0\leq \Re \sigma^{(\infty)}<1$,  it is always 
possible to choose $X\in D(M;\sigma^{(\infty)})$, provided that $|X|$ is 
big enough. But for $\Re \sigma^{(\infty)}=1$   we have a restriction on the 
argument of the points of  $D(M;\sigma^{(\infty)})$ given by a set  
$B(\sigma^{(\infty)})$ analogous to 
(\ref{come se non bastasse!!}). 
This implies that $X$ may not be chosen in  $D(M;\sigma^{(\infty)})$ for any 
value of $|X|$. In this case, we can choose  $X$ in one  of the domains 
 $
D(M;\sigma^{(\infty)})$, $ D(M;-\sigma^{(\infty)})$, 
$D(M;2-\sigma^{(\infty)})$,  
$ D(M;\sigma^{(\infty)}-2)$. See figure \ref{figure15}. This is almost always 
possible, except for the case when $\arg x$ lies in the strip  
between $B(\sigma^{(\infty)})$ 
and $B(2-\sigma^{(\infty)})$, where there may be movable 
poles (see the 
discussion about these strips at the beginning of section \ref{beyond}). 

We recall that $a^{(\infty)}$ depends on $(x_{\infty},
-x_1,x_0-x_1x_{\infty})$ but it 
is also 
affected by the choice of $\pm \sigma^{(\infty)}+2n$. Thus we write below 
$a^{(\infty)}(\pm \sigma^{(\infty)}+2n)$. We conclude that 
  the analytic continuation of $y(x;x_0,
x_1, x_{\infty})=y(x;\sigma^{(0)},a^{(0)})$ at $X$ is either 
 $y\bigl(X;~\sigma^{(\infty)},~a^{(\infty)}(\sigma^{(\infty)})~\bigr)$, or 
 $y\bigl(X;~-\sigma^{(\infty)},~a^{(\infty)}(-\sigma^{(\infty)})~\bigr)$, or 
$y\bigl(X;~2-\sigma^{(\infty)},~a^{(\infty)}(2-\sigma^{(\infty)})~\bigr)$, or 
$y\bigl(X;~\sigma^{(\infty)}-2,~a^{(\infty)}(\sigma^{(\infty)}-2)~\bigr)$, 
provided that $X$ is not in the 
strip  where there may be poles. If $X$ falls in the strip, 
this is not actually a limitation, because we 
can slightly change $\arg x$ in such a way that $x$ is still in 
$D(\sigma^{(0)})$ and $X$ falls into  $
D(M;\sigma^{(\infty)}) \cup  D(M;-\sigma^{(\infty)}) \cup 
D(M;2-\sigma^{(\infty)}) \cup  D(M;\sigma^{(\infty)}-2)$.

\vskip 0.2 cm 

 In the same way we treat the connection problem between $x=0$ and $x=1$
 We repeat the same argument keeping (\ref{favero2}) into account.  
We remark again that
 for 
$\Re \sigma^{(1)}=1 $ it is necessary to consider the union of $
 D(\sigma^{(1)})$, $D(-\sigma^{(1)})$,
 $D(2-\sigma^{(1)})$, $D(\sigma^{(1)}-2)$ to include all
 possible values of $\arg(1-x)$.



\section{ Proof of Theorem 1}\label{proof of theorem 1}

We  recall that $PVI_{\mu}$ is equivalent to the  Schlesinger equations 
for the  2$\times$2 matrices $A_0(x)$,
$A_{x}(x)$, $A_{1}(x)$ of (\ref{seihat}): 
\be 
\left.\matrix{
               {d A_0\over dx}= {[A_x,A_0]\over x} \cr
\cr               
               {dA_1\over dx}={[A_1,A_x]\over 1-x} \cr
\cr  
               {d A_{x} \over d x} = {[A_x,A_0]\over x}+{[A_1,A_x]\over 1-x}\cr
}\right.
\label{sch}
\ee   
We look for solutions   satisfying  
$$A_0(x)+A_x(x)+A_{1}(x)=\pmatrix{-\mu&0\cr
                                 0 & \mu \cr}:=-A_{\infty}
~~~\mu\in{\bf C},~~~2\mu\not\in{\bf Z}
$$
$$ \hbox{ tr}(A_i)=\det(A_i)=0$$
Now let  
$$ 
A(z,x):={A_0\over z}+{A_x\over z-x}+{A_1\over z-1}
  $$
We explained that  
$y(x)$ is a solution of $PVI_{\mu}$ if and only if $A(y(x),x)_{12}=0$.

\vskip 0.15 cm
 
 The system (\ref{sch}) is a particular case of the system 
\be 
\left. \matrix{
{d A_{\mu}\over dx}=\sum_{\nu=1}^{n_2} [A_{\mu},B_{\nu}]~f_{\mu \nu}(x)
\cr\cr
{dB_{\nu}\over dx}= -{1\over x} \sum_{ \nu^{\prime}=1 }^{ n_2 }[B_{\nu},
B_{\nu^{\prime}}]+
\sum_{\mu=1}^{n_1}[B_{\nu},A_{\mu}]~g_{\mu \nu}(x) +     
\sum_{\nu^{\prime}=1}^{n_2} [B_{\nu},B_{\nu^{\prime}}] ~h_{\nu
\nu^{\prime} }(x) 
}
\right.
\label{sch1}
\ee
where the functions $f_{\mu\nu}$, $g_{\mu \nu}$, $h_{\mu \nu}$ are
meromorphic with poles at $x=1,\infty$ and $\sum_{\nu} B_{\nu} +
\sum_{\mu} A_{\mu}=-A_{\infty}$ (here the subscript $\mu$ is a label,
not the eigenvalue of $A_{\infty}$ !). System (\ref{sch}) is obtained
for $f_{\mu \nu}=g_{\mu \nu}=b_{\nu}/(a_{\mu}-xb_{\nu})$, $h_{\mu
\nu}=0$, $n_1=1$, $n_2=2$, $a_1=b_2=1$, $b_1=0$ 
and $B_1=A_0$, $B_2=A_x$, $A_1=A_1$.

\vskip 0.15 cm 
 We prove the analogous result of \cite{SMJ}, page 262, in 
 the domain $D(\epsilon;\sigma;\theta_1,\theta_2,\tilde{\sigma})$ for $\sigma 
 \not \in (-\infty,0)\cup [1,+\infty)$:

\vskip 0.3 cm
\noindent
 {\bf Lemma 1: } {\it  
                  Consider  matrices $B_{\nu}^0$ ($\nu=1,..,n_2$), 
$A_{\mu}^0$ ($\mu=1,..,n_1$) and   $\Lambda$, independent of $x$ and such
                  that 
$$
   \sum_{\nu} B_{\nu}^0 + \sum_{\mu} A_{\mu}^0 = -A_{\infty}
$$
$$
   \sum_{\nu} B_{\nu}^0=\Lambda, ~~~~\hbox{ eigenvalues}(\Lambda)=
   {\sigma \over 2},~-{\sigma\over 2},~~~
\sigma \not\in(-\infty,0)\cup [1,+\infty). 
$$
  Suppose that  $f_{\mu \nu}$, $g_{\mu \nu}$, $h_{\mu \nu}$ are holomorphic
if  $|x|<\epsilon^{\prime}$, for some small $\epsilon^{\prime}<1 $.

For any $0<\tilde{\sigma}<1$ and  $\theta_1,\theta_2$ real there exists a
sufficiently small $0<\epsilon<\epsilon^{\prime}$ such that 
the system  (\ref{sch1}) has holomorphic 
solutions 
 $A_{\mu}({x})$, $B_{\nu}({x})$   in
$D(\epsilon;\sigma;\theta_1,\theta_2,\tilde{\sigma})$   satisfying:
$$
   || A_{\mu}({x})-A_{\mu}^0||\leq C~ |x|^{1-\sigma_1},~~~~~
           || {x}^{-\Lambda} B_{\nu}({x})~ {x}^{\Lambda} 
-B_{\nu}^0||\leq C~ |x|^{1-\sigma_1}
$$          
            Here $C$ is a positive  constant and $\tilde{\sigma}<\sigma_1<1$
}

\vskip 0.3 cm
\noindent
{\bf Important remark:} There is no need 
to assume here that $2 \mu \not \in {\bf Z}$. 
The theorem holds true for any value of $\mu$.  If in the system (\ref{sch1}) 
the functions $f_{\mu \nu}$, $g_{\mu \nu}$, $h_{\mu \nu}$ are chosen in such a 
way to yield Schlesinger equations for the fuchsian system of $PVI_{\mu}$, 
the assumption $2 \mu \not \in {\bf Z}$ is still not necessary, provided the 
matrix  $R$ in (\ref{staraRIMS}) 
 is considered as a monodromy datum independent 
of the deformation parameter $x$.

\vskip 0.3 cm
\noindent
{\it Proof:} Let $A(x)$ and $B(x)$ be $2\times 2$ matrices holomorphic
on $D(\epsilon;\sigma)$ (we omit $\theta_1,\theta_2,\tilde{\sigma}$) 
  and such that 
$$ 
  ||A(x) ||\leq C_1,~~~~||B(x)||\leq C_2~~~\hbox{ on }
  D(\epsilon;\sigma)
$$
Let $f(x)$ be a holomorphic function for $|x|<\epsilon^{\prime}$. 
Let $\sigma_2$ be a real number such that $\tilde{\sigma}<\sigma_2<1$. 
Then, there exists a 
 sufficiently small $\epsilon<\epsilon^{\prime}$ such that 
for $x\in D(\epsilon;\sigma)$ we have:
$$ 
    ||x^{\pm\Lambda}~A(x)~x^{\mp\Lambda}||\leq C_1 |x|^{-\sigma_2}
$$
$$ 
    ||x^{\pm\Lambda}~B(x)~x^{\mp\Lambda}||\leq C_2 |x|^{-\sigma_2}
$$
$$
\left|\left| x^{-\Lambda}~\int_{L(x)}ds
~A(s)~s^{\Lambda}~B(s)~s^{-\Lambda}~f(s)~x^{\Lambda}  \right|    \right|
\leq C_1C_2 ~|x|^{1-\sigma_2}
$$
$$
\left|\left| x^{-\Lambda}~\int_{L(x)}ds
~s^{\Lambda}~B(s)~s^{-\Lambda}~A(s)~f(s)~x^{\Lambda}   \right|   \right|
\leq C_1C_2 ~|x|^{1-\sigma_2}
$$
where $L(x)$ is a path in $D(\epsilon;\sigma)$ joining $0$ to $x$. To
prove the estimates, we observe that  
$$
  ||x^{\Lambda}||=||x^{\hbox{diag}({\sigma\over 2},-{\sigma\over
  2})}||= \max \{ |x^{\sigma}|^{1\over 2},
  |x^{\sigma}|^{-{1\over 2}} \}\leq e^{{\theta_1\over 2}
\Im \sigma}~|x|^{-{\tilde{\sigma}\over 2}},
~~~~~~
\hbox{ in }~D(\epsilon;\sigma)
$$
Note here the importance of the bound 
$|x^{\sigma}|\leq e^{-\theta_2 \Im \sigma}$ in the
 definition of $D(\epsilon;\sigma)$: it determines the above estimates of  $
  ||x^{\Lambda}||$ because it assures that 
  $|x^{-\sigma}|^{1\over 2}$ is dominant. If this were not true, the lemma 
would fail, and Theorem 1 could not be proved.  Now we estimate  
$$
    || x^{\Lambda} ~A(x) ~x^{-\Lambda}||\leq ||
    x^{\Lambda}||~||A(x)||~ ||~x^{-\Lambda}||\leq 
             e^{\theta_1 \Im \sigma}~C_1 ~|x|^{-\tilde{\sigma}}
   $$
$$
  = \left(    e^{\theta_1 \Im \sigma}
  ~|x|^{\sigma_2-\tilde{\sigma}} \right)~C_1~|x|^{-\sigma_2}
$$
Thus, if $\epsilon$ is small enough (we require
$\epsilon^{\sigma_2-\tilde{\sigma}} \leq e^{-\theta_1 \Im \sigma}$) we
obtain $  ||x^{\Lambda}~A(x)~x^{-\Lambda}||\leq C_1 |x|^{-\sigma_2}$.

\begin{figure}
\epsfxsize=9cm
\centerline{\epsffile{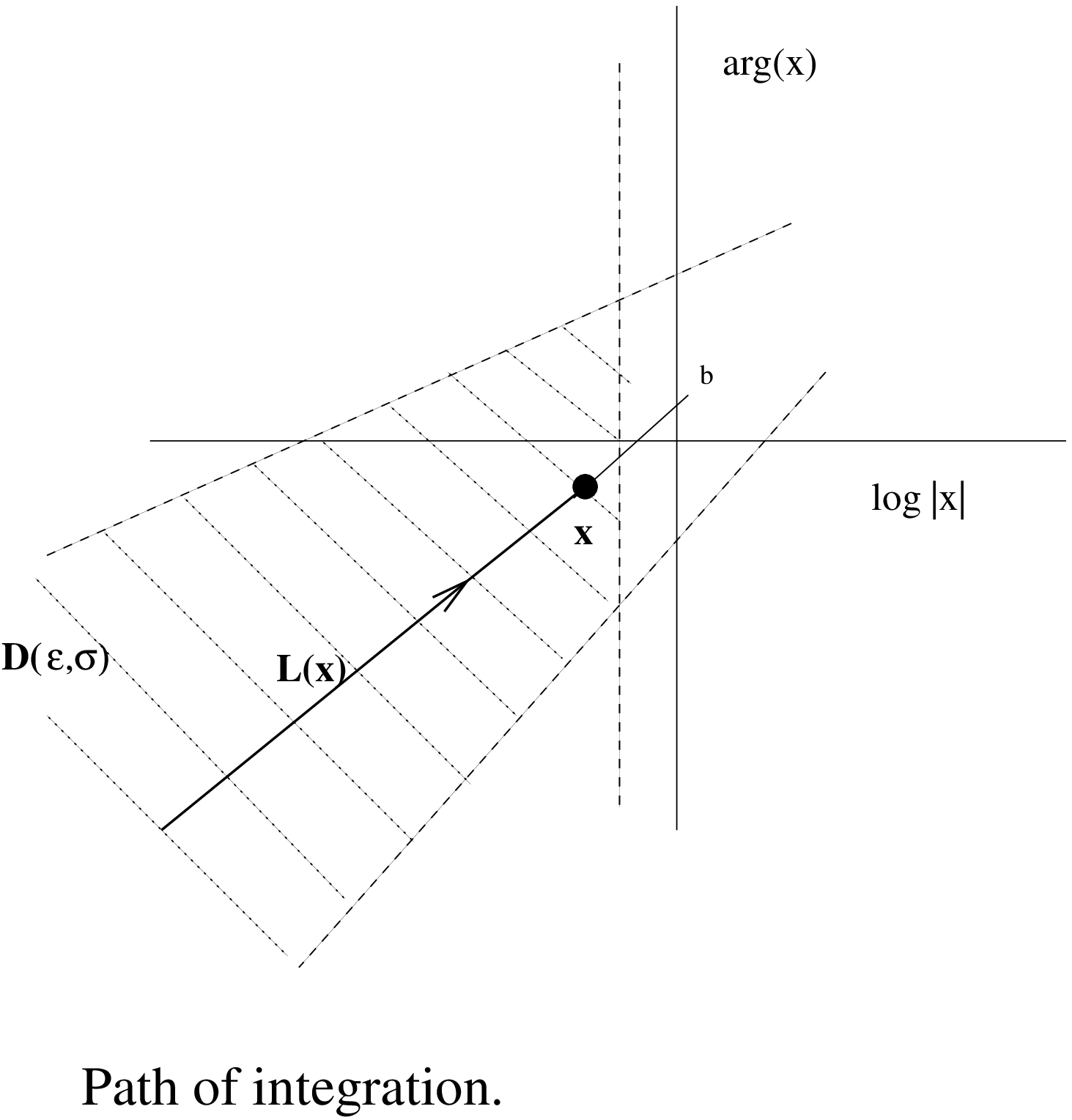}}
\caption{}
\label{figura4}
\end{figure}

We turn to the integrals.  We choose a real number $\sigma^{*}$ such that $
 0\leq \sigma^{*}\leq \tilde{\sigma}$ and 
we choose a path $L(x)$ from 0 to $x$, represented in figure \ref{figura4}. 
For $\Im \sigma \neq 0$,  $L(x)$ is given by  
$$
      \arg(s)=a \log|s| +b ,~~~~a={\Re \sigma - \sigma^{*}\over \Im
      \sigma} , ~~~~b= \arg x -{\Re \sigma - \sigma^{*}\over \Im
      \sigma} \log |x|,~~~~ |s|\leq |x|
$$
For $\Im \sigma=0$ we choose $L(x)$  with $\sigma^{*}= \sigma$ and $
 \arg(s)=\arg(x)$.  
Note that on the $L(x)$ we have
$$    |s^{\sigma}|=|x^{\sigma}| {|s|^{\sigma^*} \over |x|^{\sigma^*}}
$$
 Then
we compute 
$$
\left|\left| x^{-\Lambda}~\int_{L(x)}ds
~A(s)~s^{\Lambda}~B(s)~s^{-\Lambda}~f(s)~x^{\Lambda}  \right|    \right|
= \left|\left| \int_{L(x)} ds ~ 
   x^{-\Lambda}~A(s)~x^{\Lambda} ~\left({s\over x}\right)^{\Lambda}  
~B(s)~\left( {s\over x}   \right)^{-\Lambda} ~f(s)~                
           \right|\right|
$$
$$
\leq 
     e^{ \theta_1\Im\sigma} |x|^{-\tilde{\sigma}}
     ~C_1C_2~\max_{|x|<\epsilon} |f(x)| ~ \int_{L(x)}|ds| ~{|s|^{-\sigma^*} 
\over |x|^{-\sigma^*}}  
$$
The last step in the above inequality  follows from 
$$
  \left| \left| \left({s\over x}  \right)^{\Lambda} \right|\right|
= 
 \left| \left|
\hbox{diag}({s^{\sigma\over 2}\over x^{\sigma\over
2}},{s^{-{\sigma\over 2}} \over
x^{-{\sigma\over 2}}}) \right|\right| =  
\max_L \left\{{|s^{\sigma\over 2}|\over |x^{\sigma\over 2}|},
{|s^{-{\sigma\over 2}}|\over
|x^{-{\sigma\over 2}}|}  
\right\}  
$$
$$= 
\max \left\{  {|s|^{\sigma^*\over 2}\over |x|^{\sigma^*\over 2}},
{|s|^{-{\sigma^*\over 2}}\over |x|^{-{\sigma^*\over 2}}}  \right\}=
 {|s|^{-{ \sigma^*\over 2}}\over |x|^{-{\sigma^*\over 2}}}, 
~~~~~~|s|\leq |x|
$$
We choose the parameter $\rho =|s|$ on  $L(x)$; therefore:
$$
    s= \rho ~e^{i \left\{ \arg x +{\Re \sigma -\sigma^* \over \Im \sigma} 
          ~\log {\rho \over |x|} \right\}},~~~~0<\rho\leq|x|
$$
and we obtain 
$$
|ds|=P(\sigma,\sigma^*)~d\rho, ~~~~~
P(\sigma,\sigma^*):=\left\{ 
                    \matrix{
     \sqrt{ 1 +
\left( {\Re \sigma - \sigma^* \over \Im \sigma}\right)^2 } 
\hbox{ for } \Im \sigma \neq 0 \cr \cr
1    \hbox{ for } \Im \sigma = 0 \cr
                      } 
                         \right.
$$
$$ 
   \int_{L(x)}|ds| ~ |s|^{-\sigma^*}= P(\sigma,\sigma^*)~ \int_0^{|x|}
d\rho ~\rho^{-\sigma^*}=  {P(\sigma,\sigma^*) \over 1-\sigma^*}~ 
|x|^{1-\sigma^*}
$$
Let $P(\sigma):= \max_{\sigma^*} P(\sigma,\sigma^*)$.  
The initial integral is less or equal to 
$$
     e^{ \theta_1\Im \sigma}~\max_{|x|<\epsilon}|f(x)| 
~C_1C_2~{P(\sigma) \over 1-\tilde{\sigma}} ~|x|^{1-\tilde{\sigma}}
$$
Now, we write   
$|x|^{1-\tilde{\sigma}}=|x|^{\sigma_2-\tilde{\sigma}}~
 |x|^{1-\sigma_2}$ and we obtain, for sufficiently small $\epsilon$: 
$$
     e^{ \theta_1 \Im \sigma}~\max_{|x|<\epsilon}|f(x)| 
~C_1C_2~{P(\sigma) \over 1-\tilde{\sigma}} ~|x|^{1-\tilde{\sigma}}
\leq C_1 C_2 ~|x|^{1-\sigma_2}
$$
We remark that for $\sigma=0$ the above estimates are still valid. Actually
$
||x^{\Lambda}||\equiv||x^{\pmatrix{0&1 \cr 0 & 0\cr}}||$ diverges like
$|\log x|$, $|| x^{\Lambda} A(x) x^{-\Lambda}||$ are less or equal to 
$ C_1~|\log(x)|^2$, and finally  $|| x^{-\Lambda}~\int_{L(x)}ds
~A(s) $ $ s^{\Lambda}~B(s)~s^{-\Lambda}~f(s)~x^{\Lambda}  ||
$ is less or  equal to $C_1 C_2 \max|f| ~|\log(x)|^2 \int_{L(x)} |ds|
~|\log s|^2$.  We chose $L(x)$ to be a radial path $s=\rho \exp(i
\arg x)$, $0<\rho \leq |x|$. Then the integral is $|x|(
\log|x|^2- 2 \log|x| +2 + \alpha^2)$. The factor $|x|$ does the job,
because we rewrite it as $|x|^{\sigma_2}~|x|^{1-\sigma_2}$ (here
$\sigma_2$ is any number between 0 and 1) and we
proceed as above to choose $\epsilon$ small enough in such a way
that $\bigl(\max|f|~ |x|^{\sigma_2} \times$ function diverging like
$\log^2|x|\bigr) \leq 1$. 

The estimates above are useful  to prove the
 lemma. 

We solve the Schlesinger equations by successive
 approximations, as in \cite{SMJ}: let
 $\tilde{B}_{\nu}(x):=x^{-\Lambda}B_{\nu}(x)x^{\Lambda} $. The
 Schlesinger equations are re-written as 
\be 
 {d A_{\mu}\over dx}=\sum_{\nu=1}^{n_2}
 [A_{\mu},x^{\Lambda}\tilde{B}_{\nu} x^{-\Lambda}]~f_{\mu \nu}(x)
\label{systemdifRIMS1}
\ee
\be
{d\tilde{B}_{\nu}\over dx}= {1\over x} [\tilde{B}_{\nu},
\sum_{\mu}x^{-\Lambda}( A_{\mu}(x)-A_{\mu}^0)x^{\Lambda} ]+
\sum_{\mu=1}^{n_1}[\tilde{B}_{\nu},x^{-\Lambda}A_{\mu}x^{\Lambda} 
]~g_{\mu \nu}(x) +     
\sum_{\nu^{\prime}=1}^{n_2} [\tilde{B}_{\nu},\tilde{B}_{\nu^{\prime}}] ~h_{\nu
\nu^{\prime} }(x) 
\label{systemdifRIMS2}
\ee
We consider the following system of integral equations: 
\be
 A_{\mu}(x)= A_{\mu}^0+ \int_{L(x)}ds~\sum_{\nu}
 [A_{\mu}(s),s^{\Lambda} \tilde{B}_{\nu}(s)
 s^{-\Lambda}]~f_{\mu \nu}(s) 
\label{systemintRIMS1}
\ee
$$
\tilde{B}_{\nu}(x)= B_{\nu}^0+ \int_{L(x)}ds~ \left\{ 
         {1\over s} [\tilde{B}_{\nu}(s),\sum_{\mu}
         s^{-\Lambda} ( A_{\mu}(s)-A_{\mu}^0)s^{\Lambda}] +
         \right. 
$$
\be
\left. +
\sum_{\mu}[\tilde{B}_{\nu}(s),s^{-\Lambda}A_{\mu}(s)
s^{\Lambda} ]~g_{\mu
\nu}(s)+\sum_{\nu^{\prime}}[\tilde{B}_{\nu}(s),
 \tilde{B}_{\nu^{\prime}}(s)]~h_{\nu\nu^{\prime}}
\right\}
\label{systemintRIMS2}
\ee
We solve it  by successive approximations:
$$
 A_{\mu}^{(k)}(x)= A_{\mu}^0+ \int_{L(x)}ds~\sum_{\nu}
 [A_{\mu}^{(k-1)}(s),s^{\Lambda} \tilde{B}_{\nu}^{(k-1)}(s)
 s^{-\Lambda}]~f_{\mu \nu}(s) 
$$
$$
\tilde{B}_{\nu}^{(k)}(x)= B_{\nu}^0+ \int_{L(x)}ds~ \left\{ 
         {1\over s} [\tilde{B}_{\nu}^{(k-1)}(s),\sum_{\mu}
         s^{-\Lambda} ( A_{\mu}^{(k-1)}(s)-A_{\mu}^0)s^{\Lambda}] +
         \right. 
$$
$$
\left. +
\sum_{\mu}[\tilde{B}_{\nu}^{(k-1)}(s),s^{-\Lambda}A_{\mu}^{(k-1)}(s)
s^{\Lambda} ]~g_{\mu
\nu}(s)+\sum_{\nu^{\prime}}[\tilde{B}_{\nu}^{(k-1)}(s),
 \tilde{B}_{\nu^{\prime}}^{(k-1)}(s)]~h_{\nu\nu^{\prime}}
\right\}
$$
The functions $A_{\mu}^{(k)}(x)$, $\tilde{B}_{\nu}^{(k)}(x)$ are holomorphic in
$D(\epsilon;\sigma)$, by construction. 
Observe that $||A_{\mu}^0||\leq C$, $||B_{\nu}^0||\leq C$ for some
constant $C$.   We claim that for $|x|$ sufficiently small 
\be
\left.\matrix{
||A_{\mu}^{(k)}(x)-A_{\mu}^0||\leq C|x|^{1-\sigma_1}
\cr\cr
\left|\left| x^{-\Lambda}\left(A_{\mu}^{(k)}(x)-A_{\mu}^0 \right)x^{\Lambda} 
 \right|\right|\leq C^2 |x|^{1-\sigma_2}
\cr\cr
|| \tilde{B}_{\nu}^{(k)}(x)-B_{\nu}^0||\leq C|x|^{1-\sigma_1}
\cr
}\right.
\label{succ1}
\ee
where  $ \tilde{\sigma}<\sigma_2<\sigma_1<1$. 
Note that the above inequalities imply $||A_{\mu}^{(k)}||\leq 2C$,
$||\tilde{B}_{\nu}^{(k)}||\leq 2C$. Moreover we claim that 
\be
\left. \matrix{
||A_{\mu}^{(k)}(x)-A_{\mu}^{(k-1)}(x)||\leq C~\delta^{k-1}~|x|^{1-\sigma_1}
\cr\cr
\left|\left| x^{-\Lambda}\left(A_{\mu}^{(k)}(x)-A_{\mu}^{(k-1)}(x)
 \right )x^{\Lambda} 
 \right|\right|\leq C^2~\delta^{k-1}~ |x|^{1-\sigma_2}
\cr\cr
|| \tilde{B}_{\nu}^{(k)}(x)-\tilde{B}_{\nu}^{(k-1)}(x)||\leq
C~\delta^{k-1}|x|^{1-\sigma_1}
\cr
}
\right.
\label{succ2}
\ee
where $0<\delta<1$. 

For $k=1$ the above inequalities are proved using the simple methods
used in the estimates at the beginning of the proof. Then we proceed
by induction, still using the same estimates.
 As an example, we prove the $(k+1)^{th}$ step of the first of
(\ref{succ2}) 
supposing that the  $k^{th}$ step of (\ref{succ2}) is true.  All the other 
inequalities are proved in the same way. 
Let us consider: 
$$
||A_{\mu}^{(k+1)}(x)-A_{\mu}^{(k)}(x)||=
\left|\left| \int_{L(x)}ds~\sum_{\nu=1}^{n_2} \left( 
A_{\mu}^{(k)}s^{\Lambda}\tilde{B}_{\nu}^{(k)}s^{-\Lambda} -
A_{\mu}^{(k-1)} s^{\Lambda} \tilde{B}_{\nu}^{(k-1)} s^{-\Lambda}
+\right.\right. \right. $$
$$
\left.\left.\left. +
s^{\Lambda} \tilde{B}_{\nu}^{(k-1)} s^{-\Lambda}
A_{\mu}^{(k-1)}-s^{\Lambda} \tilde{B}_{\nu}^{(k)} s^{-\Lambda} A_{\mu}^{(k)}
\right)
         ~f_{\mu\nu}(s)
  \right|\right| \leq
$$
 $$
\leq 
 \int_{L(x)}|ds|~\sum_{\nu=1}^{n_2} \left|\left|
A_{\mu}^{(k)}s^{\Lambda}\tilde{B}_{\nu}^{(k)}s^{-\Lambda} -
A_{\mu}^{(k-1)} s^{\Lambda} \tilde{B}_{\nu}^{(k-1)} s^{-\Lambda}
\right|\right|~|f_{\mu\nu}(s)|+  $$
$$ 
+\int_{L(x)}|ds|~\sum_{\nu=1}^{n_2}
\left| \left|
s^{\Lambda} \tilde{B}_{\nu}^{(k-1)} s^{-\Lambda}
A_{\mu}^{(k-1)}-s^{\Lambda} \tilde{B}_{\nu}^{(k)} s^{-\Lambda} A_{\mu}^{(k)}
  \right|\right|~|f_{\mu\nu}(s)|
$$
Now we estimate 
$$
 \left|\left|
A_{\mu}^{(k)}s^{\Lambda}\tilde{B}_{\nu}^{(k)}s^{-\Lambda} -
A_{\mu}^{(k-1)} s^{\Lambda} \tilde{B}_{\nu}^{(k-1)} s^{-\Lambda}
\right|\right|\leq
$$
$$
\leq \left|\left|  
A_{\mu}^{(k)}s^{\Lambda}\tilde{B}_{\nu}^{(k)}s^{-\Lambda}-A_{\mu}^{(k-1)}s^{
\Lambda} \tilde{B}_{\nu}^{(k)} s^{-\Lambda} \right|\right| +
$$
$$ +\left| \left|
A_{\mu}^{(k-1)}s^{
\Lambda} \tilde{B}_{\nu}^{(k)} s^{-\Lambda} -
A_{\mu}^{(k-1)} s^{\Lambda} \tilde{B}_{\nu}^{(k-1)} s^{-\Lambda}
\right|\right|
$$
$$
 \leq || A_{\mu}^{(k)}-A_{\mu}^{(k-1)}||
 ~||s^{\Lambda}\tilde{B}_{\nu}^{(k)} s^{-\Lambda}||+
 ||A_{\mu}^{(k-1)}|| ~ || s^{\Lambda}||~||\tilde{B}_{\nu}^{(k)} 
-\tilde{B}_{\nu}^{(k-1)} || ~ ||s^{-\Lambda}||
$$
     By induction then:
$$
  \leq 
         \left( C~\delta^{k-1}~|s|^{1-\sigma_1} \right) 
~2C~e^{\theta_1\Im \sigma}|s|^{-\tilde{\sigma}} + 
      2 C~\left(C~\delta^{k-1}~|s|^{1-\sigma_1}\right) ~e^{\theta_1\Im \sigma}
|s|^{-\tilde{\sigma}}
$$ 
The other term is estimated in an analogous way. Then 
$$ 
  ||A_{\mu}^{(k+1)}-A_{\mu}^{(k)}||\leq ~{P(\sigma)\over 1-\tilde{\sigma}}~
8n_2C^2 \max |f_{\mu \nu}| ~\delta^{k-1}~e^{\theta_1\Im \sigma}
  |x|^{1-\tilde{ \sigma}}~|x|^{1-\sigma_1}
$$ 
We choose $\epsilon$  small enough to have 
${P(\sigma)\over 1-\sigma^*} ~ 8n_2C\max|f|
~e^{\theta_1\Im \sigma} |x|^{1-\tilde{\sigma}}\leq \delta$. 
Note that the choice of
$\epsilon$  is independent of
$k$.  In the case $\sigma=0$,  $|x|^{1-\tilde{\sigma}}$ is substituted by
$|x| (\log^2|x| +O(\log|x|))$. 

\vskip 0.2 cm

 The inequalities  (\ref{succ1}) (\ref{succ2}) 
 imply the convergence of the successive
 approximations to a solution of the integral equations (\ref{systemintRIMS1}),
 (\ref{systemintRIMS2})  
 satisfying the  assertion of the lemma, plus the
 additional inequality
$$ ||   x^{-\Lambda}(A_{\mu}(x)-A_{\mu}^0)x^{\Lambda}||\leq C^2
|x|^{1-\sigma_2} $$
In order to prove that the solution also solves the differential equations 
(\ref{systemdifRIMS1}), (\ref{systemdifRIMS2}) we need the following:

\vskip 0.2 cm
\noindent 
{\bf Sub-Lemma 1:}  {\it Let $f(x)$ be a holomorphic function in 
$D(\epsilon, \sigma)$ such that $f(x)= O(|x|+|x^{1-\sigma}|)$ for $x\to 0$ 
 in 
$D(\epsilon, \sigma)$. Then 
$$
  F(x):= \int_{L(x)} {1\over s} ~f(s)~ds $$
is holomorphic  in $D(\epsilon, \sigma)$ and 
${dF(x)\over dx}= {1\over x}~ f(x)$
}

\vskip 0.2 cm
We understand that the Sub-Lemma applies to our case, because the entries of 
the matrices in the integrals in (\ref{systemintRIMS1}), 
  (\ref{systemintRIMS2}) are of order $s^{-1}$, $s^{-\sigma}$ or higher. 
Thus, if we prove it, the proof of Lemma 1 will be  complete.

\vskip 0.2 cm
\noindent 
{\it Proof of Sub-Lemma 1:} Let $x+\Delta x$ be another point in 
$D(\epsilon; \sigma)$ close to $x$. To prove the Sub-Lemma it is enough to 
prove that 
$\int_{L(x+\Delta x)}{1\over s} ~f(s)~ds - \int_{L(x)} {1\over s} ~f(s)~ds
= \int_x^{x+\Delta x}    {1\over s} ~f(s)~ds$, where the last integral is 
on a segment from $x$ to $x + \Delta x$.   Namely, we prove that 
$$ 
 \left( \int_{L(x)}-\int_{L(x+\Delta x)} - \int_x^{x+\Delta x}\right)ds ~{f(s) 
\over s} =0
$$
 We consider a small disk $U_R$ centered at $x=0$ of small radius $R<
\min \{ \epsilon, |x|\}$ and the points $x_R:= L(x)\cap U_R$, $x^{\prime}_R :=
  L(x+\Delta x)\cap U_R$. Since  the integral of $f/s$ on a finite 
close curve (not containing 0)is zero we have:
\be
 \left( \int_{L(x)}-\int_{L(x+\Delta x)} - \int_x^{x+\Delta x}\right)ds ~{f(s) 
\over s}
=\left( \int_{L(x_R)}-\int_{L(x^{\prime}_R)} - \int_{\gamma(x_R,x^{\prime}_R)}
\right)ds ~{f(s) 
\over s}
\label{provcopia}
\ee
The last integral is on the arc $\gamma(x_R,x^{\prime}_R)$ 
   from $x_R$ to $x^{\prime}_{R}$ on the circle $|s|=R$. 
  We have also kept into account the obvious fact that 
$L(x_R)$ is contained in $L(x)$ and $L(x^{\prime}_R)$ is contained in $L(x+
\Delta x)$. 

We take $R\to 0$ and we prove that the r.h.s. in 
(\ref{provcopia}) vanishes. First of all we use the hypothesis, we estimate 
integrals  in the same way we did before and we obtain:
$$ 
  \left|\int_{L(x_R)} {f(s)\over s}~ds\right| \leq  \int_{L(x_R)} 
 ~{1\over |s|} O(| s| +
| s^{1-\sigma}|)~|ds|\leq {P(\sigma,\sigma^*)\over 1-\sigma^*} O(R +
O(R^{1-\sigma^*}))
$$
Therefore $ 
  \int_{L(x_R)}  {f(s)\over s}~ds \to 0$ for $R\to 0$ (recall 
that $0\leq \sigma^*<1$). In the same 
way we prove that  $ 
  \int_{L(x^{\prime}_R)}  {f(s)\over s}~ds \to 0$ for $R\to 0$. 
 We finally estimate the integral on the arc. Since $x_R\in L(x)$ and 
 $x^{\prime}_R \in L(x+\Delta x)$ we   have 
$$ 
  \arg x_R = \arg x + {\Re \sigma - \sigma^* \over \Im \sigma} 
\log {R \over |x|}, ~~~~ \arg x^{\prime}_R =
 \arg (x+\Delta x) + {\Re \sigma - \sigma^* \over \Im \sigma} 
\log {R \over |x|}.
$$
Thus $|\arg x_R - \arg x^{\prime}_R | =\left| \arg x- \arg (x+\Delta x)+ 
 {\Re \sigma - \sigma^* \over \Im \sigma} \log \left| 1+ {\Delta x \over x}
\right| \right|$ is independent of $R$. This implies that the length of 
$ \gamma(x_R,x^{\prime}_R)$ is $O(R)$. Moreover $f(x)= O(R+R^{1-\sigma^*})$ 
on the arc. Hence: 
$$
\left|\int_{\gamma(x_R,x^{\prime}_R)} {1\over s} f(s)~ds  \right| 
\leq  {1\over R} \int_{\gamma}|f(s)| |ds| = O(R^{1-\sigma^*}) \to 0 \hbox{ for } 
 R \to 0
$$
This completes the proof of Sub-Lemma 1 and Lemma 1. 

 \rightline{$\Box$}
\vskip 0.3 cm

 We observe that in the proof of Lemma 1 we imposed 
${P(\sigma)\over 1-\tilde{\sigma}} ~ 8n_2C\max|f|
~e^{\theta_1\Im \sigma} |x|^{1-\tilde{\sigma}}\leq \delta$. 
We obtain an important condition on $\epsilon$ which we used for the Remark 
in section \ref{Local Behaviour -- Theorem 1}.
\be 
  ~e^{\theta_1\Im \sigma} |\epsilon|^{1-\tilde{\sigma}}\leq c,~~~~
c:=
{\delta \over 8 n_2 C} { 1-\tilde{\sigma}\over P(\sigma)}~
 {1\over \max |f_{\mu\nu}|}
\label{condizioni su epsilon}
\ee
( here 
$C=\max\{||A_{\mu}^0||,||B_{\nu}^0||\}$ ). 

\vskip 0.3 cm  

We turn to the  case in which we are concerned: we consider three
matrices $A_0^0,A_x^0,A_1^0$ such that 
$$ 
   A_0^0+A_x^0=\Lambda,~~~A_0^0+A_x^0+A_1^0=\hbox{diag}(-\mu,\mu),~~~~
 \hbox{tr}(A_i^0)=\det(A_i^0)=0, ~~i=0,x,1
$$

\vskip 0.3 cm
\noindent
{\bf Lemma 2: } {\it Let $r$ and $s$ be two complex numbers not equal
  to 0 and
  $\infty$. 
  Let $T$ be the matrix which brings  $\Lambda$ to the  Jordan form:
$$
   T^{-1} \Lambda T = \left\{ \matrix{
                 \hbox{ diag}({\sigma \over 2}, -{\sigma \over2})
                                   ,~~~~\sigma \neq 0 \cr\cr
                         \pmatrix{0&1\cr
                                  0&0\cr},~~~~\sigma =0 \cr
}\right.
                     $$
The general solution of} 
$$A_0^0+A_x^1+A_1^0=\pmatrix{-\mu&0\cr 
                             0 & \mu \cr},~~~~\hbox{
                             tr}(A_i)=\det(A_i)=0,~~~~
A_0^0+A_x^0=\Lambda$$
{\it 
is the following:
\vskip 0.15 cm
For $\sigma\neq 0, \pm 2\mu$:
$$
\Lambda = {1\over 8 \mu} \pmatrix{
                                   -\sigma^2-(2\mu)^2 &
                                            (\sigma^2-(2\mu)^2)r \cr 
                            {(2\mu)^2-\sigma^2 \over r} &
                                   \sigma^2+(2\mu)^2 \cr }
~~~~A_1^0={\sigma^2 -(2\mu)^2\over 8\mu}\pmatrix{ 1  &  -r \cr
                                                  {1\over r} & -1 \cr}
$$
$$
  A_0^0=T \pmatrix{ {\sigma\over 4} & {\sigma \over 4} ~s \cr
                        -{\sigma\over 4}{1\over s} &  -{\sigma\over
                        4} \cr }~T^{-1},
~~~~A_x^0=T \pmatrix{ {\sigma\over 4} & -{\sigma \over 4} ~s \cr
                        {\sigma\over 4}{1\over s} &  -{\sigma\over
                        4} \cr }~T^{-1}
$$
where 
$$ T=\pmatrix{ 1 & 1 \cr
                {(\sigma + 2\mu)^2\over \sigma^2-(2\mu)^2} {1\over r}
                &
                  {(\sigma - 2\mu)^2\over \sigma^2-(2\mu)^2} {1\over
                r}\cr 
}$$
\vskip 0.15 cm 
For $\sigma=-2\mu$: $A_0^0$ and $A_x^0$ as above, but 
\be 
     \Lambda=\pmatrix{-\mu & r \cr 
                       0 & \mu \cr}~~~~~
    A_1^0=\pmatrix{0 & -r \cr
                   0 & 0 \cr }
       ~~~~~
   T = \pmatrix{1 & 1\cr 
                0 & {2 \mu \over r} \cr
           }
\label{cucu}
\ee
or 
\be
       \Lambda=\pmatrix{-\mu & 0 \cr 
                       r & \mu \cr}~~~~~
    A_1^0=\pmatrix{0 & 0 \cr
                   -r & 0 \cr }
       ~~~~~
   T = \pmatrix{1 & 0\cr 
                 -{r \over 2\mu }&1 \cr
           }
\label{cucu2}
\ee           
\vskip 0.15 cm
For $\sigma=2\mu$: $A_0^0$ and $A_x^0$ as above, but 
\be
     \Lambda=\pmatrix{-\mu & r \cr 
                       0 & \mu \cr}~~~~~
    A_1^0=\pmatrix{0 & -r \cr
                   0 & 0 \cr }
       ~~~~~
T= \pmatrix{1 & 1\cr 
                 { 2\mu\over r }&0 \cr
           }
\label{cucu1}
\ee
or 
\be
 \Lambda=\pmatrix{-\mu & 0 \cr 
                       r & \mu \cr}~~~~~
    A_1^0=\pmatrix{0 & 0 \cr
                   -r & 0 \cr }
       ~~~~~
 T = \pmatrix{0 & 1\cr 
                1 & -{r\over 2\mu} \cr
           }
\label{cucu3}
\ee
\vskip 0.15 cm 
For $\sigma=0$:
$$
   A_0^0=T \pmatrix{ 0 & s \cr 
                     0 & 0 \cr }T^{-1}
~~~~~A_x^0=T \pmatrix{0 & 1-s \cr 0 & 0 \cr }T^{-1}
$$
$$
\Lambda= \pmatrix{ -{\mu \over 2} & -{\mu^2 \over 4} r \cr
                   {1\over r}  & {\mu \over 2} \cr
}
~~~~~~A_1^0=\pmatrix{ -{\mu \over 2} & {\mu^2 \over 4} r \cr
                       -{1\over r} & {\mu \over 2} \cr
},~~
~~~~
T= \pmatrix{1 & 1 \cr
             -{2\over \mu r} & -2{\mu+2\over \mu^2} {1\over r}\cr}
$$
}

\vskip 0.2 cm
\noindent
We leave the proof as an exercise for the reader. $\Box$

\vskip 0.3 cm 

We are ready to prove Theorem 1, namely:

\vskip 0.3 cm 
\noindent 
{ \it Let $a:= -{1\over 4s}$ if $\sigma\neq 0$, or $a:=s$ if $\sigma=0$. Consider the family of  paths  
$$
\Im \sigma\arg(x) =
\Im \sigma\arg(x_0)   + 
(\Re \sigma - \Sigma) \log{|x|\over |x_0|},~~~
0\leq \Sigma\leq \tilde{\sigma},
$$ 
contained in  
 $D(\epsilon;\sigma,\theta_1,\theta_2)$, starting at $x_0$. 
If $\Im \sigma=0$ we    consider any regular path. Along these paths, 
the solutions of $PVI_{\mu}$, corresponding to
the solutions of Schlesinger  equations (\ref{sch}) 
obtained in Lemma 1, have the
following behavior for ${x}\to 0$
$$
   y(x)= a(x)~{x}^{1-\sigma} (1+O(|x|^{\delta}))
$$
where $0<\delta<1$ is a small number, and $$
a(x)=a ~~~~\hbox{ if } 0<\Sigma \leq \tilde{\sigma} ~\hbox{ or if $\sigma$ is real}
$$ 
If $\Sigma=0$, then  $x^{\sigma}= C e^{i\alpha(x)}$ ( $C$ is a constant
=$ ~|x_0^{\sigma}| \equiv |x^{\sigma}|  $ and   $\alpha(x)$ is the 
real  phase of $x^{\sigma}$) and 
\be
a(x)= a~\left(1 +{1\over 2a}~ C e^{i \alpha(x)} 
             + {1\over 16 a^2}~ 
C^2 e^{2 i \alpha(x)}  \right) = O(1). 
\label{infinitorims}
\ee
}

\vskip 0.2 cm 
\noindent
{\it Proof:}  
$y(x)$ can be computed in terms of the $A_i(x)$ from $A(y(x),x)_{12}=0$:
$$ 
  y(x)={x (A_0)_{12}\over (1+x) (A_0)_{12}+(A_x)_{12}
 +x(A_1)_{12}}\equiv {x (A_0)_{12}\over x(A_0)_{12} -(A_1)_{12} + x (A_1)_{12}}
$$ 
$$
   = -x {(A_0)_{12} \over (A_1)_{12}} ~{1\over 1-
   x(1+{(A_0)_{12}\over (A_1)_{12}}    )}
$$
As a consequence of Lemma 1 and 2 it follows that $|x~(A_1)_{12}|\leq
c ~|x|~(1+O(|x|^{1-\sigma_1}))$ and $|x~(A_0)_{12}|\leq
c ~|x|^{1-\tilde{\sigma}}~(1+O(|x|^{1-\sigma_1}))$, where $c$ is a
constant.  Then 
$$
y(x) =  -x {(A_0)_{12} \over (A_1)_{12}}~(1
+O(|x|^{1-\tilde{\sigma}}))
$$
From Lemma 2 we find, for $\sigma \neq 0, \pm 2\mu$: 
$$
(A_0)_{12}= -r {\sigma^2 - 4 \mu^2 \over 32 \mu} \left[ 
   {x^{-\sigma} \over s}
   (1+O(|x|^{1-\sigma_1}))+s~x^{\sigma}~(1+O(|x|^{1-\sigma_1})) -2
   (1+O(|x|^{1-\sigma_1}))  \right]
$$
$$ 
  (A_1)_{12}=-r {\sigma^2 - 4 \mu^2 \over 8 \mu}
  (1 +O(|x|^{1-\sigma_1}))
$$
Then  (recall that $\tilde{\sigma}<\sigma_1$)
$$
   y(x) = -{x\over 4} \left[    {x^{-\sigma} \over s}
   (1+O(|x|^{1-\sigma_1}))+s~x^{\sigma}~(1+O(|x|^{1-\sigma_1})) -2
   (1+O(|x|^{1-\sigma_1}))\right]~\bigl( 1+O(|x|^{1-\sigma_1}) \bigr)
                             $$
Now $x\to 0$ along a path 
$$
\Im \sigma \arg(x)=\Im \sigma \arg(x_0)+ (\Re \sigma -\Sigma) \log{|x|\over 
|x_0|}
$$  
for $0\leq \Sigma\leq\tilde{\sigma}$. Along this path
we rewrite  $x^{\sigma}$ in terms of its absolute value $|x^{\sigma}|=C
|x|^{\Sigma} $ ($C=|x_0^{\sigma}|/|x_0|^{\Sigma}$) and its real phase $\alpha(x)$ 
$$
  x^{\sigma}= C~|x|^{\Sigma} ~e^{i \alpha(x)},~~~~\alpha(x)=
\Re\sigma\arg(x)+\Im\sigma \ln |x|\Bigl|_{
\Im \sigma \arg(x)=\Im \sigma \arg(x_0)+ (\Re \sigma -\Sigma) 
\log{|x|\over |x_0|}}.$$
Then
$$
 y(x)  =  -{x^{1-\sigma} \over 4} \left[ 
                                     {1\over s} -2 C e^{i \alpha(x)}
                                     |x|^{\Sigma}(1+O(|x|^{1-\sigma_1})) + 
s C^2 e^{2 i
                                     \alpha(x)} ~|x|^{2\Sigma}
                                     ~(1+O(|x|^{1-\sigma_1}) )
                                                  \right]~\left( 
                     1 + O(|x|^{1-\sigma_1} )\right)
$$
For $\Sigma \neq 0$ the above expression becomes
$$ 
  y(x)= a x^{1-\sigma} ~\left( 1 + O(|x|^{1-\sigma_1}) + 
                       O(|x|^{\Sigma})
                       \right),~~~\hbox{ where } ~a:=-{1 \over 4s}
$$
We collect the two  $O(..)$ contribution in $O(|x|^{\delta})$ where
$\delta =   \min
\{ 1-\sigma_1, \Sigma \}$ is a small number between 0 and 1.  

 We take
the occasion here to remark that in the case of real $0<\sigma<1$,
 if we consider $x\to 0$ along a
radial path (i.e. $\arg(x) = \arg(x_0) $), then $\Sigma = \tilde{\sigma}
= \sigma$  and  thus:
$$   
  y(x) = \left\{ \matrix{- { 1 \over 4s} x^{1-\sigma} 
(1 +O(|x|^{\sigma})) \hbox{ for }
                                               0<\sigma <{1 \over 2}
                                               \cr
 \cr
- { 1 \over 4s} x^{1-\sigma} (1 +O(|x|^{1-\sigma_1})) \hbox{ for }
                                               {1\over 2}<\sigma <1 
                                               \cr
} \right.
 $$

Along the path with $\Sigma=0$ we have:
$$
 y(x) = -{x^{1-\sigma} \over 4}~\left({1\over s} -2 C e^{i \alpha(x)} 
             + s~ C^2 e^{2 i \alpha(x)}  \right)~\left(1 +
             O(|x|^{1-\sigma_1})  \right).  
$$
This is (\ref{infinitorims}), for $a=-{1\over 4s}$. 
We let the reader verify the theorem also in the cases $\sigma =\pm 2 \mu$
  (use the matrices (\ref{cucu}) and
(\ref{cucu1}) -- We must disregard  the matrices
(\ref{cucu2}), (\ref{cucu3}); the reason  will be
clarified in the comment following Lemma 5 and at
 the end of  the proof of Theorem 2) and in the case  $\sigma=0$. 
For $\sigma=0$ we obtain  
$$ 
  y(x) = a ~x~ \bigl(1 +O(|x|^{1-\sigma_1}) \bigr),~~~\hbox{ where } a:=s. 
$$
 
\rightline{$\Box$}

\vskip 0.2 cm 
In  the proof of Lemma 1 we imposed (\ref{condizioni su epsilon}). Hence, the  
reader may observe that $\epsilon$ depends on  $\tilde{\sigma}$, $\theta_1$
and on  $||A_0^0||$, $||A_x^0||$, $||A_1^0||$; thus it depends also 
 on  $a$.   


\section{ Proof of theorem 2}\label{dimth2}

\begin{figure}
\epsfxsize=15cm
\centerline{\epsffile{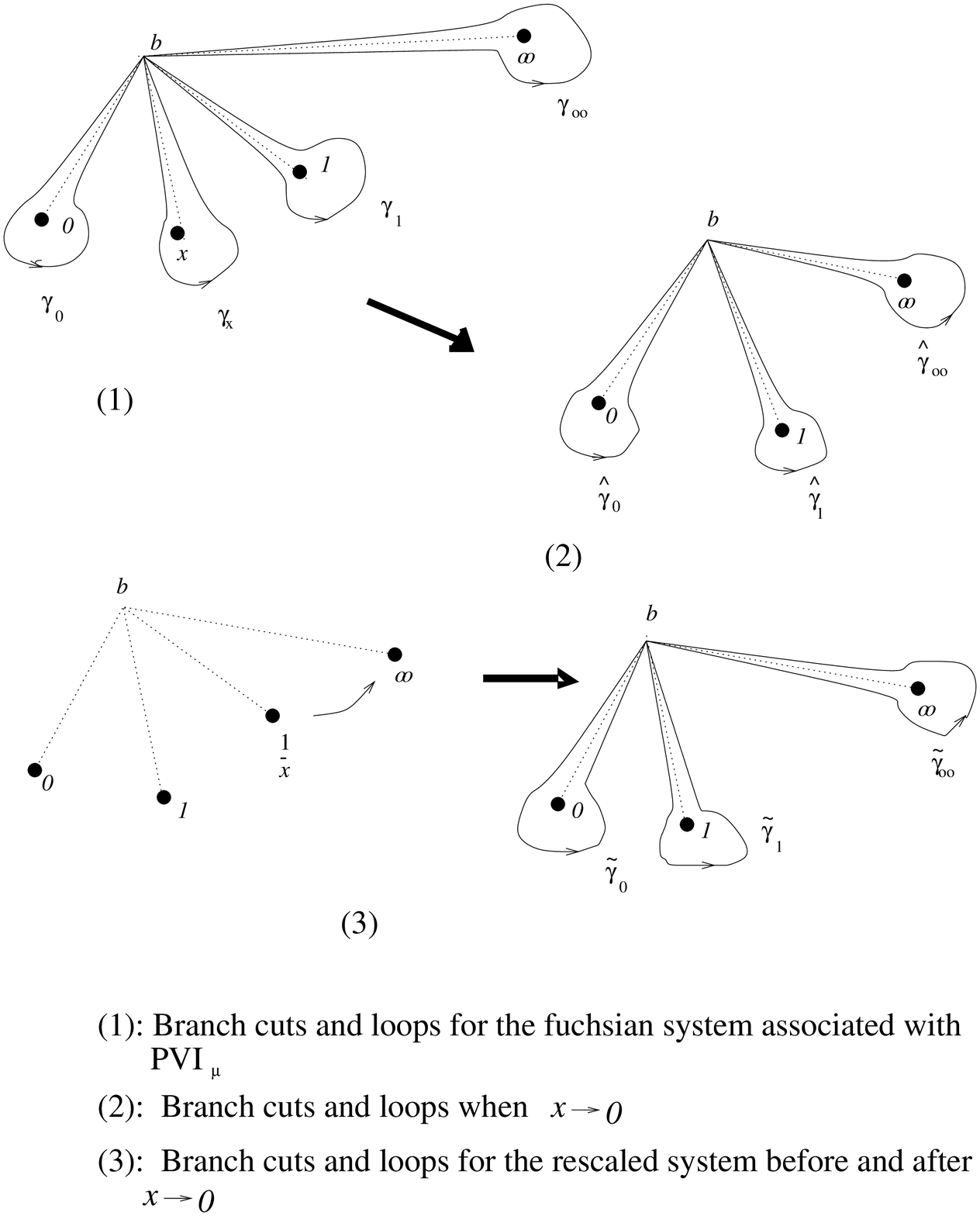}}
\caption{}
\label{figura567}
\end{figure}

We are interested in Lemma 1 when  $
  f_{\mu \nu}=g_{\mu \nu}={b_{\nu}\over a_{\mu}
  -xb_{\nu}}$, $h_{\mu\nu}=0$, $
  a_{\mu},b_{\nu}\in{\bf C}$, $a_{\mu}\neq 0~~~ \forall \mu=1,...,n_1$. 
Equations (\ref{sch1}) are the isomonodromy deformation equations for the
fuchsian system 
$$
   {dY\over dz} = \left[\sum_{\mu=1}^{n_1}{A_{\mu}(x)\over
   z-a_{\mu}}+\sum_{\nu=1}^{n_2} {B_{\nu}(x)\over z- xb_{\nu}}
   \right] ~Y
$$
As a corollary of Lemma 1, for a fundamental matrix solution $Y(z,x)$
of the fuchsian system  the limits 
$$
   \hat{Y}(z):=\lim_{{x}\to 0} Y(z,{x}),~~~~
 \tilde{Y}(z):=\lim_{{x}\to 0} ~{x}^{-\Lambda} 
Y({x}z,{x})
$$
exist when  ${x}\to 0$ in $D(\epsilon;\sigma)$. They satisfy 
$$
  {d \hat{Y}\over dz} = \left[\sum_{\mu=1}^{n_1}{A_{\mu}^0\over
   z-a_{\mu}}+{\Lambda \over z}
   \right] ~\hat{Y},~~~~
  {d\tilde{Y}\over dz} = \sum_{\nu=1}^{n_2} {B_{\nu}(x)\over z- b_{\nu}}
    ~\tilde{Y}
$$
In  our case, the last three  systems reduce to 
\be
  {dY\over dz} = \left[{A_0({x})\over z} +{A_x({x})\over z-x}+
{A_1({x})\over z-1}
   \right] ~Y
\label{stofucs}
\ee
\be
{d \hat{Y}\over dz} = \left[{A_1^0 \over z-1}+{\Lambda \over z}
   \right] ~\hat{Y}
\label{systemhat}
\ee
\be
  {d\tilde{Y}\over dz} = \left[{A_0^0\over z}+{A_x^0\over z-1}   \right]
    ~\tilde{Y}
\label{systemtilde}
\ee
Before taking the limit ${x}\to 0$, let us choose
\be
Y(z,{x})= \left(I+O\left({1\over z}\right)\right) 
~z^{-A_{\infty}}z^R ,~~~~~~z\to \infty
\label{solutionzx}
\ee 
and define as above
$$ 
   \hat{Y}(z):=\lim_{{x}\to 0} Y(z,{x}), 
~~~ \tilde{Y}(z):=\lim_{{x}\to 0} ~x^{-\Lambda} 
Y({x}z,{x})
$$
For the system (\ref{systemhat}) we choose a fundamental
matrix solution   normalized as follows
\be
   \hat{Y}_N(z)= \left(I+O\left({1\over z}\right)\right) ~z^{-A_{\infty}}z^R,
~~~~z\to \infty
\label{solutionhat}
\ee
$$
  =   (I+O(z))~z^{\Lambda}~\hat{C}_0,~~~~z\to 0 
$$
$$
   = \hat{G}_1 (I+O(z-1))~(z-1)^J ~\hat{C}_1,~~~~z \to 1 
$$
Where $\hat{G}_1^{-1} A_1^0 \hat{G}_1=J$, $J=\pmatrix{0&1\cr 0&0\cr}$ and  
$\hat{C}_0$,  $\hat{C}_1$ are invertible  connection matrices. Note that $R$
is the same of (\ref{solutionzx}), since it is independent of  $x$. For  
(\ref{systemtilde}) we choose a fundamental matrix solution normalized
as follows
 \be
   \tilde{Y}_N(z)= \left(I+O\left({1\over z}\right)\right) ~z^{\Lambda},
~~~~z\to \infty
\label{solutiontilde}
\ee
$$
  =  \tilde{G}_0 (I+O(z))~z^J~\tilde{C}_0,~~~~z\to 0 
$$
$$
   = \tilde{G}_1 (I+O(z-1))~(z-1)^J ~\tilde{C}_1,~~~~z \to 1. 
$$
Here $\tilde{G}_0^{-1}A_0^0 \tilde{G}_0 =J$, $\tilde{G}_1^{-1}A_x^0
\tilde{G}_1 =J$.
We prove that 
$$
    \hat{Y}(z) = \hat{Y}_N(z)
$$ 
\be     
 \tilde{Y}(z)=
 \tilde{Y}_N(z)~\hat{C}_0
 \label{lim}
\ee
The proof we give here uses the  technique of the proof of
 Proposition 2.1 in \cite{Jimbo}, generalized  to the domain $D(\sigma)$. 
The (isomonodromic) dependence of
 $Y(z,x)$ on $x$ is given by 
$$
  {dY(z,x)\over dx}= - {A_x(x) \over z-x} Y(z,x):=F(z,x) Y(z,x)
$$
Then 
$$ 
  Y(z,x)= \hat{Y}(z) + \int_{L(x)} dx_1~ F(z,x_1) Y(z,x_1)
$$
The integration is on a path $L(x)$ defined by  $\arg(x)=a \log|x| +b$,
$a={\Re\sigma-\sigma^* \over \Im \sigma}$ ($0\leq \sigma^*\leq
\tilde{\sigma}$),  or $\arg(x)=0$ if $\Im
\sigma=0$. The path is contained 
in $D(\sigma)$ and  joins $0$ and $x$,  like  
in the proof of Theorem 1 (figure 10). By successive approximations we have:
$$ 
   Y^{(1)}(z,x)= \hat{Y}(z) + \int_{L(x)} dx_1 F(z,x_1) \hat{Y}(z)
$$
$$
Y^{(2)}(z,x)= \hat{Y}(z) + \int_{L(x)} dx_1 F(z,x_1) Y^{(1)}(z,x_1)
$$
$$
 \vdots
$$
$$
Y^{(n)}(z,x)= \hat{Y}(z) + \int_{L(x)} dx_1 F(z,x_1) Y^{(n-1)}(z,x_1)
$$
$$
   = 
    \left[ I + \int_{L(x)} dx_1 \int_{L(x_1)} dx_2...\int_{L(x_{n-1})} dx_n
        F(z,x_1)F(z,x_2)...F(z,x_n)\right]~\hat{Y}(z)
$$
 Performing integration like in the proof of theorem 1 we evaluate 
$
  || Y^{(n)}(z,x) - Y^{(n-1)}(z,x)||$. Recall that $\hat{Y}(z)$ has 
singularities at $z=0$, $z=x$. Thus, if $|z|>|x|$ we obtain 
 $$
|| Y^{(n)}(z,x) - Y^{(n-1)}(z,x)||\leq {M C^n \over \Pi_{m=1}^n
(m-\sigma^*)} |x|^{n-\sigma^*},$$
where $M$ and $C$ are constants. 
Then $Y^{(n)}= \hat{Y}+ (Y^{(1)}-\hat{Y})+...+
(Y^{(n)}-Y^{(n-1)})$ converges for $n\to \infty$ 
 uniformly in $z$ in every compact set contained in $\{z~|~|z|>|x|\}$
and uniformly in $x\in D(\sigma)$. We can exchange limit and
integration, thus obtaining $
Y(z,x)= \lim_{n\to \infty} Y^{(n)}(z,x)$. Namely
$$
   Y(z,x)= U(z,x) \hat{Y}(z),$$
$$
   U(z,x)= I + \sum_{n=1}^{\infty} \int_{L(x)} dx_1 \int_{L(x_1)} dx_2...
         \int_{L(x_{n-1})} dx_n F(z,x_1)F(z,x_2)...F(z,x_n)
$$
being the convergence of the series uniformly in  $x\in D(\sigma)$
and in $z$ in every compact set contained in $\{z~|~|z|>|x|\}$. Of
course 
$$ 
   U(z,x)= I + O\left({1\over z}\right) 
 \hbox{ for } x\to 0 \hbox{ and } Y(z,x) \to \hat{Y}(z)
$$
But now observe that 
$$ 
  \hat{Y}(z) = U(z,x)^{-1} Y(z,x) = \left(I+O\left({1\over z}\right)
  \right) ~\left(I+O\left({1\over z}\right) \right) z^{-A_{\infty}}
  z^R, ~~~~ z\to \infty
$$
Then 
    $$ \hat{Y}(z) \equiv \hat{Y}_N(z)
$$
Finally, for $z\to 1$, 
$$
  Y(z,x) = U(x,z) \hat{Y}_N(z)= U(x,z) ~\hat{G}_1 (I +O(z-1)) (z-1)^J
  \hat{C}_1
$$
$$
  = G_1(x)(I+O(z-1))(z-1)^J \hat{C}_1$$
This implies 
             $$   C_1\equiv \hat{C}_1$$
and then
\be 
 M_1 =
\hat{C}_1^{-1} e^{2\pi i J} \hat{C}_1
\label{MONODROMY1}
\ee   
 Here we have  chosen a monodromy representation for
(\ref{stofucs}) by fixing  a base-point and a basis in the fundamental
group of 
${\bf P}^1$ as in figure \ref{figura567}. 
 $M_0$, $M_1$, $M_x$, $M_{\infty}$ are 
the monodromy matrices for the solution (\ref{solutionzx}) 
corresponding to the loops $\gamma_i$
$i=0,x,1,\infty$.  $M_{\infty}M_1M_xM_0=I$. The result
 (\ref{MONODROMY1}) may also be proved simply observing that $M_1$ becomes
$\hat{M}_1$ as $x\to 0$ in
$D(\sigma)$ because  the system (\ref{systemhat}) 
is obtained from (\ref{stofucs}) when
$z=x$ and $z=0$ merge and the singular point $z=1$ does not move. 
$x$ may converge to $0$ along spiral paths (figure \ref{figura567}). 
We recall that the braid
$\beta_{i,i+1}$ changes the monodromy matrices of ${dY\over dz} = \sum_{i=1}^n
{A_i(u)\over z-u_i} Y$ according to $M_i\mapsto M_{i+1}$, $M_{i+1}\mapsto
M_{i+1}M_i M_{i+1}^{-1}$, $M_k \mapsto M_k$ for any $k\neq i, i+1$ (see 
\cite{DM}). Therefore, 
if $\arg (x)$ increases of $2\pi$ as $x\to 0$ in (\ref{stofucs}) we have 
$$M_0 \mapsto M_x,~~~~M_x\mapsto M_x M_0 M_x^{-1},~~~M_1\mapsto M_1$$ 
If follows that $M_1$ does not change and then   
 \be 
          M_1\equiv\hat{M}_1=\hat{C}_1^{-1} e^{2\pi i J} \hat{C}_1
\label{mono1}
\ee
where $\hat{M}_1$ is the monodromy matrix of   (\ref{solutionhat}) for
the loop $\hat{\gamma}_1$ in
the basis of figure \ref{figura567}.

\vskip 0.2 cm 

Now we turn  to $\tilde{Y}(z)$. Let $\tilde{Y}(z,x):= x^{-\Lambda}
Y(xz,x)$, and by definition $\tilde{Y}(z,x) \to \tilde{Y}(z)$ as $x\to
0$. In this case 
$$ 
  {d \tilde{Y}(z,x) \over dx} = \left[{ x^{-\Lambda} (A_0+A_x)
  x^{\Lambda} - \Lambda \over x} + { x^{-\Lambda } A_1 x^{\Lambda}
  \over x - {1\over z}}\right] \tilde{Y}(z,x) := \tilde{F}(z,x)
  \tilde{Y}(z,x)
$$ 
Proceeding by successive approximations as above we get 
$$ 
   \tilde{Y}(z,x) = V(z,x) \tilde{Y}(z),$$
$$ V(z,x)=
I+\sum_{n=1}^{\infty} \int_{L(x)} dx_1...\int_{L(x_{n-1})} dx_n
\tilde{F}(z,x_1)... \tilde{F}(z,x_n) \to I \hbox{ for } x\to 0 
$$ 
uniformly in $x\in D(\sigma)$ and in $z$ in every compact subset of
$\{ z~|~ |z|<{1\over |x|} \}$. 

Let's investigate the behavior of $\tilde{Y}(z)$ as $z\to \infty$ and
compare it to the behavior of $\tilde{Y}_N(z)$. First we note that
$$ x^{-\Lambda}\hat{Y}_N(xz) = x^{-\Lambda} (I+O(xz)) (xz)^{\Lambda}
\hat{C}_0 \to z^{\Lambda} \hat{C}_0 \hbox{ for } x\to 0. 
$$
Then 
$$
    \left[ x^{-\Lambda} Y(xz,x) \right] \left[x^{-\Lambda}
    \hat{Y}_N(xz)  \right]^{-1} = x^{-\Lambda} U(xz,x) x^{\Lambda}
   \to \tilde{Y}(z) \hat{C}_0^{-1} z^{-\Lambda}.
$$ 
On the other hand, from the properties of $U(z,x)$ we know that
$x^{-\Lambda} U(xz,x) x^{\Lambda}$ is holomorphic in every compact
subset 
of $\{z~|~ |z|>1\}$ and $x^{-\Lambda} U(xz,x) x^{\Lambda}=
I+O\left({1\over z}\right)$ as $z\to \infty$. Thus 
$$
\tilde{U}(z) := \lim_{x\to 0} x^{-\Lambda} U(xz,x) x^{\Lambda}$$ 
exists uniformly in every compact subset of  $\{z~|~ |z|>1\}$ and 
$$ 
  \tilde{U}(z) = I +O\left({1\over z}  \right) ,~~~~z\to \infty
$$
Then 
$$ 
\tilde{Y}(z) = \tilde{U}(z) z^{\Lambda} \hat{C}_0 \equiv
\tilde{Y}_N(z) \hat{C}_0,
$$
as we wanted to prove. Finally, the above result implies 
$$
    Y(z,x) = x^{\Lambda} V({z\over x},x) \tilde{Y}_N({z\over x} )
    \hat{C}_0 
$$
$$
= \left\{ \matrix{
                         x^{\Lambda} V({z\over x} ,x) \tilde{G}_0 \left( I
                         + O(z/x)\right) x^{-J} z^J \tilde{C}_0
                         \hat{C}_0 ~=~ G_0(x) (I+O(z)) z^{J}  \tilde{C}_0
                         \hat{C}_0,~~~z\to 0
                                               \cr
                      \cr
                          x^{\Lambda} V({z\over x} ,x)
                         \tilde{G}_1\left( 
                                             O\left({z\over x}-1)\right) 
                                             \right) \left({z\over x}
                         -1\right)^J \tilde{C}_1 \hat{C}_0 
                     ~   =  ~G_x(x)(I+O(z-x)) (z-x)^J
                         \tilde{C}_1\hat{C}_0,~~~z\to x \cr 
} 
\right. 
$$ 
Let $\tilde{M}_0$, $\tilde{M}_1$ denote the monodromy matrices of
$\tilde{Y}_N(z)$ in the basis of figure 13. Then: 
\be
M_0=
 \hat{C}_0^{-1}\tilde{C}_0^{-1} e^{2\pi i J} \tilde{C}_0 \hat{C}_0=\hat{C}_0^{-1}\tilde{M}_0 \hat{C}_0
\label{mono2}
\ee
\be
M_x=
 \hat{C}_0^{-1}\tilde{C}_1^{-1} e^{2\pi i J} \tilde{C}_1 \hat{C}_0=\hat{C}_0^{-1}\tilde{M}_1 \hat{C}_0
\label{mono3}
\ee
The same result may be obtained observing that from  
\be
  {d (x^{-\Lambda}Y(xz,x))\over d z}= \left[
{x^{-\Lambda}A_0 x^{\Lambda}\over z} +{x^{-\Lambda}A_x
x^{\Lambda}\over z-1}+
{x^{-\Lambda}A_1 x^{\Lambda}\over z-{1\over x}}
\right]~x^{-\Lambda}Y(xz,x)
\label{prov}
\ee
we obtain the system (\ref{systemtilde})  as $z={1\over x}$ and $z=\infty$
merge (figure \ref{figura567}). The singularities  $z=0$, $z=1$, $z=1/x$
of (\ref{prov}) correspond to $z=0$, $z=x$, $z=1$ of
(\ref{stofucs}). The poles $z=0$ and $z=1$ of (\ref{prov}) do not move
as $x\to 0$ and  ${1\over x}$  converges to $\infty$, in general  along
spirals. At any turn of the spiral the system (\ref{prov}) has new
monodromy matrices according to the action of the braid group
$$ M_1 \mapsto M_{\infty},~~~M_{\infty}\mapsto M_{\infty} M_1
M_{\infty}^{-1}$$ 
but  
$$ M_0 \mapsto M_0,~~~~M_{x}\mapsto M_x$$
Hence, the limit 
$\tilde{Y}(z)$ still has monodromy $M_0$ and $M_x$ at $z=0,x$.   Since 
$\tilde{Y}=\tilde{Y}_N
\hat{C}_0$  we conclude that $M_0$ and $M_x$ are
(\ref{mono2}) and  (\ref{mono3}).

\vskip 0.3 cm 

In order to find the parameterization $y({x};\sigma,a)$ in terms
of $(x_0,x_1,x_{\infty})$ we have to compute the monodromy matrices
$M_0$, $M_1$, $M_{\infty}$ in terms of $\sigma$ and $a$  and 
then take
the traces of their products. In order to do this we use the formulae 
(\ref{mono1}), (\ref{mono2}),(\ref{mono3}).  In fact, the matrices
$\tilde{M}_i$ ($i=0,1$) and $\hat{M}_1$ can be computed explicitly
because a $2\times 2$ fuchsian system with three singular points 
can be reduced to the
hyper-geometric equation, whose monodromy is completely known. 

Before going on with the proof, we  recall that in the proof of 
Theorem 1 we defined   
$a=-{1\over 4s}$ (or $a=s$ for $\sigma=0$). 

\vskip 0.3 cm
\noindent 
{\bf Lemma 3: } {\it 
The Gauss hyper-geometric equation 
\be
   z(1-z)~ {d^2 y \over dz^2} +[\gamma_0-z(\alpha_0+\beta_0+1)]~{dy\over dz}
   -\alpha_0 \beta_0 ~y=0
\label{hyper1}
\ee
 is equivalent to the system 
\be
{d\Psi\over dz}= \left[{1\over z}\pmatrix{0 & 0 \cr
                                       -\alpha_0 \beta_0 & -\gamma_0 \cr}
                                +{1\over z-1}\pmatrix{0&1\cr
                                                     0 & \gamma_0
                                       -\alpha_0-\beta_0
                                                     \cr} \right]~\Psi
\label{hyper2}
\ee
where $\Psi=\pmatrix{y \cr (z-1) {dy\over dz}\cr}$. 
}

\vskip 0.3 cm 
\noindent
{\bf  Lemma 4: } { \it Let $B_0$ and $B_1$ be matrices  of eigenvalues $0,
1-\gamma$,  and  $0, \gamma-\alpha-\beta-1$ respectively, such
that }
$$
  B_0+B_1= \hbox{ diag}(-\alpha, -\beta),~~~\alpha\neq \beta$$
{\it Then
$$B_0=\pmatrix{
                         {\alpha(1+\beta-\gamma)\over \alpha-\beta} &
                         {\alpha (\gamma-\alpha-1)\over \alpha-\beta}~r
                         \cr
               {\beta(\beta+1-\gamma)\over \alpha-\beta}~{1\over r} & 
               {\beta(\gamma-\alpha-1)\over \alpha-\beta} \cr
}
 ,~~~
   B_1=\pmatrix{  
                   {\alpha(\gamma-\alpha-1)\over \alpha-\beta} &
                   -(B_0)_{12} \cr
                   -(B_0)_{21} & {\beta(\beta+1-\gamma)\over
                   \alpha-\beta} \cr
}
$$
for any  $r\neq 0$. 
}
\vskip 0.3 cm 
\noindent
We leave the proof as an exercise.       
The following lemma connects lemmas 3 and 4:
\vskip 0.3 cm

\noindent
{\bf Lemma 5: } {\it The system (\ref{hyper2}) with 
$$\alpha_0=\alpha,~~~\beta_0=\beta+1,~~~\gamma_0=\gamma,~~~\alpha\neq \beta$$ 
is gauge-equivalent to
the system 
\be
 {d X\over dz}= \left[{B_0\over z}+{B_1\over z-1}\right] ~X
\label{hyper3}
\ee
where $B_0$, $B_1$ are given in lemma 4. This means that there exists
a matrix 
$$
   G(z):=\pmatrix{1 & 0 \cr 
 {(\alpha-\beta)z+\beta+1-\gamma\over (1+\alpha-\gamma)r} & 
  z~{ \alpha-\beta \over \alpha(1+\alpha-\gamma)}~{1\over r}\cr
}
$$
such that $X(z)=G(z) ~\Psi(z)$. 
 It follows that (\ref{hyper3}) and the corresponding hyper-geometric
equation (\ref{hyper1}) have the same fuchsian singularities
0,1,$\infty$ and the same monodromy group.
}
\vskip 0.3 cm 
\noindent
{\it Proof:}  By direct computation. $\Box$

\vskip 0.3 cm
\noindent
 Note that the form of $G(z)$ ensures that if
$y_1$, $y_2$ are independent solutions of the hyper-geometric
equation, then a fundamental matrix of (\ref{hyper3}) may be chosen to
be 
$X(z)=
\pmatrix{ y_1(z) & y_2(z) \cr
           *  & *\cr}
$. We also observe that if we re-define $r_1:= r ~{\alpha(\gamma-\alpha-1)\over \alpha-\beta}$, the matrices $G(z)$, $B_0$, $B_1$ are not singular except for $\alpha=\beta$. Actually, we have 
$$
   B_0=   \pmatrix{{\alpha(\beta+1-\gamma)\over \alpha-\beta} & 
                     r_1 \cr 
{\alpha\beta(\beta+1-\gamma)(\gamma-1-\alpha)\over (\alpha-\beta)^2 }~{1\over r_1}
&
{\beta(\gamma-\alpha-1)\over\alpha-\beta} \cr},~~~
B_1= \pmatrix{ {\alpha(\gamma-\alpha-1)\over \alpha-\beta} & 
                - r_1  \cr
-(B_0)_{21}   &  {\beta(\beta+1-\gamma)\over \alpha-\beta} \cr
}
$$
\vskip 0.2 cm 
$$
   G(z)= \pmatrix{ 1  &   0  \cr  
{\alpha((\alpha-\beta)z+\beta+1-\gamma)\over \beta-\alpha}  ~{1\over r_1} & 
   - {z\over r_1} \cr
}
$$
The form of $B_0$, $B_1$ of Lemma 4 will correspond to the matrices define in 
Lemma 2 in general, while the form of $B_0$, $B_1$ above will correspond to
  (\ref{cucu}) and
(\ref{cucu1})  
of Lemma 2 (with $r_1 \mapsto r$). For this reason, 
we must disregard  the matrices
(\ref{cucu2}), (\ref{cucu3}) when we prove Theorem 1. 

\vskip 0.3 cm
Now we compute the monodromy matrices for the systems
(\ref{systemhat}), (\ref{systemtilde}) by reduction to an
hyper-geometric equation. 
We first study the case  $\sigma \not \in {\bf Z}$. 
 Let us start with (\ref{systemhat}). With the
gauge 
$$ Y^{(1)}(z):=z^{-{\sigma\over 2}}~\hat{Y}(z)$$ 
we transform (\ref{systemhat}) in
\be 
      {d Y^{(1)}\over dz}=\left[ {A_1^0\over
      z-1}+{\Lambda-{\sigma\over2 }I\over z} \right]~Y^{(1)}
\label{novosis}
\ee
We identify the matrices $B_0$, $B_1$ with    
$\Lambda-{\sigma\over2 }I$ and  $A_1^0$, with eigenvalues 0, $-\sigma$
and 0, 0 
respectively. Moreover $A_1^0 +\Lambda-{\sigma\over2 }I=$
diag($-\mu-{\sigma\over 2},\mu-{\sigma\over 2})$. Thus:
$$
  \alpha=\mu+{\sigma\over 2},~~~\beta=-\mu+{\sigma\over
  2},~~~\gamma=\sigma+1;~~~~ \alpha-\beta= 2\mu \neq 0 ~~~\hbox{ by hypothesis}$$
The parameters of the correspondent hyper-geometric equation are 
$$ 
 \left\{ \matrix{                                         
\alpha_0=\mu+{\sigma\over 2} \cr
\beta_0=1-\mu+{\sigma\over
  2} \cr
      \gamma_0=\sigma+1 \cr
}\right.
$$
 From them  we deduce the nature of two linearly independent solutions
  at $z=0$. Since $\gamma_0\not \in {\bf Z}$   ($\sigma \not \in
  {\bf Z}$) the solutions are expressed in terms of
  hyper-geometric functions. On the other hand, the effective 
  parameters at $z=1$ and $z=\infty$ are respectively:
$$ 
 \left\{ \matrix{
\alpha_1:=\alpha_0=\mu+{\sigma\over 2} \cr
\beta_1 := \beta_0=1-\mu+{\sigma\over
  2} \cr
     \gamma_1:=\alpha_0+\beta_0-\gamma_0+1 =1 \cr
}\right. ,~~~~
 \left\{ \matrix{
\alpha_{\infty}:=\alpha_0=\mu+{\sigma\over 2} \cr
\beta_{\infty}:= \alpha_0-\gamma_0+1=\mu-{\sigma\over
  2} \cr
     \gamma_{\infty}=\alpha_0-\beta_0+1 =2\mu \cr
}\right.
$$
Since $\gamma_1=1$, at least one  solution has a
logarithmic singularity at $z=1$. Also note that 
$\gamma_{\infty}=2\mu$, therefore  logarithmic singularities appear
at $z=\infty$ if $2\mu \in {\bf Z}\backslash \{0\}$. 

\vskip 0.15 cm 
For the derivations which follows,  we use the notations of the
fundamental paper by Norlund \cite{Norlund}. 
 To derive the connection formulae we use
the paper of Norlund when logarithms are involved. Otherwise, in the
generic case, any textbook of special functions (like \cite{Luke}) may
be used.

\vskip 0.2 cm
 \noindent
{\bf First case: $\alpha_0, \beta_0 \not \in {\bf Z}$}. This means 
$$ 
  \sigma \neq \pm 2\mu +2m ,~~~~m\in{\bf Z}
$$

We can choose the following independent solutions of the
hyper-geometric equation: 

\vskip 0.15 cm
\noindent
At $z=0$ 
$$ 
  y_1^{(0)}(z)=F(\alpha_0,\beta_0,\gamma_0;z)$$
\be
  y_2^{(0)}(z)=
  z^{1-\gamma_0}~F(\alpha_0-\gamma_0+1,\beta_0-\gamma_0+1,2-\gamma_0;z) 
\label{ffirr1}
\ee
where $F(\alpha,\beta,\gamma;z)$ is the well known hyper-geometric
  function (see \cite{Norlund}).

\vskip 0.15 cm
\noindent
At $z=1$
$$ 
 y_1^{(1)}(z)=F(\alpha_1,\beta_1,\gamma_1;1-z),~~~~
y_2^{(1)}(z)=g(\alpha_1,\beta_1,\gamma_1;1-z)
$$
Here $g(\alpha,\beta,\gamma;z)$ is a logarithmic solution introduced in
\cite{Norlund}, and $\gamma\equiv \gamma_1=1$. 

\vskip 0.15 cm
\noindent
At $z=\infty$, we consider first the case $2\mu \not \in {\bf Z}$, while the
resonant case will be considered later. Two independent solutions are:  
$$
 y_1^{(\infty)}=z^{-\alpha_0}
 ~F(\alpha_{\infty},\beta_{\infty},\gamma_{\infty} ;{1\over z}),~~~~
y_2^{(\infty)}=z^{-\beta_0}~F(\beta_0,\beta_0-\gamma_0+1,\beta_0-\alpha_0
+1; {1\over z})
$$
Then, from the connection formulas between $F(...;z)$ and $g(...;z)$ 
of \cite{Luke} and \cite{Norlund} we
derive 
$$
  [y_1^{(\infty)},y_2^{(\infty)}]=[y_1^{(0)},y_2^{(0)}] ~C_{0\infty}
$$
$$ 
  C_{0\infty}= \pmatrix{ 
                 e^{-i\pi \alpha_0}
                 {\Gamma(1+\alpha_0-\beta_0)\Gamma(1-\gamma_0) \over
                 \Gamma(1-\beta_0) \Gamma(1+\alpha_0-\gamma_0)}& 
                 e^{-i\pi \beta_0}{
                 \Gamma(1+\beta_0-\alpha_0)\Gamma(1-\gamma_0)\over
                 \Gamma(1-\alpha_0) \Gamma(1+\beta_0-\gamma_0)} \cr
                 e^{i\pi(\gamma_0-\alpha_0-1)}{\Gamma(1+\alpha_0-\beta_0)
\Gamma(\gamma_0 -1)\over \Gamma(\alpha_0)\Gamma(\gamma_0-\beta_0)}&
                      e^{i\pi(\gamma_0-\beta_0-1)}{\Gamma(1+\beta_0-\alpha_0)
\Gamma(\gamma_0-1)\over \Gamma(\beta_0)\Gamma(\gamma_0-\alpha_0)} \cr
}
$$

\vskip 0.2 cm
$$
[y_1^{(0)},y_2^{(0)}]=[y_1^{(1)},y_2^{(1)}] ~C_{01}
$$
$$
   C_{01}=\pmatrix{
                   0 & -{\pi \sin(\pi(\alpha_0+\beta_0))\over
                   \sin(\pi\alpha_0) \sin(\pi
                   \beta_0)}{\Gamma(2-\gamma_0)\over
                   \Gamma(1-\alpha_0) \Gamma(1-\beta_0) } \cr
                   -{ \Gamma(\gamma_0) \over \Gamma(\gamma_0
                   -\alpha_0) \Gamma(\gamma_0-\beta_0)} & -{
                   \Gamma(2-\gamma_0) \over \Gamma(1-\alpha_0)
                   \Gamma(1- \beta_0)} \cr
}
$$
We observe that 
$$Y^{(1)}(z)= \left(I+{F\over
z}+O\left({1\over z^2}\right)\right)~z^{\hbox{diag}(-\mu-{\sigma\over2},\mu-{\sigma\over 2})},~~~~z\to
\infty
$$
$$ 
 =
\hat{G}_0(I+O(z))~z^{\hbox{diag}(0,-\sigma)}~\hat{G}_0^{-1}\hat{C_0},~~~~z\to
0
$$
$$= \hat{G}_1(I+O(z-1))~(z-1)^J~\hat{C}_1,~~~~z\to 1
$$
where $\hat{G}_0\equiv T$ of lemma 2; namely $\hat{G}_0^{-1}\Lambda
\hat{G}_0=$ diag$({\sigma\over 2},-{\sigma\over 2})$. 
By direct substitution in the
differential equation we compute the coefficient $F$ 
$$
   F= -\pmatrix{ (A_1^0)_{11} & {(A_1^0)_{12}\over 1-2\mu} \cr
                  {(A_1^0)_{21}\over 1+2\mu} &
         (A_1^0)_{22}
}
,~~~\hbox{where}
~~A_1^0={\sigma^2 -(2\mu)^2\over 8\mu}\pmatrix{ 1  &  -r \cr
                                                  {1\over r} & -1 \cr}
$$
Thus, from the asymptotic behavior of the hyper-geometric function
($F(\alpha,\beta, \gamma;{1\over z})\sim 1$, $z\to \infty$ ) we
derive 
$$
  Y^{(1)}(z)= \pmatrix{ y_1^{(\infty)}(z) & r{\sigma^2-(2\mu)^1\over
  8\mu(1-2\mu) } ~y_2^{(\infty)} \cr
  * & * \cr
}
$$
From  
\be
 Y^{(1)}(z) \sim \pmatrix{ 1 & z^{-\sigma} \cr
 * & * \cr}
~\hat{G}_0^{-1} \hat{C}_0,~~~z\to 0
\label{ffirr2}
\ee
we derive 
 $$
 Y^{(1)}(z)= 
          \pmatrix{ y_1^{(0)}(z) & y_2^{(0)}(z) \cr
* & * \cr
} 
\hat{G}_0^{-1} \hat{C}_0
$$
Finally, observe that $\hat{G}_1=\pmatrix{ u & {u\over \omega} +v r \cr
                                 {u\over r} & v \cr
}
$ for arbitrary $u,v\in{\bf C}$, $u\neq 0$, and $\omega:= 
{\sigma^2-(2\mu)^2\over 8\mu}$. We recall 
 that $y_2^{(1)}=g(\alpha_1,\beta_1,1;1-z)\sim
                                 \psi(\alpha_1)+\psi(\beta_1)-2\psi(1)
                                 -i\pi +\log(z-1)$, $|\arg(1-z)|<
\pi$, as $z\to 1$. We can choose  $u=1$ and a suitable $v$, in such a way
                                 that the asymptotic behavior of
                                 $Y^{(1)}$ for $z\to 1$ is precisely
                                 realized by 
$$
  Y^{(1)}(z)= \pmatrix{ y_1^{(1)}(z) & y_2^{(1)}(z) \cr
                              * & * \cr
}~\hat{C}_1
$$
Therefore  we conclude that the connection matrices are: 
$$
  \hat{C}_0= \hat{G}_0 \pmatrix{
(C_{0\infty})_{11} & r~{\sigma^2-(2\mu)^2\over 8\mu
(1-2\mu)}(C_{0\infty})_{12}\cr
(C_{0\infty})_{21}& r~ {\sigma^2-(2\mu)^2\over 8\mu
(1-2\mu)}(C_{0\infty})_{22}\cr
}
$$
$$
\hat{C}_1= C_{01}~(\hat{G}_0^{-1}\hat{C}_0)= C_{01}   
\pmatrix{
            (C_{0\infty})_{11} & r~{\sigma^2-(2\mu)^2\over 8\mu
(1-2\mu)}(C_{0\infty})_{12} \cr
(C_{0\infty})_{21}& r~ {\sigma^2-(2\mu)^2\over 8\mu
(1-2\mu)}(C_{0\infty})_{22} \cr
}
$$

\vskip 0.2 cm
 It's now time to consider the resonant case $2\mu \in {\bf
Z}\backslash \{0\}$. The behavior of $Y^{(1)}$ at $z=\infty$ is 
$$
   Y^{(1)}(z) = \left(I +{F\over z} +O\left({1\over z^2} \right) \right) 
                   ~z^{\hbox{diag}(-\mu -{\sigma\over 2}, \mu -{\sigma\over
                   2}) } z^R$$
$$ 
   R= \pmatrix{ 0 & R_{12} \cr 
                0 & 0      \cr}, ~~~\hbox{ for } \mu ={1\over 2},1,{3\over
                2}, 2, {5\over 2},...
$$
$$ 
   R= \pmatrix{ 0 & 0 \cr 
                R_{21} & 0      \cr}, ~~~\hbox{ for } \mu =-{1\over 2},-1,
                - {3\over2}, -2, -{5\over 2},...
$$
and the entry $R_{12}$ is determined by the entries of $A^{0}_1$.
 For example,
if $\mu={1\over 2}$ we can compute $R_{12}= (A^0_1)_{12}= -r {\sigma^2-1 \over
4} $ (and $F_{12}$ arbitrary); if $\mu=-{1\over 2}$ we have $ R_{21}=
(A^0_1)_{21}=-{1\over r}
{\sigma^2 -1 \over 4}$ (and $F_{21}$ arbitrary); if $\mu =1$  we have 
$R_{12}= -r {\sigma^2(\sigma^2-4)\over 32}$.

Since $\sigma\not \in {\bf Z}$, $R\neq 0$. 
This is true for any $2\mu  \in {\bf Z}\backslash
\{0\}$. Note that the $R$ computed here coincides (by
isomonodromicity) to the $R$ of the system(\ref{stofucs}).

Therefore, there is a logarithmic solution at $\infty$. 
Only
$C_{0\infty}$ 
and thus
$\hat{C}_0$  and $\hat{C}_1$ change with respect to the non-resonant
case. 
We will see in a while that such matrices
disappear in the computation of tr($M_i M_j)$, $i,j=0,1,x$. Therefore, it is
not necessary to know them explicitly, 
 the only important matrix
to know being $C_{01}$, which is not affected by resonance of
$\mu$. This is the reason why the formulae of theorem 2 hold true also
in the resonant case.


\vskip 0.3 cm
\noindent
{\bf Second case:} $\alpha_0, \beta_0 \in {\bf Z}$, namely
$$ 
\sigma=\pm 2\mu+2m ,~~~~~~m\in {\bf Z}
$$

The formulae are almost identical to the first case, but $C_{01}$
changes. To see this, we need to distinguish four  cases.

i) 
  $\sigma=2\mu+2m$, $m=-1,-2,-3,...$. We choose
$$
   y_2^{(1)}(z)= g_0(\alpha_1,\beta_1,\gamma_1;1-z)
$$
Here $g_0(z)$ is another logarithmic solution of \cite{Norlund}. 
Thus 
$$C_{01}= \pmatrix{ {\Gamma(-m)\Gamma(-2\mu-m+1)\over \Gamma(-2\mu-2m)}
& 0 \cr
 0 & -{\Gamma(1-2\mu -2m)\over \Gamma(1-m-2\mu)\Gamma(-m)} \cr
}
$$
As usual, the matrix is computed from the connection formulas between
the hyper-geometric functions and $g_0$ that the reader can find in
\cite{Norlund}. 

ii) $\sigma=2\mu+2m$, $m=0,1,2,...$. We choose
$$
y_1^{(2)}=g(\alpha_1,\beta_1,\gamma_1;1-z)
$$
Thus
$$C_{01}= \pmatrix{0 & {\Gamma(m+1)\Gamma(2\mu+m)\over
\Gamma(2\mu+2m)}\cr
-{\Gamma(2\mu+2m+1)\over \Gamma(2\mu+m)\Gamma(m+1)} & 0 
\cr}
$$

iii)  $\sigma=-2\mu+2m$, $m=0,-1,-2,...$. We choose
$$
   y_2^{(1)}(z)=  g_0(\alpha_1,\beta_1,\gamma_1;1-z)
$$
Thus
$$C_{01}= \pmatrix{ {\Gamma(1-m)\Gamma(2\mu-m)\over \Gamma(2\mu-2m)}
& 0 \cr
 0 & -{\Gamma(1+2\mu -2m)\over \Gamma(2\mu-m)\Gamma(1-m)} \cr
}
$$

iv)  $\sigma=-2\mu+2m$, $m=1,2,3,...$. We choose
 $$
  y_2^{(1)}(z)= g(\alpha_1,\beta_1,\gamma_1;1-z)
$$
Thus
$$C_{01}= \pmatrix{0 & {\Gamma(m)\Gamma(m+1-2\mu)\over
\Gamma(2m-2\mu)}\cr
-{\Gamma(2m+1-2\mu)\over \Gamma(m+1-2\mu)\Gamma(m)} & 0 
\cr}
$$

Note that this time $F=\pmatrix{0 & {r\over 1-2\mu} \cr
 0 & 0 \cr}$ in the case $\sigma=\pm 2\mu$ (i.e. $m=0$) because
 $A_1^0$ has a special form in this case.  Then in
 $\hat{C}_0$ the elements $ {\sigma^2-(2\mu)^2\over 8\mu
(1-2\mu)}(C_{0\infty})_{12}$, $ {\sigma^2-(2\mu)^2\over 8\mu
(1-2\mu)}(C_{0\infty})_{22}$ must be substituted, for $m=0$,  with 
$ {1\over 
1-2\mu}(C_{0\infty})_{12}$, 
$ {1\over 
1-2\mu}(C_{0\infty})_{22}$.

\vskip 0.3 cm 
We turn to the system (\ref{systemtilde}). Let $\tilde{Y}$ be the
fundamental matrix (\ref{solutiontilde}). With the gauge 
$$ Y^{(2)}(z):=\hat{G}_0^{-1}~\left(\tilde{Y}_N(z)\hat{G}_0 \right)$$ 
we have 
$$
  {d Y^{(2)} \over d z} = \left[ {\tilde{B}_0 \over
  z}+{\tilde{B}_1\over z-1} \right] Y^{(2)}
$$
  $$ \tilde{B}_0=\hat{G}^{-1} A_0^0 \hat{G}_0 = \pmatrix{{\sigma\over
4} & {\sigma\over 4} s \cr
-{\sigma \over 4 s} & -{\sigma \over 4} \cr
  },~~~~~ \tilde{B}_1=\hat{G}^{-1} A_x^0 \hat{G}_0 = \pmatrix{{\sigma\over
4} & -{\sigma\over 4} s \cr
{\sigma \over 4 s} & -{\sigma \over 4} \cr
  }$$
This time the effective parameters at $z=0,1,\infty$ are 
$$
\left\{\matrix{\alpha_0=-{\sigma\over 2} \cr
                \beta_0={\sigma\over 2} +1 \cr
\gamma_0=1
}\right.  ,~~~~
\left\{\matrix{\alpha_1=-{\sigma\over 2} \cr
                \beta_1={\sigma\over 2} +1\cr
\gamma_1=1
}\right.  ,~~~~
\left\{\matrix{\alpha_{\infty}=-{\sigma\over 2} \cr
                \beta_{\infty}={\sigma\over 2} \cr
\gamma_{\infty}=\sigma
}\right.  
$$
If follows that both at $z=0$ and $z=1$ there are logarithmic
solutions. We skip  the derivation of the connection formulae,
which is done as in the previous cases, with some more technical
complications. Before giving the results we observe that 
$$
   Y^{(2)}(z)=\left(I+O\left({1\over z}\right)\right)
 z^{\hbox{diag}({\sigma\over
   2},-{\sigma\over 2})} ,~~~~z\to \infty
$$
$$
  = \hat{G}_0^{-1}\tilde{G}_0~(1+O(z))z^J~C_0^{\prime},~~~~z\to 0
$$
$$
 = \hat{G}_0^{-1}\tilde{G}_1~(1+O(z-1))(z-1)^J~C_1^{\prime},~~~~z\to 1
$$
where 
$$  
   C_i^{\prime}:= \tilde{C}_i\hat{G}_0,~~~~~i=0,1
$$
Then
 $$
        \tilde{M}_0=\hat{G}_0 ~(C_0^{\prime})^{-1} ~
\pmatrix{1&2\pi i \cr
0& 1}
~  C_0^{\prime}~\hat{G}_0^{-1},~~~~~
        \tilde{M}_1=\hat{G}_0 ~(C_1^{\prime})^{-1} ~
\pmatrix{1&2\pi i \cr
0& 1}
~  C_1^{\prime}~\hat{G}_0^{-1}
$$
So, we need to compute 
$C_i^{\prime}$, $i=0,1$. The result is

$$
C_0^{\prime} = 
\pmatrix{
           (C_{0\infty}^{\prime})_{11} & {\sigma\over \sigma+1}{s\over
                                                                      4}
           (C_{0\infty}^{\prime})_{12}             \cr
      (C_{0\infty}^{\prime})_{21}&{\sigma\over \sigma+1}{s\over
                                       4} (C_{0\infty}^{\prime})_{22}\cr
            }
, ~~~~
   C_1^{\prime}=C_{01}^{\prime} C_0^{\prime}$$

\noindent
where 
     $$ 
(C_{0\infty}^{\prime})^{-1}= 
\pmatrix{ 
{\Gamma(\beta_0-\alpha_0)\over
\Gamma(\beta_0)\Gamma(1-\alpha_0)} e^{i\pi\alpha_0} & 0 \cr
 { \Gamma(\alpha_0-\beta_0)\over
\Gamma(\alpha_0)\Gamma(1-\beta_0)} e^{i\pi\beta_0}&
-{\Gamma(1-\alpha_0)\Gamma(\beta_0)\over
\Gamma(\beta_0-\alpha_0+1)}~e^{i\pi \beta_0} \cr
}
,~~~
C_{01}^{\prime}= \pmatrix{
0 & -{\pi \over \sin(\pi\alpha_0)} \cr
                                 -{\sin(\pi\alpha_0)\over \pi} &
                                 -e^{-i\pi \alpha_0} \cr
                          }$$

\vskip 0.3 cm 

The case $\sigma\in{\bf Z}$ interests us only if $\sigma=0$, otherwise
$\sigma\not \in {\bf C}\backslash \{ (-\infty,0)\cup [1,+\infty)\}$.  
We observe that the system (\ref{systemhat})
is precisely the system for $Y^{(2)}(z)$ with the substitution
$\sigma\mapsto -2\mu$. In the formulae for $x_i^2$, $i=0,1,\infty$ 
 we only need
$C_{01}$, which is obtained from $C_{01}^{\prime}$ with
$\alpha_0=\mu$. 

As for the system (\ref{systemtilde}), the gauge
$Y^{(2)}=\hat{G}_0^{-1}\tilde{Y} \hat{G}_0$   
yields  $\tilde{B}_0=\pmatrix{0&s\cr 0 & 0\cr}$, $\tilde{B}_1= \pmatrix{
0 & 1-s \cr
0 & 0 \cr}$. Here $\hat{G}_0$ is the matrix such that $\hat{G}_0^{-1}
\Lambda \hat{G}_0= \pmatrix{0 & 1 \cr 0 & 0 \cr}$. The behavior of
$Y^{(2)}(z)$ is now:  
$$
 Y^{(2)}(z)=(I+O({1\over z}))~z^J ~~~~~~~~~~~z\to \infty
$$  
$$
  = \tilde{\tilde{G}}_0^{-1}~(1+O(z))z^J~C_0^{\prime},~~~~z\to 0
$$
$$
 = \tilde{\tilde{G}}_1~(1+O(z-1))(z-1)^J~C_1^{\prime},~~~~z\to 1
$$ 
Here  $\tilde{\tilde{G}}_i$ is  the matrix that puts 
$\tilde{B}_i$ in Jordan form, for $i=0,1$.  
$Y^{(2)}$ can be computed explicitly: 
$$ Y^{(2)}(z)= \pmatrix{1 & s\log(z)+(1-s)\log(z-1) \cr
                       0 & 1\cr
}
$$
If we choose   $\tilde{\tilde{G}}_0=$diag$(1,1/s)$, then  
$$
  C_0^{\prime}=\pmatrix{1&0\cr
                        0 & s\cr}
$$
In the same way we find 
$$
  C_1^{\prime}=\pmatrix{1&0\cr
                        0 & 1-s\cr}
$$

\vskip 0.3 cm 
To prove Theorem 2 it is now enough to compute 
$$
   2- x_0^2= \hbox{ tr}(M_0M_x)\equiv  \hbox{ tr}(e^{2\pi i J} 
(C_{01}^{\prime})^{-1} e^{2\pi i J} 
C_{01}^{\prime} )
$$
$$
2-x_1^2=  \hbox{ tr}(M_xM_1)\equiv  
 \hbox{ tr}((C_1^{\prime})^{-1}e^{2\pi i J} C_1^{\prime} C_{01}^{-1}
 e^{2\pi i J} C_{01})
$$

$$
2-x_{\infty}^2=  \hbox{ tr}(M_0M_1)\equiv  
 \hbox{ tr}((C_0^{\prime})^{-1}e^{2\pi i J} C_0^{\prime} C_{01}^{-1}
 e^{2\pi i J} C_{01})
$$
Note the remarkable simplifications obtained from the cyclic property
of the trace (for example, $\hat{C}_0$, $\hat{C}_1$   and $\hat{G}_0$ 
disappear). 
The fact that
$\hat{C}_0$ and 
$\hat{C}_1$  disappear implies that the formulae of Theorem 2 
are derived for any
$\mu\neq 0$, including the resonant cases. Thus, the connection
formulae in the resonant case $2\mu \in {\bf Z}\backslash \{0\}$ are
the same of the non-resonant case. 
The final result of the computation of the traces is: 

\vskip 0.15 cm
\noindent
 I) Generic case: 
\be
   \left\{ \matrix{2(1-cos(\pi \sigma))=x_0^2 \cr\cr
                  {1\over
                  f(\sigma,\mu)}\left(2+F(\sigma,\mu)~s+{1\over
                  F(\sigma,\mu)~s}\right) =x_1^2
                         \cr\cr
                  {1\over f(\sigma,\mu)}\left(
                         2-F(\sigma,\mu)e^{-i\pi\sigma}~s-{1 \over
                  F(\sigma,\mu)
e^{-i\pi\sigma}~s }
                         \right) =x_{\infty}^2\cr
                         }\right. 
\label{20ii}
\ee
where   
$$
f(\sigma,\mu)={2 \cos^2({\pi \over 2} \sigma) \over \cos(\pi \sigma)- 
                  \cos(2\pi\mu)}\equiv{4-x_0^2\over
    x_1^2+x_{\infty}^2-x_0x_1x_{\infty}}, ~~~
~~~ F(\sigma,\mu)=f(\sigma,\mu) {16^{\sigma}\Gamma({\sigma+1\over 2})^4\over
                  \Gamma(1-\mu+ {\sigma\over
                  2})^2\Gamma(\mu+{\sigma\over 2})^2}$$

\vskip 0.15 cm 
\noindent
 II) $\sigma\in 2{\bf Z}$, $x_0=0$. 
           $$
             \left\{
                     \matrix{ 2(1-cos(\pi \sigma))=0\cr\cr
                              4\sin^2(\pi\mu) ~(1-s)=x_1^2 \cr\cr 
                             4 \sin^2(\pi\mu) ~s=x_{\infty}^2               
                 }
 \right.
$$

\vskip 0.15 cm 
\noindent
III) $x_0^2=4 \sin^2(\pi\mu)$.  Then (\ref{10ii}) implies
 $x_{\infty}^2=-x_1^2 ~\exp(\pm 2\pi i \mu)$ . 
 Four cases which yield the values of $\sigma$ non included in I)
and II) must be considered 

\vskip 0.15 cm
   III1)  $x_{\infty}^2=-x_1^2 e^{- 2
\pi i \mu}$
$$\sigma=  2\mu + 2m,~~~~m=0,1,2,...$$ 
 $$ s=
               { \Gamma(m+1)^2 \Gamma(2\mu+m)^2 \over 16^{2\mu+2m}
 \Gamma(\mu+m+{1\over 2})^4}~x_1^2 
$$

\vskip 0.15 cm
  III2) $x_{\infty}^2=-x_1^2 e^{2\pi i \mu}$
$$\sigma=2\mu+2m,~~~~m=-1,-2,-3,...$$
    $$
            s={\pi^4\over \cos^4(\pi\mu)}\left[
                 16^{2\mu+2m} \Gamma(\mu+m+{1\over 2})^4
\Gamma(-2\mu-m+1)^2 \Gamma(-m)^2~ x_1^2 \right]^{-1} 
$$

\vskip 0.15 cm
  III3) $x_{\infty}^2=-x_1^2e^{2\pi i\mu}$
$$\sigma=-2\mu+2m,~~~~m=1,2,3,...$$
$$s=
         {
\Gamma(m-2\mu+1)^2 \Gamma(m)^2 \over   16^{-2\mu+2m} \Gamma(-\mu+m+{1\over 2})^4}~ x_1^2 
$$

\vskip 0.15 cm
III4)  $x_{\infty}^2=-x_1^2 e^{-2\pi i \mu}$
$$\sigma=-2\mu+2m,~~~~m=0,-1,-2,-3,...$$
$$
  s=             {\pi^4\over \cos^4(\pi\mu)} \left[16^{-2\mu+2m}
\Gamma(-\mu+m+{1\over 2})^4\Gamma(2\mu-m)^2 \Gamma(1-m)^2 ~ x_1^2
\right]^{-1} 
$$ 

\vskip 0.2 cm 
We recall that $a$ in $y(x;\sigma,a)$ is $a=-{1\over 4s}$ in general, 
and $a=s$ for $\sigma=0$. 

\vskip 0.2 cm 
 To compute $\sigma$ and $s$ in  the generic case $I$) for a given 
triple $(x_0,x_1,x_{\infty})$, 
 we solve the system (\ref{20ii}). It has two unknowns  and three equations and
we need to prove that it is 
compatible. Actually, the first
equation $ 2(1-cos(\pi \sigma))=x_0^2$ has always solutions. Let us choose
a solution 
$\sigma_0$  ($\pm \sigma_0 + 2 n$, $\forall n\in {\bf Z}$
are also solutions). Substitute it in the last two equations. We need
to verify they are compatible.  
Instead of $s$ and ${1\over s}$ write $X$ and $Y$. We have the linear
system in two variable $X$, $Y$
$$ 
                 \pmatrix{ F(\sigma_0) & {1 \over F(\sigma_0)} \cr\cr
                            F(\sigma_0)~e^{-i \pi \sigma_0}& {1\over
                           F(\sigma_0)} e^{-i\pi\sigma_0} \cr
}
 \pmatrix{X\cr \cr Y\cr}=
                          \pmatrix{
                                    f(\sigma_0) ~x_1^2-2 \cr
\cr
                                    2 - f(\sigma_0) ~x_{\infty}^2 \cr
}$$
 The system has a unique solution if and only if $2 i \sin(\pi \sigma_0)=\det 
     \pmatrix{ F(\sigma_0) & {1 \over F(\sigma_0)} \cr\
                           F(\sigma_0)~e^{-i \pi \sigma_0}& {1\over
                           F(\sigma_0)} e^{-i\pi\sigma_0} \cr
}\neq 0$. This happens for $\sigma_0 \not\in {\bf Z}$. The condition
                           is
                            not restrictive, because for $\sigma$ even
                           we
                           turn to the case $II$), 
 and $\sigma$ odd is not in 
${\bf C}\backslash [ (-\infty,0)\cup [1,+\infty)]$. 
   The solution
                           is then 
$$ X= {2 (1+e^{-i \pi\sigma_0})-f(\sigma_0)(x_1^2+x_{\infty}^2 e^{-i \pi
\sigma_0})\over F(\sigma_0)(e^{-2 \pi i \sigma_0} -1)}$$
$$  
   Y=
   F(\sigma_0)~{f(\sigma_0)e^{-i\pi\sigma_0}(e^{-i\pi\sigma_0}x_1^2+x_{\infty}^2)-2 e^{-i\pi\sigma_0}(1+e^{-i\pi\sigma_0})\over e^{-2\pi i \sigma_0}-1}
$$
Compatibility of the system means that $X~Y\equiv 1$. This is
verified  by direct computation.
 $$ 
   XY= {e^{-i\pi\sigma}~\left[
   2(1+e^{-i\pi\sigma})-(x_1^2+x_{\infty}^2e^{-i\pi\sigma})f(\sigma) 
\right]~\left[
           (x_1^2e^{-i\pi\sigma}+x_{\infty}^2)f(\sigma)- 
 2(1+e^{-i\pi\sigma}) 
\right] \over (e^{-2i\pi\sigma}-1)^2}$$
$$
    = { 8 \cos^2({\pi\sigma\over
    2})(x_1^2+x_{\infty}^2)f(\sigma)
    -4(4-\sin^2({\pi\sigma\over
    2}))-((x_1^2+x_{\infty}^2)^2-x_0^2x_1^2x_{\infty}^2) f(\sigma)^2
    \over -4 \sin^2(\pi\sigma)}$$
Using the relations $\cos^2({\pi \sigma \over 2})=1-x_0^2/4$,
    $\cos(\pi\sigma) =1-x_0^2/2$ and $f(\sigma)= {4-x_0^2\over
    x_1^2+x_{\infty}^2-x_0x_1x_{\infty}}$ we obtain 
$$
 XY= {1\over x_0^2}~\left(
    -2(x_1^2+x_{\infty}^2)f(\sigma)+4+{(x_1^2+x_{\infty}^2)^2- 
        (x_0x_1x_{\infty})^2 \over
    x_1^2+x_{\infty}^2-x_0x_1x_{\infty}} f(\sigma)  
                         \right)
$$
$$ 
  ={1 \over x_0^2}~\left(4 -(x_1^2+x_{\infty}^2-x_0x_1x_{\infty})f(\sigma)
  \right)= {1\over x_0^2} \left(4-(4-x_0^2)\right)=1
 $$

It follows from this construction that for any $\sigma$ solution of
the first equation of (\ref{20ii}), there always exists a unique $s$
which solves the last two equations.

\vskip 0.2 cm 

To complete the proof of Theorem 2 (points $i)$, $ii)$, $iii)$), 
 we just have to compute the square
roots of the $x_i^2$ ($i=0,1,\infty$) in such a way that (\ref{10ii}) is
satisfied.  For example, the square root of I) satisfying (\ref{10ii}) is 
$$
   \left\{ \matrix{x_0=2 \sin({\pi \over 2}\sigma) \cr\cr
                  x_1={1\over \sqrt{f(\sigma,\mu)}}
                    \left(\sqrt{F(\sigma,\mu)~s}+{1\over
 \sqrt{F(\sigma,\mu)~s}
                 }\right)
                      \cr\cr
                  x_{\infty}={i\over
 \sqrt{f(\sigma,\mu)}}\left(\sqrt{F(\sigma,\mu)~s}~e^{-i{\pi \sigma\over
 2}}-{1\over \sqrt{F(\sigma,\mu)~s}~e^{-i{\pi \sigma\over 2}}}\right)\cr
                         }\right. 
$$
which yields i), with $F(\sigma,\mu)=f(\sigma,\mu) (2 G(\sigma,\mu))^2$. . 

\vskip 0.2 cm

We remark that in case II) only $\sigma =0$ is in ${\bf C}\backslash 
\{ (-\infty,0)\cup[1,+\infty)\}$. If $\mu$ integer
in II), 
the formulae  give $(x_0,x_1,x_{\infty})=(0,0,0)$. The
triple is not admissible, and direct computation gives $R=0$ for the
system (\ref{novosis}). This is the case of commuting monodromy
matrices with a 1-parameter family of rational solutions of
$PVI_{\mu}$. 

\vskip 0.2 cm 
The last remark concerns the choice of (\ref{cucu}), (\ref{cucu1})
instead of (\ref{cucu2}), (\ref{cucu3}). The reason is that at $z=0$
the system  (\ref{novosis})
has solution  
corresponding to (\ref{ffirr1}). This is true for any
$\sigma \neq 0$ in  ${\bf C}\backslash 
\{ (-\infty,0)\cup[1,+\infty)\}$, also for $\sigma \to \pm 2\mu$. 
Its behavior is (\ref{ffirr2}), which is obtainable from
the $\hat{G}_0 = T$ of (\ref{cucu}), (\ref{cucu1}) but not of 
(\ref{cucu2}), (\ref{cucu3}). See also the comment following Lemma 5. 

\rightline{$\Box$}

\vskip 0.3 cm
\noindent
{\it Remark:}
In the proof of Theorem 2 we take the limits 
of the system and of the rescaled system for $x\to 0$ in $D(\sigma)$. 
At $x$ we assign the
monodromy $M_0, M_1, M_x$ characterized by $(x_0,x_1,x_{\infty})$ and then we
take the limit proving the theorem. 
If we  start from another point $x^{\prime}\in D(\sigma)$ 
 we have to  choose 
 the same monodromy $M_0, M_1, M_x$, because what we are 
doing is  the limit for $x\to 0$ in $D(\sigma)$  of the matrix coefficient 
$A(z,x;x_0,x_1,x_{\infty})$ of the system (\ref{stofucs})  
considered as a function defined on the universal
covering of ${\bf C}_0\cap \{|x|<\epsilon\}$.

 \vskip 0.3 cm 
\noindent
{\it  Proof of Remark 2 of section \ref{Parametrization of a branch through Monodromy Data --
Theorem 2}:}

We prove that $a(\sigma)={1\over 16 a(-\sigma)}$, namely 
$s(\sigma)={1\over s(-\sigma)}$  ($a=-{1\over 4s}$). 
 Given monodromy data $(x_0, x_1,
 x_{\infty})$ the parameter  $s$
 corresponding to $\sigma$ is uniquely determined by  
 $$
  {1\over
                  f(\sigma)}\left(2+F(\sigma)~s+{1\over
                  F(\sigma)~s}\right) =x_1^2
  $$
 $$                      
                  {1\over f(\sigma)}\left(
                         2-F(\sigma)e^{-i\pi\sigma}~s-{1 \over
                  F(\sigma)
 e^{-i\pi\sigma}~s }
                          \right) =x_{\infty}^2$$
   We observe that $f(\sigma)=f(-\sigma)$ and that  
 the properties of the Gamma function 
 $$
   \Gamma(1-z)\Gamma(z)= {\pi \over \sin(\pi z)},~~~~~
  \Gamma(z+1)=z\Gamma(z)$$
 imply
 $$F(-\sigma)= {1 \over F(\sigma)}$$
 Then the value of $s$ corresponding to $-\sigma$ is (uniquely) determined by 
$${1\over
                   f(\sigma)}\left(2+{s\over F(\sigma)}+{F(\sigma)\over
                  s}\right) =x_1^2
$$                      
 $$
     {1\over f(\sigma)}\left(
                         2-{s \over F(\sigma)e^{-i\pi\sigma}}-{
                  F(\sigma)
 e^{-i\pi\sigma}\over s }
                         \right) =x_{\infty}^2
$$
 We conclude that $s(-\sigma)= -{1\over s(\sigma)}$.

\vskip 0.3 cm

 \vskip 0.3 cm 
\noindent
{\it  Proof of formula (\ref{analy}): }

 We are ready to prove formula (\ref{analy}), namely: 
 $$\beta_1^2:~~(\sigma,a)\mapsto(\sigma,a e^{-2\pi i \sigma}) $$
 For $\sigma=0$ we have $x_0=0$ and $\beta_1^2:(0,x_1,x_{\infty})
 \mapsto ( 0,x_1,x_{\infty})$. Thus 
$$ 
a={x_{\infty}^2\over x_1^2+x_{\infty}^2}\mapsto{x_{\infty}^2\over
 x_1^2+x_{\infty}^2}\equiv a$$

\noindent
 For $\sigma=\pm 2 \mu +2 m$, we consider the example
 $\sigma=2 \mu+2m$, $m=0,1,2,...$. The other cases are analogous. We have $s=
 x_1^2~H(\sigma)=-x_{\infty}^2H(\sigma)e^{2\pi i \mu}$, where the
 function $H(\sigma)$ is explicitly given in theorem 2, $III$). Then
$$ 
  \beta_1:~~s= -x_{\infty}^2 H(\sigma)e^{2\pi i \mu} \mapsto -x_1^2
   H(\sigma) e^{2\pi i \mu}= -s e^{2\pi i \mu}
$$
Then 
$$ \beta_1^2:~~s\mapsto se^{4\pi i \mu}~~~ \Longrightarrow~~~ a\mapsto
ae^{-4\pi i \mu} \equiv a e^{-2\pi i \sigma}$$

\vskip 0.15 cm
\noindent 
 For  the generic case $I$) ($\sigma\not\in{\bf Z}$, $\sigma \neq \pm
 2\mu +2m$) recall that 
$$
\left\{ \matrix{
                 F(\sigma)~s+{1\over F(\sigma)~s} = x_1^2 f(\sigma)-2
                 \cr\cr
               F(\sigma)e^{-i\pi\sigma}~s+{1\over F(\sigma)e^{-i\pi
                 \sigma}~s } = 2-x_{\infty}^2 f(\sigma)
}\right.
$$
has a unique solution $s$. 
Also observe that $\beta_1: x_{\infty}\mapsto x_1$. Then 
  the transformed parameter $\beta_1:~s\mapsto s^{\beta_1}$
satisfies the equation 
$$ 
     F(\sigma)e^{-i\pi\sigma}~s^{\beta_1}+{1\over F(\sigma)e^{-i\pi
                 \sigma}~s^{\beta_1} }=  2-x_{1}^2 f(\sigma)
                $$
$$ \equiv  
                 -\left(    F(\sigma)~s+{1\over F(\sigma)~s} \right)$$
Thus $s^{\beta_1}=-e^{i\pi\sigma}~ s$. This implies 
$$
  \beta_1^2:~s\mapsto s e^{2 \pi i \sigma}~~~\Longrightarrow~~~ a \mapsto
  a~e^{-2\pi i \sigma}
$$

\rightline{$\Box$}

\vskip 0.3 cm

 We finally  prove the Proposition stated at the end of  Section \ref{Parametrization of a branch through Monodromy Data --
Theorem 2}. 

\vskip 0.2 cm
\noindent
{\it  Proof:}  Observe that both $y(x)$ and $y(x;\sigma,a)$ have the same
asymptotic behavior for $x\to 0$ in $D(\sigma)$. Let  
 $A_0(x)$, $A_1(x)$, $A_x(x)$ be the matrices constructed from
 $y(x)$ and $A^*_0(x)$, $A^*_1(x)$, $A^*_x(x)$ constructed from
 $y(x;\sigma,a)$ by means of   
 the formulae (\ref{definite in extremis}). 
 It follows that $A_i(x)$ and $A^*_i(x)$, $i=0,1,x$, 
 have the same asymptotic behavior as
 $x\to 0$. This is the
behavior of Lemma 1 of section \ref{proof of theorem 1} 
(adapted to our case). From the proof of Theorem 2 if follows that $A_0(x)$,
$A_1(x)$, $A_x(x)$  and $A^*_0(x)$, $A^*_1(x)$, $A^*_x(x)$ produce   the same
 triple 
 $(x_0,x_1,x_{\infty})$. 
The solution of the Riemann-Hilbert problem for such a triple 
is unique, up to conjugation of the fuchsian systems.  
Therefore $A_i(x)$ and $ A^*_i(x)$, $i=0,1,x$ are conjugated. If $2\mu \not \in {\bf Z}$ the conjugation is diagonal. If $2\mu \in {\bf Z}$ and $R\neq 0$, then  $A_i(x) =  A^*_i(x)$  (see notes 2 and 3 for details).  Putting $[A(z;x)]_{12}=0$ and $ [A^*(z;x)]_{12}=0$ we conclude that 
$y(x)\equiv y(x;\sigma,a)$.  \rightline{$\Box$}


\section{ Proof of Theorem 3}\label{provaellittica}

The elliptic representation was derived by R. Fuchs in \cite{fuchs}. In the 
case of $PVI_{\mu}$ the representation is discussed at the beginning of 
sub-section   \ref{Elliptic Representation}. 
Here we study the solutions of (\ref{difficilissima}).

 To start with, we
 derive the elliptic form for the general Painlev\'e 6 equation. We
follow \cite{fuchs}. 
 We put
\be
u=\int_{\infty}^y ~{d\lambda \over \sqrt{\lambda(\lambda-1)(\lambda-x)}}
\label{elliptint1}
\ee
We observe that
$$
  {du \over dx}= {\partial u \over \partial y }~{dy \over dx} + {\partial u
  \over \partial x } ={1\over \sqrt{y(y-1)(y-x)} } ~
{dy \over dx} + {\partial u
  \over \partial x } 
$$
from which we compute
$$
{d^2 u\over dx^2}+{2x-1\over x(x-1)} {du\over dx} +{u\over 4 x(x-1)}=$$
$$
  = {1\over  \sqrt{y(y-1)(y-x)}}\left[{d^2 y\over dx^2}+\left({1\over x}
  +{1\over x-1}+{1\over y-x}\right) ~{dy \over dx}-{1\over 2}\left({1\over y}
  +{1\over y-1}+{1\over y-x}\right)~\left({dy \over dx}\right)^2 \right]
$$
$$
+{\partial^2 u \over \partial x^2}+{2x-1\over x(x-1)}~{\partial u \over
  \partial x}+ {u\over 4x(x-1)} 
$$
By direct calculation we have:
$$
{\partial^2 u \over \partial x^2}+{2x-1\over x(x-1)}~{\partial u \over
  \partial x}+ {u\over 4x(x-1)} 
= -{1\over 2} {\sqrt{y(y-1)(y-x)} \over x(x-1)}
  ~{1\over (y-x)^2}
$$
 Therefore, $y(x)$ satisfies the Painlev\'e 6 equation if and only if 
\be
 {d^2 u\over dx^2}+{2x-1\over x(x-1)}~ {du\over dx} +{u\over 4
x(x-1)}={\sqrt{y(y-1)(y-x)} \over 2 x^2 (1-x)^2}~\left[2 \alpha+2\beta{x\over
y^2} +\gamma {x-1\over (y-1)^2} +\left(\delta-{1\over2}\right){x(x-1)\over
(y-x)^2} \right]
\label{ellipainleve1}
\ee

\vskip 0.2 cm
We invert the function $u=u(y)$ by observing that we are dealing with an
elliptic integral. Therefore, we write
$$ 
  y=f(u,x)$$
where $f(u,x)$ is an elliptic function of $u$. This implies that 
$$
  {\partial y\over \partial u}=\sqrt{y(y-1)(y-x)}
$$
 The above equality allows us to  rewrite (\ref{ellipainleve1}) in the
 following way: 
\be
x(1-x)~{d^2 u\over dx^2}+(1-2x)~{du\over dx}-{1\over 4}~u = {1\over 2x(1-x)}
 ~{\partial \over \partial u}\psi(u,x),
\label{ellipainleve2}
\ee
where 
$$
\psi(u,x):= 2 \alpha f(u,x)-2\beta{x\over f(u,x)} +2\gamma{1-x\over f(u,x)-1}
+(1-2\delta){x(x-1)\over f(u,x)-x}
$$

\vskip 0.2 cm
The last step concerns the form of $f(u,x)$. We observe that 
$4\lambda(\lambda-1)(\lambda-x)$ is not in Weierstrass canonical form. We change variable:
$$\lambda= t +{1+x\over 3},$$
and we get the Weierstrass form:
$$
    4\lambda(\lambda-1)(\lambda-x)= 4 t^3-g_2 t-g_3,~~~~ 
 g_2={4\over 3}(1-x+x^2),~~~~g_3:={4\over
    27}(x-2)(2x-1)(1+x)
$$
Thus
$$
  {u\over 2}=\int_{\infty}^{y-{1+x\over 3}}~{dt\over  \sqrt{4 t^3-g_2 t-g_3}} 
$$
which implies 
$$
   y(x)= \wp\left({u\over 2}; \omega_1,\omega_2\right)+{1+x\over 3}
$$
 We still need to explain what are the {\it half periods} $\omega_1$,
 $\omega_2$. In order to do that, we first observe that the Weierstrass form is 
$$
   4 t^3-g_2 t-g_3= 4(t-e_1)(t-e_2)(t-e_3)
$$
where
$$ 
    e_1= {2-x\over 3},~~~e_2={2x-1\over 3},~~~e_3=-{1+x\over 3}.
$$
Therefore
$$
   g:=\sqrt{e_1-e_2}=1,~~~{\kappa}^2:={e_2-e_3\over e_1-e_3} =
   x,~~~{{\kappa}^{\prime}}^2 := 1- {\kappa}^2=1-x
$$
and the half-periods are 
$$
\omega_1= {1\over g} \int_0^1 ~{d\xi\over \sqrt{(1-\xi^2)(1-\kappa^2 \xi^2)}} 
= \int_0^1 ~{d\xi\over \sqrt{(1-\xi^2)(1-x \xi^2)}} ={\bf K}(x)
$$
$$
 \omega_2= {i\over g} \int_0^1 ~{d\xi\over
\sqrt{(1-\xi^2)(1-{\kappa^{\prime}}^2  
\xi^2)}} 
= i~\int_0^1 ~{d\xi\over \sqrt{(1-\xi^2)(1-(1-x) \xi^2)}} =i{\bf
K}(1-x) 
$$
The elliptic integral ${\bf K}(x)$ is known:
$$
{\bf K}(x)= {\pi \over 2} ~F\left({1\over 2},{1\over 2}, 1; x\right),~~~~
F\left({1\over 2},{1\over 2},1;x\right)=  
\sum_{n=0}^{\infty}{ \left[\left({1 \over 2}\right)_n\right]^2
  \over (n!)^2 } x^n.
$$
${\bf K}(x)$ and ${\bf K}(1-x)$ 
 are two linearly independent solutions of the hyper-geometric equation 
$$
 x(1-x) \omega^{\prime\prime}+(1-2x) \omega^{\prime} -{1\over 4} \omega=0.
$$
Observe that for $|\hbox{arg}(x)|<\pi$: 
$$-\pi F\left({1\over 2},{1\over 2}, 1;1- x\right) =  F\left({1\over
2},{1\over 2},1;x\right)  \ln(x) +
F_1(x)  
$$ 
where 
$$
F_1(x):=  
\sum_{n=0}^{\infty}{ \left[\left({1 \over 2}\right)_n\right]^2
  \over (n!)^2 } 2\left[ \psi(n+{1\over 2}) - \psi(n+1)\right]
x^n,~~~~\psi(z) = 
{d \over dz}\ln \Gamma(z).
$$
 Therefore $\omega_2(x)= -{i\over 2}[F(x)\ln(x)+F_1(x)]
$ where $F(x)$ is a abbreviation for $ F\left({1\over 2},{1\over 2},
1;x\right)$. The series of $F(x)$ and $F_1(x)$ 
 converge for $|x|<1$. 
We conclude that $PVI_{\mu}$ is equivalent to (\ref{difficilissima}). 

 Incidentally, we observe that 
$$
  y(x)= \wp\left({u(x)\over 2}; \omega_1(x),\omega_2(x)\right)-e_3= 
{1\over \hbox{sn}^2\left({u(x)\over 2},\kappa^2=x\right)}
$$

\vskip 0.3 cm 
\noindent
{\it Proof of theorem 3:}  We let $x\to 0$. If $ \Im \tau>0$  and
\be 
\left|\Im \left({u\over 4\omega_1}\right)\right|<\Im \tau
\label{speriamobene}
\ee
we expand the elliptic function in Fourier series 
(\ref{rondine}).  The first condition $\Im \tau>0$ is always satisfied for
$x\to 0$ because 
$$ 
    \Im \tau(x)= -{1\over \pi} \ln|x| +{4\over \pi} \ln 2 +O(x), ~~~~x\to 0.
$$
Therefore, in the following we assume that $|x|<\epsilon<1$ for a sufficiently 
small $\epsilon$. 
We look for a solution $u(x)$ of (\ref{difficilissima}) of the form 
$$
  u(x)= 2 \nu_1 \omega_1(x)+2 \nu_2 \omega_2(x) + 2 v(x)
$$
where $v(x)$ is a (small) perturbation 
 to be determined from (\ref{difficilissima}).  
 We observe that  
$$
   {u(x)\over 4 \omega_1(x)} = {\nu_1\over 2} +{\nu_2\over 2} \tau(x) +
   {v(x)\over 2\omega_1(x)} 
   = {\nu_1\over 2} +{\nu_2\over 2} \left[-{i\over \pi} \ln x - {i\over \pi}
 {F_1(x)\over F(x)}\right]+{v(x)\over 2\omega_1(x)}
$$ 
Note that for $x\to 0$ 
  ${F_1(x)\over F(x)} =-4\ln 2 +g(x)$, 
where $g(x)=O(x)$ is a convergent Taylor series starting with $x$.  Thus,  the 
 condition (\ref{speriamobene}) becomes 
\be
(2+\Re \nu_2) \ln|x|-{\cal C}(x,\nu_1,\nu_2) -8 \ln 2  < \Im \nu_2
\arg(x)< (\Re\nu_2-2) \ln|x|- {\cal C}(x,\nu_1,\nu_2)  +8 \ln 2, 
\label{speriamobene1}
\ee
where $ {\cal C}(x,\nu_1,\nu_2)=
[\Im {\pi v\over \omega_1}+4\ln 2~ \Re \nu_2 + \pi ~\Im
\nu_1 +O(x)]$. 
We expand the derivative of $\wp$ appearing in (\ref{difficilissima})
$$
   {\partial \over \partial u} \wp\left({u\over 2};
   \omega_1,\omega_2\right) = \left({\pi\over \omega_1}\right)^3
   \sum_{n=1}^{\infty} {n^2 e^{2\pi i n \tau}\over 1-e^{2\pi i n\tau}} \sin
   \left( {n
  \pi u \over 2 \omega_1}\right) - \left({\pi \over 2
   \omega_1}\right)^3{\cos\left({\pi u 
   \over 4\omega_1}\right) \over \sin^3\left({\pi u \over 4 \omega_1}\right)}
$$
$$
= {1\over 2i}  \left({\pi\over \omega_1}\right)^3\sum_{n=1}^{\infty} {
n^2 e^{2\pi i n \tau}\over 1-e^{2\pi
i n\tau}} \left(e^{in{\pi u \over 2 \omega_1}}-e^{-in{\pi u \over 2 \omega_1}}
\right)+4 i  \left({\pi \over 2
   \omega_1}\right)^3 {e^{i {\pi u \over 4 \omega_1}} + e^{-i {\pi u \over 4
   \omega_1}}\over \left( e^{i {\pi u \over 4 \omega_1}} - e^{-i {\pi u \over 4
   \omega_1}}  \right)^3 }
$$
Now we come to a crucial step in the construction: we collect $ e^{-i {\pi u \over 4
   \omega_1}} $ in the last term, which becomes
$$ 
4 i  \left({\pi \over 2
   \omega_1}\right)^3 {e^{4 \pi i { u \over 4 \omega_1}} + e^{2\pi i
 { u \over 4
   \omega_1}}\over \left( e^{2 \pi i { u \over 4 \omega_1}} - 1  \right)^3 }.
$$
The denominator {\it does not vanish} if $ \left|e^{2 \pi i { u \over 4
\omega_1}}\right|<1$. 
From now on, this condition is added to (\ref{speriamobene})
and reduces the domain (\ref{speriamobene1}). The expansion of $ {\partial
\over \partial u} \wp$ becomes
$$
 {\partial \over \partial u} \wp\left({u\over 2};
   \omega_1,\omega_2\right) ={1\over 2i}  \left({\pi\over
   \omega_1}\right)^3\sum_{n=1}^{\infty} {n^2 e^{i \pi   n
   \left[-\nu_1+(2-\nu_2) \tau-{v\over \omega_1}\right]}\over 1-e^{2\pi i n
   \tau}} \left(e^{2 i \pi n \left[ \nu_1+\nu_2 \tau +{v\over  \omega_1}
   \right]} -1\right)
$$
$$
+ 4i \left({\pi \over 2
   \omega_1}\right)^3{ e^{2\pi i \left[  \nu_1+\nu_2 \tau +{v\over  \omega_1}
   \right]}+ e^{\pi i  \left[  \nu_1+\nu_2 \tau +{v\over  \omega_1}
   \right]} \over \left(
 e^{\pi i  \left[  \nu_1+\nu_2 \tau +{v\over  \omega_1}
   \right]} -1\right)^3}
$$
We observe that 
$$e^{i \pi C\tau}= {x^C\over 16^C}
~e^{C~g(x)}= {x^C\over 16^C} (1+O(x)), ~~~x\to 0,~~
 ~~~~\hbox{ for any } C\in {\bf C}. 
$$ 
Hence
$$
{\partial \over \partial u} \wp\left({u\over 2};
   \omega_1,\omega_2\right)= {\cal F}\left(x, {e^{-i\pi \nu_1} \over
   16^{2-\nu_2} } 
   x^{2-\nu_2} e^{-i\pi {v\over \omega_1}},{e^{i\pi \nu_1}\over 16^{\nu_2}
   }x^{\nu_2} e^{i \pi {v\over \omega_1}} \right)
$$
where 
$$ 
  {\cal F}(x,y,z)=  {1\over 2i} \left({\pi \over 
   \omega_1(x)}\right)^3 \sum_{n=1}^{\infty} {n^2 e^{n(2-\nu_2)g(x)}\over 
                  1- \left[{1\over 16} e^{g(x)}\right]^{2n} x^{2n}} y^n(
e^{2n \nu_2 g(x)}  z^{2n}-1)
   +4i\left({\pi \over 2\omega_1(x)} \right)^3 {e^{2\nu_2 g(x)}z^2+
  e^{\nu_2 g(x)} z\over (e^{\nu_2 g(x) }z-1)^3}
$$
The series converges for $|x|<\epsilon$ and for $|y|<1$, $|yz|<1$;
this is precisely (\ref{speriamobene}). However, we require that the last term
is holomorphic, so we have to further impose $|e^{\nu_2 g(x)}~z|<1$. 
On the resulting domain
$|x|<\epsilon$, $|y|<1$, $|e^{\nu_2 g(x)}~z|<1$,  ${\cal F}(x,y,z)$ is 
holomorphic and satisfies 
$$
  {\cal F}(0,0,0)=0. 
$$
   The condition  $|y|<1$, $|e^{\nu_2 g(x)}~z|<1$ is 
$\left|{e^{-i\pi \nu_1} \over
   16^{2-\nu_2} } 
   x^{2-\nu_2} e^{-i\pi {v\over \omega_1}}\right|<1$, $\left|e^{\nu_2 g(x)}~
{e^{i\pi \nu_1}
\over 16^{\nu_2}
   }x^{\nu_2} e^{i \pi {v\over \omega_1}}\right|<1$  , namely
\be
\Re \nu_2 \ln |x| -{\cal C}(x,\nu_1,\nu_2) 
< \Im \nu_2 \arg(x)< (\Re\nu_2-2) \ln |x|
   -{\cal C}(x,\nu_1,\nu_2) +8 \ln 2,
\label{speriamobene2}
\ee
which is more restrictive that (\ref{speriamobene1}). For 
$\Im \nu_2=0$ any value of $\arg(x)$ is allowed, but 
$\left|{e^{-i\pi \nu_1} \over
   16^{2-\nu_2} } 
   x^{2-\nu_2} e^{-i\pi {v\over \omega_1}}\right|<1$, $\left|e^{\nu_2 g(x)}~
{e^{i\pi \nu_1}
\over 16^{\nu_2}
   }x^{\nu_2} e^{i \pi {v\over \omega_1}}\right|<1$ imply 
$$
     0<\nu_2<2.
$$
Thus, $\nu_2=0$ is not allowed.

The function ${\cal F}$ can be decomposed as follows: 
$$
 {\cal F}={\cal F}\left(x, {e^{-i\pi \nu_1} \over
   16^{2-\nu_2} } 
   x^{2-\nu_2},{e^{i\pi \nu_1}\over 16^{\nu_2}
   }x^{\nu_2} \right)+
$$
$$
+\left[ {\cal F}\left(x, {e^{-i\pi \nu_1} \over
   16^{2-\nu_2} } 
   x^{2-\nu_2} e^{-i\pi {v\over 2\omega_1}},{e^{i\pi \nu_1}\over 16^{\nu_2}
   }x^{\nu_2} e^{i \pi {v\over 2\omega_1}} \right)- {\cal F}\left(x, {e^{-i\pi \nu_1} \over
   16^{2-\nu_2} } 
   x^{2-\nu_2},{e^{i\pi \nu_1}\over 16^{\nu_2}
   }x^{\nu_2} \right)\right]
$$
$$
  =: {\cal F}\left(x, {e^{-i\pi \nu_1} \over
   16^{2-\nu_2} } 
   x^{2-\nu_2},{e^{i\pi \nu_1}\over 16^{\nu_2}
   }x^{\nu_2} \right)+{\cal G} \left(x, {e^{-i\pi \nu_1} \over
   16^{2-\nu_2} } 
   x^{2-\nu_2},{e^{i\pi \nu_1}\over 16^{\nu_2}
   }x^{\nu_2},v(x) \right)
$$
The above defines ${\cal G}(x,y,z,v)$. It 
 is holomorphic for $|x|$,
$|y|$, $|z|$,  $|v|$ less than a sufficiently small $\epsilon<1$. Moreover
$$ 
  {\cal G}(0,0,0,v)={\cal G}(x,y,z,0)=0.
$$

 Let us put $u=u_0+2v$, where $u_0= 2 \nu_1
 \omega_1+2\nu_2\omega_2$. Therefore ${\cal L}(u_0)=0$ and ${\cal L}(u_0+2v)=
 {\cal L}(u_0)+{\cal L}(2v)\equiv 2{\cal L}(v)$. Hence (\ref{difficilissima})
 becomes 
\be
   {\cal L}(v)= {\alpha \over 2x(1-x)} ({\cal F}+{\cal G}), 
\label{zazen}
\ee
where ${\cal F}=  {\cal F}\left(x, {e^{-i\pi \nu_1} \over
   16^{2-\nu_2} } 
   x^{2-\nu_2},{e^{i\pi \nu_1}\over 16^{\nu_2}
   }x^{\nu_2} \right)$, ${\cal G}= 
{\cal G} \left(x, {e^{-i\pi \nu_1} \over
   16^{2-\nu_2} } 
   x^{2-\nu_2},{e^{i\pi \nu_1}\over 16^{\nu_2}
   }x^{\nu_2},v(x) \right)
$. 
 We put 
$$
   w:= x v^{\prime} ~~~(\hbox{where } v^{\prime}= {d v\over dx}),
$$ 
and the equation (\ref{zazen}) becomes 
$$
   w^{\prime} = {1\over x} \left[ {\alpha \over 2(1-x)^2} {\cal F} + {x
   (w+{1\over 4} v) \over 1-x} +{\alpha \over 2(1-x)^2 } {\cal G} \right]
$$
Now, let us define
$$
    \Phi(x,y,z):=  {\alpha \over 2(1-x)^2} {\cal F}(x,y,z),$$
$$
   \Psi(x,y,z,v,w):= {x
   (w+{1\over 4} v) \over 1-x} +{\alpha \over 2(1-x)^2 } {\cal G}(x,y,z,v).
$$
They are holomorphic for $|x|,|y|,|z|,|v|,|w|$ less than $\epsilon$ 
 and 
$$
   \Phi(0,0,0)=0,~~~~\Psi(0,0,0,v,w)=\Psi(x,y,z,0,0)=0.
$$
Equation (\ref{difficilissima}) becomes the system 
$$ 
    x {d v \over dx} = w,
$$
$$
    x {dw\over dx}= \Phi (x, {e^{-i\pi \nu_1}\over 16^{2-\nu_2}} x^{2-\nu_2}, 
{e^{i\pi \nu_1} \over 16^{\nu_2}} x^{\nu_2} ) + \Psi(x, {e^{-i\pi \nu_1}\over 16^{2-\nu_2}} x^{2-\nu_2}, 
{e^{i\pi \nu_1} \over 16^{\nu_2}} x^{\nu_2},v(x),w(x)).
$$
 We reduce it to a system of integral equations 
$$
w(x)=\int_{L(x)} {1\over s} \left\{
 \Phi (s, {e^{-i\pi \nu_1}\over 16^{2-\nu_2}} s^{2-\nu_2}, 
{e^{i\pi \nu_1} \over 16^{\nu_2}} s^{\nu_2} ) + \Psi(s, {e^{-i\pi \nu_1}\over 16^{2-\nu_2}} s^{2-\nu_2}, 
{e^{i\pi \nu_1} \over 16^{\nu_2}} s^{\nu_2},v(s),w(s))\right\}~ds
$$
$$
v(x)=\int_{L(x)}{1\over s}~\int_{L(s)} {1\over t}\left\{
  \Phi (t, {e^{-i\pi \nu_1}\over 16^{2-\nu_2}} t^{2-\nu_2}, 
{e^{i\pi \nu_1} \over 16^{\nu_2}} t^{\nu_2} ) + \Psi(t, {e^{-i\pi \nu_1}\over 16^{2-\nu_2}} t^{2-\nu_2}, 
{e^{i\pi \nu_1} \over 16^{\nu_2}} t^{\nu_2},v(t),w(t))\right\}~dt~ds
$$
The point $x$ and 
the path of integration  are chosen to belong to the domain where $|x|$, 
$|{e^{-i\pi \nu_1}\over 16^{2-\nu_2} } x^{2-\nu_2}|$, $|{e^{i\pi \nu_1} \over 16^{\nu_2}} x^{\nu_2}|$, $|v(x)|$, $|w(x)|$ are less than $\epsilon$,
 in such a way that $\Phi$ and $\Psi$ are holomorphic. 
That such a domain is not empty will be shown below. In particular, 
we'll show that if we require that  
 $|x|<r$, $|{e^{-i\pi \nu_1}\over 16^{2-\nu_2}} x^{2-\nu_2}|<r$, 
$|{e^{i\pi \nu_1} \over 16^{\nu_2}} x^{\nu_2}|<r$, where $r<\epsilon$ is 
small enough,  
also $|v(x)|$ and $|w(x)|$ are less than $\epsilon$. Such a domain is 
precisely the domain of Theorem 3, which is contained in 
(\ref{speriamobene2}). 

\vskip 0.2 cm 

We choose the path of integration $L(x)$ connecting 0 to $x$, defined by  
 $\arg(s)= {\Re \nu_2-\nu^* \over \Im \nu_2} \log|s| +b$, 
 where $b=\arg x - {\Re \nu_2-\nu^* \over \Im \nu_2} \log|x|$. Namely:
$$
\arg(s)=\arg(x) +{\Re \nu_2-\nu^* \over \Im \nu_2} \log{|s|\over |x|}
$$
If $x$ belongs to the domain 
(\ref{speriamobene2}) (or to  ${\cal D}(r;\nu_1,\nu_2)$ ) 
than  the path  does not leave the domain  when $s\to 0$, provided that 
$$  0 <\nu^* <2. 
$$
If $\Im \nu_2=0$ we take the path $\arg s = \arg x$, namely $
\nu^*=\nu_2$. 
The parameterization of the path is 
$$
   s=\rho ~e^{i\left\{ \arg x + {\Re \nu_2 - \nu^* \over \Im \nu_2} \log {\rho 
\over |x|}\right\}},~~~~0<\rho\leq|x|
$$
therefore
$$
|ds|=P(\nu_2,\nu^*)~d\rho,~~~~~P(\nu_2,\nu^*):= \sqrt{1+\left(
 {\Re \nu_2- \nu^* \over \Im \nu_2}\right)^2}
$$
We  observe that for any complex numbers $A$, $B$ we have
\be
\int_{L(x)} {1\over |s|} \left( |s|+ |A s^{2-\nu_2}|+ |B s^{\nu_2}| \right)^n 
~|ds| 
 \leq 
{P(\nu_2,\nu^*) \over n \min(\nu^*,2-\nu^*)} \left(|x|+ 
|A x^{2-\nu_2}|+ |B x^{\nu_2}| \right)^n
\label{CONDIZRIMS1}
\ee
This follows from the consideration that on $L(x)$ we have
$$
   |s^{\nu_2}|=|x^{\nu_2}| {|s|^{\nu^*} \over |x|^{\nu^*}}.
$$
Therefore
$$
\int_{L(x)} {1\over |s|} |s|^i|As^{2-\nu_2}|^j |Bs^{\nu_2}|^k  ~|ds| 
= 
  { |A x^{2-\nu_2}|^j~|Bx^{\nu_2}|^k \over|x|^{(2-\nu^*)j} ~|x|^{\nu^* k}} 
~P(\nu_2,\nu^*)\int_0^{|x|} d\rho~\rho^{i-1+(2-\nu^*)j+\nu^* k} 
$$
$$
= {P(\nu_2,\nu^*)\over i+j(2-\nu^*)+k\nu^*} 
|x|^i ~|A x^{2-\nu_2}|^j~|B x^{\nu_2}|^k\leq 
 {P(\nu_2,\nu^*)\over (i+j+k) \min(\nu^*,2-\nu^*)} 
|x|^i ~|A x^{2-\nu_2}|^j~|B x^{\nu_2}|^k
$$
from which (\ref{CONDIZRIMS1}) follows, provided that $0<\nu^*<2$. For 
$\Im \nu_2=0$ this brings again $0<\nu_2<2$.

\vskip 0.2 cm 
We observe that a solution of the integral equations is also a solution of the 
differential equations, by virtue of 
the analogous of Sub-Lemma 1 of section \ref{proof of theorem 1}: 

\vskip 0.2 cm 
\noindent
{\bf Sub-Lemma 2:} {\it Let $f(x)$ be a holomorphic function in the domain 
$|x|<\epsilon$,  $|A x^{2-\nu_2}|<\epsilon$, 
$|B x^{\nu_2}|<\epsilon$, such that 
$f(x)=O(|x|+|Ax^{2-\nu_2}|+|Bx^{\nu_2}|)$, $A,B\in {\bf C}$. Let $L(x)$ 
be the path of integration  define above for $0<\nu^*<2$ and 
$$
 F(x):= \int_{L(x)} {1\over s} f(s) ~ds
$$
Then, $F(x)$ is holomorphic on the domain and ${dF(x)\over dx}= {1\over x} f(x)
$
}

\vskip 0.2 cm
\noindent
{\it Proof:} We repeat exactly the argument of the proof of  Sub-Lemma 1, 
section \ref{proof of theorem 1}. We choose the point $x+\Delta x$ close to $x$ and we prove that $\int_{L(x)}- \int_{L(x+\Delta x)} = \int_x^{x+\Delta x} $, where the last integral is on a segment. Again, we reduce to the evaluation 
of the integral in the small portion of $L(x)$, $L(x+\Delta x)$ contained in 
the disc $U_R$ of radius $R<|x|$ and on the arc $\gamma(x_R,x^{\prime}_R)$ on 
the circle $|s|=R$.
Taking into account that $f(x)=O(|x|+|Ax^{2-\nu_2}|+|Bx^{\nu_2}|)$ and 
(\ref{CONDIZRIMS1}) we have 
$$ 
 \left| \int_{L(x_R)} {1\over s} f(s) ~ds\right|\leq \int_{L(x_R)} 
{1\over |s|} O(|s| +|As^{2-\nu_2}|+|Bs^{\nu_2}| )~ |ds|
$$
$$
\leq {P(\nu_2,\nu^*)  \over 
\min(\nu^*,2-\nu^*)} O(|x_R| + |Ax_R^{2-\nu_2}|+|Bx_R^{\nu_2}|)= 
{P(\nu_2,\nu^*)
  \over 
\min(\nu^*,2-\nu^*)}~O(R^{\min\{\nu^*,2-\nu^*\}})
$$
The last step follows from  $|x_R^{\nu_2}|= {|x^{\nu_2}| \over |x|^{\nu^*}}
 R^{\nu^*}$. So the integral vanishes for $R\to 0$. The same is proved for 
$\int_{L(x+\Delta x)}$. As for the integral on the arc we have 
$$
 |\arg x_R-\arg x^{\prime}_R|= \left| \arg x -\arg(x+\Delta x) +{\Re \nu_2 -
\nu^* \over \Im \nu_2} \log\left| 1 +{\Delta x\over x}\right|
\right|$$
or $
 |\arg x_R-\arg x^{\prime}_R|= \left| \arg x -\arg(x+\Delta x)\right|$ if $\Im \nu_2=0$. This is independent of $R$, 
therefore the length of the arc is $O(R)$ and 
$$
 \left|\int_{\gamma(x_R,x^{\prime}_R){1\over| s|} |f(s)|}~|ds| \right|= O(R^{
\min\{\nu^*,2-\nu^*\}}) \to 0 \hbox{ for } x\to 0
$$
\rightline{$\Box$}

\vskip 0.2 cm
Now we prove a fundamental lemma:

\vskip 0.2 cm 
\noindent
{\bf Lemma 6:} 
{\it 
For any complex $\nu_1$,  $\nu_2$  such that 
$$\nu_2 \not\in (-\infty,0]\cup[2,+\infty)$$ there exists a sufficiently small
$r<1$ such that 
the system of integral equations 
 has a solution $v(x)$ holomorphic in 
$$
  {\cal D}(r;\nu_1,\nu_2):= \left\{ x\in \tilde{\bf C}_0~ \hbox{ such that }
  |x|<r, \left|{e^{-i\pi \nu_1}\over 16^{2-\nu_2}} x^{2-\nu_2} \right|<r,
\left| 
{e^{i\pi \nu_1} \over 16^{\nu_2}} x^{\nu_2}\right|<r \right\}
$$
Moreover, there exists a constant $M(\nu_2)$ depending on $\nu_2$ such that 
 $v(x)\leq M(\nu_2) \left(|x|+\left|{e^{-i\pi \nu_1}\over 16^{2-\nu_2}} x^{2-\nu_2} \right|+\left| 
{e^{i\pi \nu_1} \over 16^{\nu_2}} x^{\nu_2}\right| \right)$ in  ${\cal
D}(r;\nu_1,\nu_2)$ . 
}

\vskip 0.2 cm

To prove Lemma 6 we need some sub-lemmas
\vskip 0.2 cm
\noindent
{\bf Sub-Lemma 3:} 
{\it Let $\Phi(x,y,z)$ and $\Psi(x,y,z,v,w)$ be two holomorphic 
functions of their arguments for $|x|,|y|,|z|,|v|,|w|< \epsilon$, satisfying
$$
   \Phi(0,0,0)=0,~~~~\Psi(0,0,0,v,w)=\Psi(x,y,z,0,0)=0
$$
Then, there exists a constant $c>0$ such that:
\be 
     \left| \Phi(x,y,z) \right|\leq c~\left(|x|+|y|+|z|\right)
\label{primacopia1}
\ee
\be
   \left| \Psi(x,y,z,v,w) \right|\leq c~\left(|x|+|y|+|z|\right)
\label{primacopia2}
\ee
\be 
\left| \Psi(x,y,z,v_2,w_2)- \Psi(x,y,z,v_1,w_1)\right|\leq c
\left(|x|+|y|+|z|\right)~\left( |v_2-v_1|+|w_2-w_1|\right)
\label{primacopia3}
\ee
 for $|x|,|y|,|z|,|v|,|w|< \epsilon$.
}

\vskip 0.2 cm
\noindent
{\it Proof:} Let's proof (\ref{primacopia2}). 
$$
 \Psi(x,y,z,v,w)=\int_0^1 ~{d\over d\lambda}
\Psi(\lambda x, \lambda y, \lambda z,v,w)~d\lambda
$$
$$
  x~\int_0^1 {\partial \Psi \over \partial x}(\lambda x, \lambda y, 
\lambda z,v,w) ~d\lambda~+~y\int_0^1 {\partial \Psi \over \partial y}
(\lambda x, \lambda y, 
\lambda z,v,w) ~d\lambda~+~z\int_0^1 {\partial \Psi \over \partial z}
(\lambda x, \lambda y, 
\lambda z,v,w) ~d\lambda
$$
Moreover, for $\delta$ small:
$$ 
{\partial \Psi \over \partial x}(\lambda x, \lambda y, 
\lambda z,v,w)= \int_{|\zeta-\lambda x|=\delta}~
{\Psi(\zeta,\lambda y,\lambda z,
v,w) \over (\zeta-\lambda x)^2}~{d\zeta \over 2 \pi i} 
$$
which implies that ${\partial \Psi \over \partial x}$ is holomorphic and 
bounded when its arguments are less than $\epsilon$. The same holds true for 
 ${\partial \Psi \over \partial y}$ and  ${\partial \Psi \over \partial z}$. This proves (\ref{primacopia2}), $c$ being a constant which bounds 
$\left|{\partial \Psi \over \partial x}\right|$, $\left|
{\partial \Psi \over \partial y}\right|$ $\left|{\partial \Psi 
\over \partial z}\right|$. The inequality (\ref{primacopia1}) is proved in 
the same way.  
We turn to (\ref{primacopia3}). First we prove that for 
$|x|,|y|,|z|,|v_1|,|w_1|,|v_2|,|w_2| < \epsilon$ there exist two 
holomorphic and bounded functions $\psi_1(x,y,z,v_1,w_1,v_2,w_2)$,  
$\psi_2(x,y,z,v_1,w_1,v_2,w_2)$ such that
$$
    \Psi(x,y,z,v_2,w_2)-\Psi(x,y,z,v_1,w_1)
$$
\be
= (v_2-v_1)~ \psi_1(x,y,z,v_1,w_1,v_2,w_2)+ (w_1-w_2)~
\psi_2(x,y,z,v_1,w_1,v_2,w_2)
\label{tmpRIMS1}
\ee 
In order  to prove this, we write 
$$
 \Psi(x,y,z,v_2,w_2)-\Psi(x,y,z,v_1,w_1)=
$$
$$
=
\int_0^1 {d\over d\lambda} \Psi(x,y,z,\lambda v_2 +(1-\lambda)v_1,
\lambda w_2 +(1-\lambda)w_1)~ d\lambda
$$
$$
=(v_2-v_1) \int_0^1 {\partial \Psi \over \partial v} 
 (x,y,z,\lambda v_2 +(1-\lambda)v_1,
\lambda w_2 +(1-\lambda)w_1)~d\lambda~+
$$
$$
+~(w_2-w_1) \int_0^1
 {\partial \Psi \over \partial w} 
 (x,y,z,\lambda v_2 +(1-\lambda)v_1,
\lambda w_2 +(1-\lambda)w_1)~d\lambda
$$
$$
=: (v_2-v_1)~ \psi_1(x,y,z,v_1,w_1,v_2,w_2)+ (w_2-w_1) ~
\psi_2(x,y,z,v_1,w_1,v_2,w_2)
$$
Moreover, for small $\delta$,  
$$
{\partial \Psi \over \partial v}(x,y,z,v,w)= \int_{|\zeta-v|=\delta} 
{\Psi(x,y,z,\zeta,w) \over (\zeta-v)^2} ~{dz\over 2\pi i} 
$$
which implies that $\psi_1$ is holomorphic and bounded for its arguments less 
 than $\epsilon$. We also obtain 
${\partial \Psi \over \partial v}(0,0,0,v,w)=0$, then  
$\psi_1(0,0,0,v_1,w_1,v_2,w_2)=0$. The proof for $\psi_2$ is analogous.
 We use (\ref{tmpRIMS1}) to complete the proof of (\ref{primacopia3}). 
Actually, we observe that 
$$ 
  \psi_i(x,y,z,v_1,w_1,v_2,w_2)= \int_0^1 {d \over d\lambda} 
\psi_i(\lambda x,\lambda y, \lambda z, v_1,w_1,v_2,w_2) ~d\lambda 
$$
$$
 = x\int_0^1 {\partial \psi_i \over \partial x} d\lambda + 
 y\int_0^1 {\partial \psi_i \over \partial y} d\lambda + 
z\int_0^1 {\partial \psi_i \over \partial z} d\lambda 
$$
and we conclude as in the proof of (\ref{primacopia2}).

\rightline{$\Box$}

\vskip 0.2 cm

We solve the system of integral equations 
by successive approximations.  We can choose any path $L(x)$ 
such that $0<\nu^*<2$. Here we choose $\nu^*=1$.  For convenience, we put  
$$
A:= {e^{-i\pi \nu_1}\over 16^{2-\nu_2}},~~~~
B:= {e^{i\pi \nu_1} \over 16^{\nu_2}}
$$
Therefore, for any $n\geq 1$ the successive approximations are: 
$$
    v_0=w_0=0
$$
\vskip 0.15 cm 
\be
w_n(x)=\int_{L(x)} {1\over t} \left\{
 \Phi (s,A s^{2-\nu_2}, 
B s^{\nu_2} ) + \Psi(s,A s^{2-\nu_2}, B s^{\nu_2},v_{n-1}(s),w_{n-1}(s))
\right\}~ds
\label{RIMSRIMS1}
\ee
\be
v_n(x)=\int_{L(x)}{1\over s}~w_{n}(s)~ ds
\label{RIMSRIMS2}
\ee

\vskip 0.2 cm
\noindent
{\bf Sub-Lemma 4:} {\it There exists a sufficiently small $\epsilon^{\prime}
<\epsilon$ such that for 
any $n\geq 0$  the functions 
$v_n(x)$ and $w_n(x)$ are holomorphic   in the domain 
$$
   {\cal D}(\epsilon^{\prime};\nu_1,\nu_2):= \left\{ x\in \tilde{\bf C}_0~ 
\hbox{ such that }
  |x|<\epsilon^{\prime}, \left|A x^{2-\nu_2} \right|<\epsilon^{\prime},
\left| B  x^{\nu_2}\right|<\epsilon^{\prime} \right\}
$$
They are also correctly bounded, namely $|v_n(x)|<\epsilon$, 
$|w_n(x)|<\epsilon$ for any
 $n$. They  satisfy 
\be 
 |v_n-v_{n-1}| \leq { (2c)^n P(\nu_2)^{2n} \over n!} \left(
|x|+|Ax^{2-\nu_2}|+|Bx^{\nu_2}|
\right)^n
\label{InEqAlItY1}
\ee
\be
 |w_n-w_{n-1}| \leq { (2c)^n P(\nu_2)^{2n} \over n!} \left(
|x|+|Ax^{2-\nu_2}|+|Bx^{\nu_2}|
\right)^n
\label{InEqAlItY}
\ee
where $P(\nu_2) :=P(\nu_2,\nu^*=1)$ and $c$ is the constant appearing in 
Sub-Lemma 3. Moreover 
$$
    x {dv_n \over dx}=w_n
$$
}

\vskip 0.2 cm 
\noindent
{\bf Proof:} We proceed by induction. 
$$ 
  w_1 = \int_{L(x)} {1\over s} 
 \Phi (s,A s^{2-\nu_2}, 
B s^{\nu_2} ) ~ds,
~~~~v_1=\int_{L(x)}{1\over s} w_1(s) ds 
$$
 It follows from Sub-Lemma 2 and (\ref{primacopia1}) that $w_1(x)$ is 
holomorphic for $|x|,|Ax^{2-\nu_2}|, |Bx^{\nu_2}|<\epsilon$. 
From (\ref{CONDIZRIMS1})  and (\ref{primacopia1}) we have
$$
   |w_1(x)|\leq \int {1\over |s| } |\Phi(s,As^{2-\nu_2},Bs^{\nu_2})|~|ds| 
$$
$$
\leq c P(\nu_2) (|x|+|Ax^{2-\nu_2}|+|Bx^{\nu_2}|) \leq 3cP(\nu_2) \epsilon^{\prime}<\epsilon
$$
 on 
 $ {\cal D}(\epsilon^{\prime};\nu_1,\nu_2)$, provided 
that $\epsilon^{\prime}$ is small enough.   
By Sub-Lemma 2, also $v_1(x)$ is holomorphic   
for $|x|,|Ax^{2-\nu_2}|, |Bx^{\nu_2}|<\epsilon$ and
$$
 x {d v_1 \over dx} =w_1
$$
By (\ref{CONDIZRIMS1}) we also have 
$$
   |v_1(x) | \leq cP(\nu_2)^2  (|x|+|Ax^{2-\nu_2}|+|Bx^{\nu_2}|)
\leq 3cP(\nu_2)^2 \epsilon^{\prime}<\epsilon
$$ 
on $ {\cal D}(\epsilon^{\prime};\nu_1,\nu_2)$. 
Note that $P(\nu_2)\geq 1$, so (\ref{InEqAlItY}) (\ref{InEqAlItY1}) are
 true for $n=1$. Now we 
 suppose that the statement of the sub-lemma 
is true for $n$ and we prove it for $n+1$. Consider: 
$$
|w_{n+1}(x) - w_n(x)| = \left| 
 \int_{L(x)}{1\over s}~ \left[
            \Psi(s,As^{2-\nu_2},Bs^{\nu_2},v_n,w_n)- 
\Psi(s,As^{2-\nu_2},Bs^{\nu_2},v_{n-1},w_{n-1}) 
\right] 
~ds
\right|
$$
By (\ref{primacopia3}) the above is  
$$
\leq 
     c~\int_{L(x)} {1\over |s|} (|s|+|As^{2-\nu_2}|+|Bs^{\nu_2}|)~
(|v_n-v_{n-1}| + |w_n-w_{n-1}|) ~|ds|
$$
By induction this is 
$$
 \leq 2c~ {(2c)^n P(\nu_2)^{2n} \over n!} ~ \int_{L(x)} {1\over |s|} 
(|s|+|As^{2-\nu_2}|+|Bs^{\nu_2}|)^{n+1} ~|ds|
$$
$$
\leq 
      2c~  {(2c)^n P(\nu_2)^{2n} \over n!} ~ {P(\nu_2)\over n+1} 
~(|x|+|Ax^{2-\nu_2}|+|Bx^{\nu_2}|)^{n+1}
$$
$$
\leq 
{(2c)^{n+1} P(\nu_2)^{2(n+1)} \over (n+1)!}
(|x|+|Ax^{2-\nu_2}|+|Bx^{\nu_2}|)^{n+1}
$$
This proves (\ref{InEqAlItY}). Now we estimate 
$$
   |v_{n+1}(x)-v_n(x)|\leq \int_{L(x)} |w_{n+1}(s)-w_n(s)|~|ds|
$$
$$
\leq {(2c)^{n+1} P(\nu_2)^{2n+1} \over (n+1)!}~\int_{L(x)} 
{1\over |s| } 
(|s|+|As^{2-\nu_2}|+|Bs^{\nu_2}|)^{n+1}~|ds|
$$
$$
\leq 
{(2c)^{n+1} P(\nu_2)^{2(n+1)} \over (n+1)~(n+1)!}~
(|x|+|Ax^{2-\nu_2}|+|Bx^{\nu_2}|)^{n+1}
$$
$$
\leq 
{(2c)^{n+1} P(\nu_2)^{2(n+1)} \over ~(n+1)!}~
(|x|+|Ax^{2-\nu_2}|+|Bx^{\nu_2}|)^{n+1}
$$
This proves (\ref{InEqAlItY1}). From Sub-Lemma 2 we also conclude that 
$w_n$ and $v_n$ are holomorphic in ${\cal D}(\epsilon^{\prime},\nu_1,\nu_2)$ 
and 
$$
  x{dv_n \over dx}=w_n
$$
Finally we see that  
$$ 
|v_n(x)|\leq \sum_{k=1}^n |v_k(x)-v_{k-1}(x)|\leq \exp\{ 2cP^2(\nu_2)
 (|x|+|Ax^{2-\nu_2}|+|Bx^{\nu_2}|)\}-1 \leq \exp\{ 6cP^2(\nu_2)
                                  \epsilon^{\prime}\}-1
$$
and the same for $|w_n(x)|$. Therefore, if $\epsilon^{\prime}$ is small 
enough  we have $|v_n(x)|<\epsilon$, $|w_n(x)|<\epsilon$ on 
 ${\cal D}(\epsilon^{\prime},\nu_1,\nu_2)$. 

\rightline{$\Box$}

\vskip 0.2 cm

Let's define 
$$
 v(x):=\lim_{n\to \infty} v_n(x),~~~~~w(x):=\lim_{n\to \infty} w_n(x) 
$$
if they exist. 
We can also rewrite
$$
 v(x) = \lim_{n\to \infty} v_n(x)= \sum_{n=1}^{\infty} (v_n(x)-v_{n-1}(x)). 
$$
 We see that the series converges uniformly in 
 ${\cal D}(\epsilon^{\prime},\nu_1,\nu_2)$ because 
$$
|\sum_{n=1}^{\infty} (v_n(x)-v_{n-1}(x)) | 
$$
$$
\leq 
\sum_{n=1}^{\infty} 
{(2c)^{n} P(\nu_2)^{2n} \over ~n!}~
(|x|+|Ax^{2-\nu_2}|+|Bx^{\nu_2}|)^{n}
$$
$$
 = \exp\{ 2 cP^2(\nu_2)
 (|x|+|Ax^{2-\nu_2}|+|Bx^{\nu_2}|) \} -1
$$
The same holds for $w_n(x)$. Therefore, $v(x)$ and $w(x)$ define holomorphic functions in  
 ${\cal D}(\epsilon^{\prime},\nu_1,\nu_2)$. From Sub-Lemma 4 we also have 
$$
    x {d v(x) \over dx } = w(x)
$$
in ${\cal D}(\epsilon^{\prime},\nu_1,\nu_2)$. 

\vskip 0.2 cm
We show that $v(x),w(x)$ solve the initial integral equations. The l.h.s. of 
(\ref{RIMSRIMS1}) converges to $w(x)$ for $n\to \infty$. Let's prove that the 
r.h.s.  also converges to 
$$\int_{L(x)} {1\over s} \left\{
 \Phi (s,A s^{2-\nu_2}, 
B s^{\nu_2} ) + \Psi(s,A s^{2-\nu_2}, B s^{\nu_2},v(s),w(s))
\right\}~ds. $$ 
 We have to evaluate the following difference: 
$$
   \left|
\int_{L(x)} {1\over s} \Psi(s,A s^{2-\nu_2}, B s^{\nu_2},v(s),w(s))
~ds ~-~ 
\int_{L(x)} {1\over s} \Psi(s,A s^{2-\nu_2}, B s^{\nu_2},v_n(s),w_n(s))
~ds
\right|
$$
By (\ref{primacopia3}) the above is 
\be
c~\leq \int_{L(x)} {1\over |s|} (|s|+|As^{2-\nu_2}|+|Bs^{\nu_2}|)~(
|v-v_n|+|w-w_n|) ~|ds|
\label{palla}
\ee
Now we observe that
$$
|v(x)-v_n(x)| \leq \sum_{k=n+1}^{\infty} |v_k-v_{k-1}| 
$$
$$
  = \sum_{k=n+1}^{\infty} {(2c)^k P(\nu_2)^{2k} \over k!} (|x|+|Ax^{2-\nu_2} |
+ |Bx^{\nu_2}| )^k
$$
$$
\leq 
     (|x|+|Ax^{2-\nu_2} |
+ |Bx^{\nu_2}| )^{n+1} 
  ~\sum_{k=0}^{\infty} {(2c)^{k+n+1} P(\nu_2)^{2(k+n+1)} \over (k+n+1)!} 
(|x|+|Ax^{2-\nu_2} |
+ |Bx^{\nu_2}| )^k
$$
The series converges. Its sum is less than some constant $S(\nu_2)$ 
 independent of $n$. We obtain 
$$
  |v(x)-v_n(x)| \leq S(\nu_2)   (|x|+|Ax^{2-\nu_2} |
+ |Bx^{\nu_2}| )^{n+1}. $$
The same holds for $|w-w_n|$. Thus, (\ref{palla}) is 
$$
\leq 
  2c~S(\nu_2) 
~\int_{L(x)} {1\over |s| } (|s|+|As^{2-\nu_2}|+|Bs^{\nu_2}|)^{n+2}
~|ds| 
\leq 
  {2cS(\nu_2) P(\nu_2)\over n+2}~  (|x|+|Ax^{2-\nu_2}|+|Bx^{\nu_2}|)^{n+2}
$$
Namely: 
$$
   \left|
\int_{L(x)} {1\over s} \Psi(s,A s^{2-\nu_2}, B s^{\nu_2},v(s),w(s))
~ds ~-~ 
\int_{L(x)} {1\over s} \Psi(s,A s^{2-\nu_2}, B s^{\nu_2},v_n(s),w_n(s))
~ds
\right|
$$
$$
\leq 
{2cS(\nu_2) P(\nu_2)\over n+2}~(3\epsilon^{\prime})^{n+2}
$$
In a similar way, the r.h.s. of (\ref{RIMSRIMS2}) is 
$$
    \left|\int {1\over s} (w(s)-w_n(s))~ds\right|\leq {S(\nu_2) P(\nu_2) \over n+1} 
(3\epsilon^{\prime})^{n+1}
$$
Therefore, the r.h. sides  of  (\ref{RIMSRIMS1}) (\ref{RIMSRIMS2}) 
converge on the domain ${\cal D}(r,\nu_1,\nu_2)$ for 
$r< \min\{\epsilon^{\prime}, 1/3\}$.  
We finally observe that $|v(x)|$ and $|w(x)|$ are bounded on ${\cal D}(r)$. 
For example  
$$
|v(x)| \leq  (|x|+|Ax^{2-\nu_2} |
+ |Bx^{\nu_2}| ) 
  ~\sum_{k=0}^{\infty} {(2c)^{k+1} P(\nu_2)^{2(k+1)} \over (k+1)!} 
(|x|+|Ax^{2-\nu_2} |
+ |Bx^{\nu_2}| )^k
$$
$$
=: M(\nu_2)  (|x|+|Ax^{2-\nu_2} |
+ |Bx^{\nu_2}| )
$$
where  the sum of the series is less than a constant $M(\nu_2)$. 
We have proved Lemma 6.

\rightline{$\Box$}

\vskip 0.3 cm 

We note that the proof of Lemma 6 only makes use of the properties of 
$\Phi$ and $\Psi$, regardless of how these functions have been constructed. 
The structure of the integral equations implies that $v(x)$ is bounded 
(namely $|v(x)|=O(r)$). 
Now, we come back to our case, where $\Phi$ and $\Psi$ have been constructed from the Fourier expansion of elliptic functions. We need to check  if 
(\ref{speriamobene2}) and ${\cal D}(r,\nu_1,\nu_2)$ have non-empty 
intersection. This is true, indeed 
  ${\cal D}(r)$   is contained in (\ref{speriamobene2}), 
because in (\ref{speriamobene2}) the term $\Im {\pi v \over 2 \omega_1}$ is 
$O(r)$, while in ${\cal D}(r,\nu_1,\nu_2)$ 
the term $\ln r$ appear, and $r$ is small.  

\vskip 0.2 cm 
 To conclude the proof of Theorem 3 we have to work out the explicit series 
(\ref{vdix}). In order to do this we observe that $w_1$ and $v_1$ 
are series of the type 
\be 
  \sum_{p,q,r\geq 0} c_{pqr}(\nu_2) ~x^p~(Ax^{2-\nu_2})^q~(Bx^{\nu_2})^r
\label{formMrims}
\ee
where $c_{pqr}(\nu_2)$ is rational in $\nu_2$. This follows from 
$$
   w_1(x)= \int_{L(x)} \Phi(s,As^{2-\nu_2},Bs^{\nu_2})~ds
$$
and from the fact that $\Phi(x,Ax^{2-\nu_2},Bx^{\nu_2})$ itself 
 is a series 
 (\ref{formMrims}) by construction, with  coefficients  
 $c_{pqr}(\nu_2)$ which are  rational functions of  $\nu_2$.  
The same holds true for  $\Psi$.  We conclude that 
  $w_n(x)$ and $v_n(x)$ have the form (\ref{formMrims}) for any $n$. This 
implies that the limit $v(x)$ is also a series of type (\ref{formMrims}). We 
can reorder such a series to obtain (\ref{vdix}). Consider the term 
$$
   c_{pqr}(\nu_2) ~  x^p~(Ax^{2-\nu_2})^q~(Bx^{\nu_2})^r,$$
and recall that by definition $B={1\over  16^2 ~A}$. We absorb $16^{-2r}$ into 
$ c_{pqr}(\nu_2)$ and we study the factor 
$$
    A^{q-r} x^{p+(2-\nu_2)q+\nu_2 r} =  A^{q-r} x^{p+2q+(r-q)\nu_2}
$$
We have three cases: 

1) $r=q$, then we have $x^{p+2q}=: x^n$, $n=p+2q$. 

2) $r>q$, then we have $x^{p+2q}~\left[{1\over A} x^{\nu_2}\right]^{r-q}=:
x^n~\left[{1\over A} x^{\nu_2}\right]^m$, $n=p+2q$, $m=q-r$.

3) $r<q$,  then we have $A^{q-r} ~x^{p+2r} ~\left[ A x^{2-\nu_2}\right]^{q-r}
 =: x^n~\left[ A x^{2-\nu_2}\right]^m$, $n= {p+2r}$, $m=q-r$.. 

\noindent
This brings a series of the type  (\ref{formMrims}) to the form 
(\ref{vdix}). The proof of Theorem 3 is complete.

\vskip 0.3 cm 
A system of integral equations similar to the one we considered here was 
 first  
studied by S.Shimomura in \cite{Sh}
 and \cite{IKSY}.


\vskip 1 cm
{\noindent}
{\it NOTES:}
\vskip 0.2 cm

\noindent
1. (Section \ref{Introduction})

There are only some exceptions to the one-to-one correspondence above, which
are already treated in \cite{M}. In order to rule them out we require 
that at most one of the entries $x_i$ of the triple may be zero and that
$(x_0,x_1,x_{\infty})\not \in \{(2,2,2)$ $(-2,-2,2)$, $(2,-2,-2)$, $(-2,2,-2)
\}$. See also Note 2 below. 
 Two triples which differ by the change of two signs identify the same transcendent. They 
are called {\it equivalent triples}. The one to one correspondence is between transcendents and classes of equivalence. 

\vskip 0.2 cm 
\noindent
2. (Section \ref{Monodromy Data and Review of Previous Results})

 The proof  is 
 done in the following way: consider
two solutions $C$ and $\tilde{C}$ of the equations (\ref{conn1}),
(\ref{conn2}). Then 
$$ 
    (C_i \tilde{C}_i^{-1})^{-1} e^{2\pi i J} (C_i \tilde{C}_i^{-1})=
    e^{2\pi i J}
$$
$$ 
 (C_{\infty} \tilde{C}_{\infty}^{-1})^{-1} e^{-2\pi i
 \hat{\mu}}e^{2\pi i R} 
 (C_{\infty} 
\tilde{C}_{\infty}^{-1})=
    e^{-2\pi i \hat{\mu}} e^{2\pi i R}
$$
We find  
$$
C_i \tilde{C}_i^{-1}= \pmatrix{a & b \cr
                               0 & a  \cr}  ~~~~~~a,b\in{\bf C},~~~a\neq 0
$$
Note that this matrix commutes with $J$, then 
$$ 
  (z-u_i)^J C_i = (z-u_i)^J \pmatrix{a & b \cr
                               0 & a  \cr}\tilde{C}_i=\pmatrix{a & b \cr
                               0 & a  \cr}(z-u_i)^J \tilde{C_i}
$$
We also find 
\be
C_{\infty} 
\tilde{C}_{\infty}^{-1}
=
    \left\{ \matrix{ 
                     i)~~~~~~~\hbox{diag}(\alpha, \beta),~~~\alpha\beta\neq 0;
                 ~~~~~~~~~~~~~~~~~~~~~~~~~~~~~~~~~~~~ 
                    ~\hbox{if } 2\mu \not \in {\bf Z} \cr
                   ii) ~ \pmatrix{\alpha & \beta \cr
                                  0  &  \alpha \cr} ~~(\mu>0),~~~ 
 \pmatrix{\alpha & 0 \cr
                                  \beta  &  \alpha \cr} ~~(\mu<0),
~~~\alpha\neq 0, ~~~\hbox{if } 2\mu \in {\bf Z}, R\neq 0 \cr
iii)~ \hbox{ Any invertible matrix} ~~~~~~~~~~~~~~~~~~~~~~~~~~~~~~~~~~~~
~\hbox{ if }  2\mu \in {\bf Z}, R= 0 \cr
}\right.
\label{alga}
\ee
Then 
$$ i)~ z^{-\hat{\mu}}C_{\infty}= z^{-\hat{\mu}}\hbox{diag}(\alpha,\beta)
\tilde{C}_{\infty} = \hbox{diag}(\alpha,\beta)z^{\hat{\mu}}\tilde{C}_{\infty}
$$
$$
ii)~z^{-\hat{\mu}}z^{-R} C_{\infty} = ~...~= \left[\alpha I+ {1\over
z^{|2\mu|}} Q\right] 
     ~z^{-\hat{\mu}}z^{-R} \tilde{C}_{\infty}
$$
where $Q=\pmatrix{0 & \beta \cr 0 & 0\cr }$, or $Q= \pmatrix{0 &  0\cr
\beta  & 0 \cr}$. 
$$
iii)~   z^{-\hat{\mu}}C_{\infty}=~...~= \left[
              {Q_1\over z^{|2\mu|}} + Q_0 + Q_{-1}z^{|2\mu|}  
\right]
             z^{-\hat{\mu}} \tilde{C}_{\infty}
$$
where $Q_{0}=$ diag$(\alpha, \beta)$, $Q_{\pm 1}$ are respectively
upper and lower triangular (or lower and upper triangular, depending
on the sign of $\mu$), and $C_{\infty}
\tilde{C}_{\infty}^{-1}=Q_1+Q_0+Q_{-1}$  

This implies that the  two solutions $Y(z;u)$, $\tilde{Y}(z;u)$ of the
form  (\ref{solrh}) with $C_{\nu}$ and $\tilde{C}_\nu$ respectively 
($\nu=1,2,3,\infty$), are such
that $Y(z;u)~\tilde{Y}(z;u)^{-1}$ is holomorphic at each $u_i$, while at
$z=\infty$ it is 
$$
   Y(z;u)~\tilde{Y}(z;u)^{-1} \to 
\left\{ \matrix{
                 i)~ G_{\infty} \hbox{diag}(\alpha,\beta)
G_{\infty}^{-1} \cr
ii)~ \alpha I \cr
iii)~ \hbox{ divergent }
    }
\right.
$$
Thus the two fuchsian systems are conjugated only in the   cases
$i)$ and $ii)$, because in those cases $Y\tilde{Y}^{-1}$ is
holomorphic everywhere on ${\bf P}^1$, and then it is a constant.  
In other words {\it the
R.H. has a unique solution, up to conjugation, for $2 \mu \not \in {\bf
Z}$ or for $2 \mu \in {\bf Z}$ and $R\neq 0$}.

\vskip 0.2 cm
\noindent
3. (Section \ref{Monodromy Data and Review of Previous Results})

Note that if $G_{\infty}=I$, then $\sum_{i=1}^3
A_i$ is  already  diagonal. Therefore, there is no loss of generality if,  for $2\mu \not \in {\bf Z}$, we solve the Riemann-Hilbert problem for given $M_1$, $M_2$, $M_3$ with 
the choice of normalization $Y(z;u)z^{\hat{\mu}}
 \to I$ if $z\to \infty$. This  determines uniquely $A_1$, $A_2$,
$A_3$ up to diagonal conjugation. Note that  for 
 any diagonal invertible matrix $D$, the sum of  $D^{-1} A_i D$  is still
 diagonal.

On the other hand, if $2 \mu\in {\bf Z}$ and $R\neq 0$, then $Y(z;u)~
\tilde{Y}(z;u)^{-1}=\alpha$, where $\alpha$ appears in (\ref{alga}), 
case ii). Therefore the two fuchsian systems obtained from $A(z;u):=dY/dz~Y^{-1}$ 
and  $\tilde{A}(z;u)=d\tilde{Y}/dz~\tilde{Y}^{-1}$ coincide. 

 In both cases, $A_{12}(z,u)$ changes at most for the multiplication by a constant, therefore the same $y(x)$ is defined by  $A_{12}(z,u)=0$ and $\tilde{A}_{12}(z,u)=0$

\vskip 0.2 cm 
\noindent
4. (Section \ref{Monodromy Data and Review of Previous Results})

$R=0$ only 
in  the case 2) of
commuting monodromy matrices and $\mu$ integer.


\baselineskip 12pt
\vskip 1 cm
\noindent
{\bf Acknowledgments (spring 2001)} I am grateful to B.Dubrovin 
 for many  discussions  and advice. 
  I would
 like to thank  A.Bolibruch, 
  A.Its,  M. Jimbo, 
  M.Mazzocco, S.Shimomura  for fruitful discussions.  
The paper was completed while the  author was
 supported by a fellowship of  the Japan Society for the 
Promotion of Science (JSPS) in RIMS, Kyoto University.

\end{document}